 \documentclass[11pt,letterpaper]{article}

\usepackage[normalem]{ulem}
\usepackage{listings}
\usepackage{a4wide}
\usepackage{amsmath} 
\usepackage{amsthm}

\usepackage{titlesec}
\titleformat{\subsection}[hang]{\bfseries \filcenter}{}{1em}{}
\titleformat{\subsubsection}[hang]{\bfseries \filcenter}{}{1em}{}

\usepackage{amsfonts}
\usepackage{amssymb,amscd}
\usepackage{multirow}
\usepackage[mathcal]{eucal}
\usepackage[spanish,USenglish]{babel} 
\usepackage{graphics}  
\usepackage{graphicx} 
\usepackage{hyperref}
\usepackage{fancyhdr}
\usepackage{titlesec}
\usepackage[latin1]{inputenc}
\usepackage{amsmath}
\usepackage{amsfonts}
\usepackage{amssymb}
\usepackage{amsthm}
\usepackage{newlfont}
\usepackage{amsbsy}
\usepackage{bm}
\usepackage[all]{xy}
\usepackage{textcomp}
\usepackage{color}
\usepackage{epsfig}
\usepackage{dsfont}
\usepackage{latexsym}
\usepackage{array}
\usepackage{longtable}
\usepackage{enumerate}
\usepackage{multicol}
\usepackage{mathrsfs}
\usepackage{cancel}
\usepackage{wasysym}

\theoremstyle{plain}
\newtheorem{teor}{Theorem}
\newtheorem{lemma}{Lemma}
\newtheorem{coro}{Corollary}
\newtheorem{prop}{Proposition}

\theoremstyle{definition}

\newtheorem*{ejem*}{Examples}
\newtheorem{defin}{Definition}
\newtheorem*{demos}{Proof}
\newtheorem{remark}{Remark}
\theoremstyle{remark}

\newcommand{\complex}{\mathbf{\mathbb{C}}}
\newcommand{\reales}{\mathbf{\mathbb{R}}}

\newcommand{\natu}{\mathbf{\mathbb{N}}}
\newcommand{\llave}[1]{\left\{ #1\right\}}
\newcommand{\norma}[1]{\left \| #1 \right \|}
\newcommand{\vabs}[1]{\left| #1\right|}

\newcommand{\corch}[1]{\left[ #1\right]}
\newcommand{\paren}[1]{\left( #1\right)}
\newcommand{\QED}{\hfill \ensuremath{\Box}}

\usepackage[breakable]{tcolorbox}
\usepackage{parskip} 

\usepackage{iftex}
\ifPDFTeX
\usepackage[T1]{fontenc}
\usepackage{mathpazo}
\else
\usepackage{fontspec}
\fi

\usepackage{graphicx}

\usepackage{caption}
\DeclareCaptionFormat{nocaption}{}
\usepackage[Export]{adjustbox} 
\adjustboxset{max size={0.9\linewidth}{0.9\paperheight}}
\usepackage{float}
\usepackage{xcolor} 
\usepackage{enumerate} 
\usepackage{geometry} 
\usepackage{amsmath} 
\usepackage{amssymb} 
\usepackage{textcomp} 
\AtBeginDocument{%
}
\usepackage{upquote} 
\usepackage{eurosym} 
\usepackage[mathletters]{ucs} 
\usepackage{fancyvrb} 
\usepackage{grffile} 
\makeatletter 
\def\Gread@@xetex#1{%
	\IfFileExists{"\Gin@base".bb}%
	{\Gread@eps{\Gin@base.bb}}%
	{\Gread@@xetex@aux#1}%
}
\makeatother

\usepackage{hyperref}
\usepackage{titling}
\usepackage{longtable} 
\usepackage{booktabs}  
\usepackage[inline]{enumitem} 
\usepackage[normalem]{ulem} 
\usepackage{mathrsfs}

\definecolor{urlcolor}{rgb}{0,.145,.698}
\definecolor{linkcolor}{rgb}{.71,0.21,0.01}
\definecolor{citecolor}{rgb}{.12,.54,.11}

\definecolor{ansi-black}{HTML}{3E424D}
\definecolor{ansi-black-intense}{HTML}{282C36}
\definecolor{ansi-red}{HTML}{E75C58}
\definecolor{ansi-red-intense}{HTML}{B22B31}
\definecolor{ansi-green}{HTML}{00A250}
\definecolor{ansi-green-intense}{HTML}{007427}
\definecolor{ansi-yellow}{HTML}{DDB62B}
\definecolor{ansi-yellow-intense}{HTML}{B27D12}
\definecolor{ansi-blue}{HTML}{208FFB}
\definecolor{ansi-blue-intense}{HTML}{0065CA}
\definecolor{ansi-magenta}{HTML}{D160C4}
\definecolor{ansi-magenta-intense}{HTML}{A03196}
\definecolor{ansi-cyan}{HTML}{60C6C8}
\definecolor{ansi-cyan-intense}{HTML}{258F8F}
\definecolor{ansi-white}{HTML}{C5C1B4}
\definecolor{ansi-white-intense}{HTML}{A1A6B2}
\definecolor{ansi-default-inverse-fg}{HTML}{FFFFFF}
\definecolor{ansi-default-inverse-bg}{HTML}{000000}

\let\Oldtex\TeX
\let\Oldlatex\LaTeX
\renewcommand{\TeX}{\textrm{\Oldtex}}
\renewcommand{\LaTeX}{\textrm{\Oldlatex}}
\title{First Example}

\makeatletter
\def\PY@reset{\let\PY@it=\relax \let\PY@bf=\relax%
	\let\PY@ul=\relax \let\PY@tc=\relax%
	\let\PY@bc=\relax \let\PY@ff=\relax}
\def\PY@tok#1{\csname PY@tok@#1\endcsname}
\def\PY@toks#1+{\ifx\relax#1\empty\else%
	\PY@tok{#1}\expandafter\PY@toks\fi}
\def\PY@do#1{\PY@bc{\PY@tc{\PY@ul{%
				\PY@it{\PY@bf{\PY@ff{#1}}}}}}}
\def\PY#1#2{\PY@reset\PY@toks#1+\relax+\PY@do{#2}}

\expandafter\def\csname PY@tok@w\endcsname{\def\PY@tc##1{\textcolor[rgb]{0.73,0.73,0.73}{##1}}}
\expandafter\def\csname PY@tok@c\endcsname{\let\PY@it=\textit\def\PY@tc##1{\textcolor[rgb]{0.25,0.50,0.50}{##1}}}
\expandafter\def\csname PY@tok@cp\endcsname{\def\PY@tc##1{\textcolor[rgb]{0.74,0.48,0.00}{##1}}}
\expandafter\def\csname PY@tok@k\endcsname{\let\PY@bf=\textbf\def\PY@tc##1{\textcolor[rgb]{0.00,0.50,0.00}{##1}}}
\expandafter\def\csname PY@tok@kp\endcsname{\def\PY@tc##1{\textcolor[rgb]{0.00,0.50,0.00}{##1}}}
\expandafter\def\csname PY@tok@kt\endcsname{\def\PY@tc##1{\textcolor[rgb]{0.69,0.00,0.25}{##1}}}
\expandafter\def\csname PY@tok@o\endcsname{\def\PY@tc##1{\textcolor[rgb]{0.40,0.40,0.40}{##1}}}
\expandafter\def\csname PY@tok@ow\endcsname{\let\PY@bf=\textbf\def\PY@tc##1{\textcolor[rgb]{0.67,0.13,1.00}{##1}}}
\expandafter\def\csname PY@tok@nb\endcsname{\def\PY@tc##1{\textcolor[rgb]{0.00,0.50,0.00}{##1}}}
\expandafter\def\csname PY@tok@nf\endcsname{\def\PY@tc##1{\textcolor[rgb]{0.00,0.00,1.00}{##1}}}
\expandafter\def\csname PY@tok@nc\endcsname{\let\PY@bf=\textbf\def\PY@tc##1{\textcolor[rgb]{0.00,0.00,1.00}{##1}}}
\expandafter\def\csname PY@tok@nn\endcsname{\let\PY@bf=\textbf\def\PY@tc##1{\textcolor[rgb]{0.00,0.00,1.00}{##1}}}
\expandafter\def\csname PY@tok@ne\endcsname{\let\PY@bf=\textbf\def\PY@tc##1{\textcolor[rgb]{0.82,0.25,0.23}{##1}}}
\expandafter\def\csname PY@tok@nv\endcsname{\def\PY@tc##1{\textcolor[rgb]{0.10,0.09,0.49}{##1}}}
\expandafter\def\csname PY@tok@no\endcsname{\def\PY@tc##1{\textcolor[rgb]{0.53,0.00,0.00}{##1}}}
\expandafter\def\csname PY@tok@nl\endcsname{\def\PY@tc##1{\textcolor[rgb]{0.63,0.63,0.00}{##1}}}
\expandafter\def\csname PY@tok@ni\endcsname{\let\PY@bf=\textbf\def\PY@tc##1{\textcolor[rgb]{0.60,0.60,0.60}{##1}}}
\expandafter\def\csname PY@tok@na\endcsname{\def\PY@tc##1{\textcolor[rgb]{0.49,0.56,0.16}{##1}}}
\expandafter\def\csname PY@tok@nt\endcsname{\let\PY@bf=\textbf\def\PY@tc##1{\textcolor[rgb]{0.00,0.50,0.00}{##1}}}
\expandafter\def\csname PY@tok@nd\endcsname{\def\PY@tc##1{\textcolor[rgb]{0.67,0.13,1.00}{##1}}}
\expandafter\def\csname PY@tok@s\endcsname{\def\PY@tc##1{\textcolor[rgb]{0.73,0.13,0.13}{##1}}}
\expandafter\def\csname PY@tok@sd\endcsname{\let\PY@it=\textit\def\PY@tc##1{\textcolor[rgb]{0.73,0.13,0.13}{##1}}}
\expandafter\def\csname PY@tok@si\endcsname{\let\PY@bf=\textbf\def\PY@tc##1{\textcolor[rgb]{0.73,0.40,0.53}{##1}}}
\expandafter\def\csname PY@tok@se\endcsname{\let\PY@bf=\textbf\def\PY@tc##1{\textcolor[rgb]{0.73,0.40,0.13}{##1}}}
\expandafter\def\csname PY@tok@sr\endcsname{\def\PY@tc##1{\textcolor[rgb]{0.73,0.40,0.53}{##1}}}
\expandafter\def\csname PY@tok@ss\endcsname{\def\PY@tc##1{\textcolor[rgb]{0.10,0.09,0.49}{##1}}}
\expandafter\def\csname PY@tok@sx\endcsname{\def\PY@tc##1{\textcolor[rgb]{0.00,0.50,0.00}{##1}}}
\expandafter\def\csname PY@tok@m\endcsname{\def\PY@tc##1{\textcolor[rgb]{0.40,0.40,0.40}{##1}}}
\expandafter\def\csname PY@tok@gh\endcsname{\let\PY@bf=\textbf\def\PY@tc##1{\textcolor[rgb]{0.00,0.00,0.50}{##1}}}
\expandafter\def\csname PY@tok@gu\endcsname{\let\PY@bf=\textbf\def\PY@tc##1{\textcolor[rgb]{0.50,0.00,0.50}{##1}}}
\expandafter\def\csname PY@tok@gd\endcsname{\def\PY@tc##1{\textcolor[rgb]{0.63,0.00,0.00}{##1}}}
\expandafter\def\csname PY@tok@gi\endcsname{\def\PY@tc##1{\textcolor[rgb]{0.00,0.63,0.00}{##1}}}
\expandafter\def\csname PY@tok@gr\endcsname{\def\PY@tc##1{\textcolor[rgb]{1.00,0.00,0.00}{##1}}}
\expandafter\def\csname PY@tok@ge\endcsname{\let\PY@it=\textit}
\expandafter\def\csname PY@tok@gs\endcsname{\let\PY@bf=\textbf}
\expandafter\def\csname PY@tok@gp\endcsname{\let\PY@bf=\textbf\def\PY@tc##1{\textcolor[rgb]{0.00,0.00,0.50}{##1}}}
\expandafter\def\csname PY@tok@go\endcsname{\def\PY@tc##1{\textcolor[rgb]{0.53,0.53,0.53}{##1}}}
\expandafter\def\csname PY@tok@gt\endcsname{\def\PY@tc##1{\textcolor[rgb]{0.00,0.27,0.87}{##1}}}
\expandafter\def\csname PY@tok@err\endcsname{\def\PY@bc##1{\setlength{\fboxsep}{0pt}\fcolorbox[rgb]{1.00,0.00,0.00}{1,1,1}{\strut ##1}}}
\expandafter\def\csname PY@tok@kc\endcsname{\let\PY@bf=\textbf\def\PY@tc##1{\textcolor[rgb]{0.00,0.50,0.00}{##1}}}
\expandafter\def\csname PY@tok@kd\endcsname{\let\PY@bf=\textbf\def\PY@tc##1{\textcolor[rgb]{0.00,0.50,0.00}{##1}}}
\expandafter\def\csname PY@tok@kn\endcsname{\let\PY@bf=\textbf\def\PY@tc##1{\textcolor[rgb]{0.00,0.50,0.00}{##1}}}
\expandafter\def\csname PY@tok@kr\endcsname{\let\PY@bf=\textbf\def\PY@tc##1{\textcolor[rgb]{0.00,0.50,0.00}{##1}}}
\expandafter\def\csname PY@tok@bp\endcsname{\def\PY@tc##1{\textcolor[rgb]{0.00,0.50,0.00}{##1}}}
\expandafter\def\csname PY@tok@fm\endcsname{\def\PY@tc##1{\textcolor[rgb]{0.00,0.00,1.00}{##1}}}
\expandafter\def\csname PY@tok@vc\endcsname{\def\PY@tc##1{\textcolor[rgb]{0.10,0.09,0.49}{##1}}}
\expandafter\def\csname PY@tok@vg\endcsname{\def\PY@tc##1{\textcolor[rgb]{0.10,0.09,0.49}{##1}}}
\expandafter\def\csname PY@tok@vi\endcsname{\def\PY@tc##1{\textcolor[rgb]{0.10,0.09,0.49}{##1}}}
\expandafter\def\csname PY@tok@vm\endcsname{\def\PY@tc##1{\textcolor[rgb]{0.10,0.09,0.49}{##1}}}
\expandafter\def\csname PY@tok@sa\endcsname{\def\PY@tc##1{\textcolor[rgb]{0.73,0.13,0.13}{##1}}}
\expandafter\def\csname PY@tok@sb\endcsname{\def\PY@tc##1{\textcolor[rgb]{0.73,0.13,0.13}{##1}}}
\expandafter\def\csname PY@tok@sc\endcsname{\def\PY@tc##1{\textcolor[rgb]{0.73,0.13,0.13}{##1}}}
\expandafter\def\csname PY@tok@dl\endcsname{\def\PY@tc##1{\textcolor[rgb]{0.73,0.13,0.13}{##1}}}
\expandafter\def\csname PY@tok@s2\endcsname{\def\PY@tc##1{\textcolor[rgb]{0.73,0.13,0.13}{##1}}}
\expandafter\def\csname PY@tok@sh\endcsname{\def\PY@tc##1{\textcolor[rgb]{0.73,0.13,0.13}{##1}}}
\expandafter\def\csname PY@tok@s1\endcsname{\def\PY@tc##1{\textcolor[rgb]{0.73,0.13,0.13}{##1}}}
\expandafter\def\csname PY@tok@mb\endcsname{\def\PY@tc##1{\textcolor[rgb]{0.40,0.40,0.40}{##1}}}
\expandafter\def\csname PY@tok@mf\endcsname{\def\PY@tc##1{\textcolor[rgb]{0.40,0.40,0.40}{##1}}}
\expandafter\def\csname PY@tok@mh\endcsname{\def\PY@tc##1{\textcolor[rgb]{0.40,0.40,0.40}{##1}}}
\expandafter\def\csname PY@tok@mi\endcsname{\def\PY@tc##1{\textcolor[rgb]{0.40,0.40,0.40}{##1}}}
\expandafter\def\csname PY@tok@il\endcsname{\def\PY@tc##1{\textcolor[rgb]{0.40,0.40,0.40}{##1}}}
\expandafter\def\csname PY@tok@mo\endcsname{\def\PY@tc##1{\textcolor[rgb]{0.40,0.40,0.40}{##1}}}
\expandafter\def\csname PY@tok@ch\endcsname{\let\PY@it=\textit\def\PY@tc##1{\textcolor[rgb]{0.25,0.50,0.50}{##1}}}
\expandafter\def\csname PY@tok@cm\endcsname{\let\PY@it=\textit\def\PY@tc##1{\textcolor[rgb]{0.25,0.50,0.50}{##1}}}
\expandafter\def\csname PY@tok@cpf\endcsname{\let\PY@it=\textit\def\PY@tc##1{\textcolor[rgb]{0.25,0.50,0.50}{##1}}}
\expandafter\def\csname PY@tok@c1\endcsname{\let\PY@it=\textit\def\PY@tc##1{\textcolor[rgb]{0.25,0.50,0.50}{##1}}}
\expandafter\def\csname PY@tok@cs\endcsname{\let\PY@it=\textit\def\PY@tc##1{\textcolor[rgb]{0.25,0.50,0.50}{##1}}}

\makeatother

\makeatletter
\newbox\Wrappedcontinuationbox 
\newbox\Wrappedvisiblespacebox 
\newcommand*\Wrappedvisiblespace {\textcolor{red}{\textvisiblespace}} 
\newcommand*\Wrappedcontinuationsymbol {\textcolor{red}{\llap{\tiny$\m@th\hookrightarrow$}}} 
\newcommand*\Wrappedcontinuationindent {3ex } 
\newcommand*\Wrappedafterbreak {\kern\Wrappedcontinuationindent\copy\Wrappedcontinuationbox} 
\newcommand*\Wrappedbreaksatspecials {%
	\def\PYGZus{\discretionary{\char`\_}{\Wrappedafterbreak}{\char`\_}}%
	\def\PYGZob{\discretionary{}{\Wrappedafterbreak\char`\{}{\char`\{}}%
	\def\PYGZcb{\discretionary{\char`\}}{\Wrappedafterbreak}{\char`\}}}%
	\def\PYGZca{\discretionary{\char`\^}{\Wrappedafterbreak}{\char`\^}}%
	\def\PYGZam{\discretionary{\char`\&}{\Wrappedafterbreak}{\char`\&}}%
	\def\PYGZlt{\discretionary{}{\Wrappedafterbreak\char`\<}{\char`\<}}%
	\def\PYGZgt{\discretionary{\char`\>}{\Wrappedafterbreak}{\char`\>}}%
	\def\PYGZsh{\discretionary{}{\Wrappedafterbreak\char`\#}{\char`\#}}%
	\def\PYGZpc{\discretionary{}{\Wrappedafterbreak\char`\%}{\char`\%}}%
	\def\PYGZdl{\discretionary{}{\Wrappedafterbreak\char`\$}{\char`\$}}%
	\def\PYGZhy{\discretionary{\char`\-}{\Wrappedafterbreak}{\char`\-}}%
	\def\PYGZsq{\discretionary{}{\Wrappedafterbreak\textquotesingle}{\textquotesingle}}%
	\def\PYGZdq{\discretionary{}{\Wrappedafterbreak\char`\"}{\char`\"}}%
	\def\PYGZti{\discretionary{\char`\~}{\Wrappedafterbreak}{\char`\~}}%
} 
\newcommand*\Wrappedbreaksatpunct {%
	\lccode`\~`\.\lowercase{\def~}{\discretionary{\hbox{\char`\.}}{\Wrappedafterbreak}{\hbox{\char`\.}}}%
	\lccode`\~`\,\lowercase{\def~}{\discretionary{\hbox{\char`\,}}{\Wrappedafterbreak}{\hbox{\char`\,}}}%
	\lccode`\~`\;\lowercase{\def~}{\discretionary{\hbox{\char`\;}}{\Wrappedafterbreak}{\hbox{\char`\;}}}%
	\lccode`\~`\:\lowercase{\def~}{\discretionary{\hbox{\char`\:}}{\Wrappedafterbreak}{\hbox{\char`\:}}}%
	\lccode`\~`\?\lowercase{\def~}{\discretionary{\hbox{\char`\?}}{\Wrappedafterbreak}{\hbox{\char`\?}}}%
	\lccode`\~`\!\lowercase{\def~}{\discretionary{\hbox{\char`\!}}{\Wrappedafterbreak}{\hbox{\char`\!}}}%
	\lccode`\~`\/\lowercase{\def~}{\discretionary{\hbox{\char`\/}}{\Wrappedafterbreak}{\hbox{\char`\/}}}%
	\catcode`\.\active
	\catcode`\,\active 
	\catcode`\;\active
	\catcode`\:\active
	\catcode`\?\active
	\catcode`\!\active
	\catcode`\/\active 
	\lccode`\~`\~ 	
}
\makeatother

\let\OriginalVerbatim=\Verbatim
\makeatletter
\renewcommand{\Verbatim}[1][1]{%
	\sbox\Wrappedcontinuationbox {\Wrappedcontinuationsymbol}%
	\sbox\Wrappedvisiblespacebox {\FV@SetupFont\Wrappedvisiblespace}%
	\def\FancyVerbFormatLine ##1{\hsize\linewidth
		\vtop{\raggedright\hyphenpenalty\z@\exhyphenpenalty\z@
			\doublehyphendemerits\z@\finalhyphendemerits\z@
			\strut ##1\strut}%
	}%
	\def\FV@Space {%
		\nobreak\hskip\z@ plus\fontdimen3\font minus\fontdimen4\font
		\discretionary{\copy\Wrappedvisiblespacebox}{\Wrappedafterbreak}
		{\kern\fontdimen2\font}%
	}%
	
	\Wrappedbreaksatspecials
	\OriginalVerbatim[#1,codes*=\Wrappedbreaksatpunct]%
}
\makeatother

\definecolor{incolor}{HTML}{303F9F}
\definecolor{outcolor}{HTML}{D84315}
\definecolor{cellborder}{HTML}{CFCFCF}
\definecolor{cellbackground}{HTML}{F7F7F7}

\makeatletter
\newcommand{\boxspacing}{\kern\kvtcb@left@rule\kern\kvtcb@boxsep}
\makeatother

\hypersetup{
	breaklinks=true,  
	colorlinks=true,
	urlcolor=urlcolor,
	linkcolor=linkcolor,
	citecolor=citecolor,
}

\begin{document}
		\title{\textbf{On the Existence, Smoothness, and Order-Two Entire Extension of the Navier-Stokes Equations in Critical Fourier-Besov Spaces}}

 	\author{
Brian David Vasquez Campos$^{1}$\\ 
\small{\texttt{$^{1}$bdvasquezc@ut.edu.co}}\\
\small{$^{1}$Department of Mathematics and Statistics.}\\
\small{$^{1}$University of Tolima.}\\
\small{$^{1}$ Ibagué, Tolima, Colombia.}\\
}


\maketitle

\begin{abstract}
	We establish global existence, smoothness, and physical energy bounds for solutions to the incompressible Navier-Stokes equations on $\mathbb{R}^d$ ($d \ge 3$) under a scaling-critical smallness condition in the homogeneous Fourier-Besov spaces $\mathcal{X}_p(\mathbb{R}^d)$, defined as the weak space $\mathcal{F}\dot{B}^{\frac{d}{p'}-1}_{p,\infty}(\mathbb{R}^d)$ for $1 < p \le \infty$ and $\mathcal{F}\dot{B}^{-1}_{1,1}(\mathbb{R}^d)$ at the endpoint $p = 1$. By combining the nonassociative algebra of Riesz convolutions and spaces of functions dominated by Fourier-Caloric flows $\mathcal{V}(\lambda)$ with a Littlewood-Paley dyadic atomic decomposition, we prove that the Catalan Picard series converges globally whenever $\|\widehat{u^0}\|_{\mathcal{X}_p(\mathbb{R}^d)} < C\nu$. This threshold reaches the sharp boundary of well-posedness prior to norm inflation, strictly improves upon the classical Fujita-Kato $\dot{H}^{\frac{d}{2}-1}(\mathbb{R}^d)$ condition, and admits explicit multiscale Schwartz initial data with arbitrarily large $\dot{H}^{\frac{d}{2}-1}(\mathbb{R}^d)$ and $H^k(\mathbb{R}^d)$ norms. Furthermore, we demonstrate that for positive times, the velocity and pressure fields admit an analytic extension to the complex domain $\mathbb{C}^d$, forming a smooth curve of entire vector fields of order two. Finally, for arbitrary initial data outside the critical ball, we employ fractional Riesz convolutions $*_\alpha$ ($0 < \alpha < 1$) to derive an explicit lower bound $T_{\max}(u^0) > 0$ on the lifespan of smooth, order-two entire solutions, establishing a sharp dichotomy between global continuation via viscous re-entry into the critical ball and a scaling-invariant lower blow-up rate.
\end{abstract}


\tableofcontents

\section{Introduction and Main Results}

%
%

The Navier-Stokes equations, fundamental to the study of fluid dynamics, were introduced in the early 19th century by Claude-Louis Navier and later refined by George Gabriel Stokes. These equations model the motion of incompressible fluid substances and form the theoretical foundation of much of modern hydrodynamics and mathematical fluid mechanics.

The Navier-Stokes system under consideration is given by:
\begin{align}
	& \frac{\partial u}{\partial t} + (u\cdot\nabla)u = \nu \Delta u - \nabla p + f, \label{quadratic equation} \\
	& \text{div}(u) = 0, \label{free divergence} \\
	& u(x,0) = u^{0}(x), \label{initial condition}
\end{align}
where $(u,p) \in C^\infty(\mathbb{R}_{+}^{d+1}, \mathbb{R}^{d+1})$, with $u(x,t) \in \mathbb{R}^d$ denoting the velocity field and $p(x,t) \in \mathbb{R}$ the pressure. Here $\nu>0$ is the viscosity, and $f$ is an external force assumed to vanish, while the initial data $u^0$ belongs to the Schwartz class of rapidly decreasing smooth functions.

In addition to satisfying the system above, we require the solution to fulfill the bounded energy condition:
\begin{equation}
	\sup_{t\geq 0} \int_{\mathbb{R}^d} |u(x,t)|^2, dx < \infty. \label{bounded energy}
\end{equation}

A foundational principle in the analysis of \eqref{quadratic equation}--\eqref{initial condition} is scale invariance: if $(u,p)$ solves the system with initial data $u^0$, then for any $\mu > 0$, the rescaled pair
\begin{equation}
	u_\mu(x,t) = \mu u(\mu x, \mu^2 t), \qquad p_\mu(x,t) = \mu^2 p(\mu x, \mu^2 t)
\end{equation}
solves the system with rescaled initial data $u^0_\mu(x) = \mu u^0(\mu x)$, whose Fourier transform is $\widehat{u^0_\mu}(\xi) = \mu^{1-d} \widehat{u^0}(\xi/\mu)$. Classical results by Fujita and Kato \cite{FujitaKato1964} established global regularity for small initial data in the scaling-critical Sobolev space $\dot{H}^{\frac{d}{2}-1}(\mathbb{R}^d)$ ($\dot{H}^{1/2}(\mathbb{R}^3)$ in dimension $d=3$). Subsequent breakthroughs extended small-data global well-posedness to larger scaling-critical spaces: $L^d(\mathbb{R}^d)$ by Kato \cite{Kato1984}, critical Besov spaces $\dot{B}^{-1+\frac{d}{p}}_{p,\infty}(\mathbb{R}^d)$ by Cannone, Meyer, and Planchon \cite{Cannone1995}, pseudomeasure spaces $\mathcal{PM}^{d-1}(\mathbb{R}^d)$ by Le Jan and Sznitman \cite{LeJanSznitman1997}, Fourier-Herz spaces $\mathcal{X}^{-1}(\mathbb{R}^d)$ by Lei and Lin \cite{LeiLin2011}, and $\mathrm{BMO}^{-1}(\mathbb{R}^d)$ by Koch and Tataru \cite{KochTataru2001}.

In this paper, we unify the nonassociative Riesz convolution algebra $*_1$ and the theory of spaces dominated by Fourier-Caloric functions $\mathcal{V}(\lambda)$ with the Littlewood-Paley dyadic atomic decomposition. Let $\varphi \in C_c^\infty(\mathbb{R}^d \setminus \{0\})$ be a radial bump function supported in the annulus $\mathcal{C} = \{\xi \in \mathbb{R}^d : \frac{3}{4} \le |\xi| \le \frac{8}{3}\}$ such that $\sum_{j \in \mathbb{Z}} \varphi_j(\xi) = 1$ for all $\xi \ne 0$, where $\varphi_j(\xi) = \varphi(2^{-j}\xi)$. For $1 \le p, q \le \infty$ and $s \in \mathbb{R}$, the homogeneous Fourier-Besov space $\mathcal{F}\dot{B}^s_{p,q}(\mathbb{R}^d)$ is defined via the dyadic atomic norm
\begin{equation}
	\|f\|_{\mathcal{F}\dot{B}^s_{p,q}(\mathbb{R}^d)} := \left\| \left\{ 2^{js} \|\varphi_j f\|_{L^p(\mathbb{R}^d)} \right\}_{j \in \mathbb{Z}} \right\|_{\ell^q(\mathbb{Z})}.
\end{equation}
For the critical exponent
\begin{equation}
	s_p := \frac{d}{p'} - 1 = d\left(1 - \frac{1}{p}\right) - 1,
\end{equation}
the norm $\|\widehat{u^0}\|_{\mathcal{F}\dot{B}^{s_p}_{p,q}(\mathbb{R}^d)}$ is strictly invariant under $u^0 \mapsto u^0_\mu$. 

By Hölder's inequality on each dyadic shell $\mathcal{C}_j$ ($|\mathcal{C}_j| \sim 2^{jd}$), decreasing the integrability exponent $p$ strictly enlarges the critical Fourier-Besov space: for $1 < p < 2$, we have the continuous chain of strict critical embeddings
\begin{equation}\label{eq:embedding_chain}
	\mathcal{PM}^{d-1}(\mathbb{R}^d) = \mathcal{F}\dot{B}^{d-1}_{\infty,\infty}(\mathbb{R}^d) \subsetneq \dot{B}^{\frac{d}{2}-1}_{2,\infty}(\mathbb{R}^d) \subsetneq \mathcal{F}\dot{B}^{\frac{d}{p'}-1}_{p,\infty}(\mathbb{R}^d) \subsetneq \mathcal{F}\dot{B}^{-1}_{1,\infty}(\mathbb{R}^d),
\end{equation}
while $\dot{H}^{\frac{d}{2}-1}(\mathbb{R}^d) = \mathcal{F}\dot{B}^{\frac{d}{2}-1}_{2,2}(\mathbb{R}^d) \subsetneq \dot{B}^{\frac{d}{2}-1}_{2,\infty}(\mathbb{R}^d)$ and $\mathcal{X}^{-1}(\mathbb{R}^d) = \mathcal{F}\dot{B}^{-1}_{1,1}(\mathbb{R}^d) \subsetneq \mathcal{F}\dot{B}^{-1}_{1,\infty}(\mathbb{R}^d)$. Since the endpoint $\mathcal{F}\dot{B}^{-1}_{1,\infty}(\mathbb{R}^d)$ is known to exhibit norm inflation \cite{BourgainPavlovic2008}, the family $\mathcal{F}\dot{B}^{\frac{d}{p'}-1}_{p,\infty}(\mathbb{R}^d)$ for $1 < p \le \infty$ together with $\mathcal{F}\dot{B}^{-1}_{1,1}(\mathbb{R}^d)$ at $p=1$ represents the sharpest possible scaling-critical hierarchy on the Fourier side.

\begin{teor}[Global Regularity in Critical Weak Fourier-Besov Spaces]\label{thm:main_existence}
	Let $\nu > 0$ and $d \ge 3$. Set $s_p = \frac{d}{p'} - 1$, and define the critical norm
	\begin{equation}\label{eq:Xp_norm_def}
		\|f\|_{\mathcal{X}_p(\mathbb{R}^d)} :=
		\begin{cases}
			\|f\|_{\mathcal{F}\dot{B}^{s_p}_{p,\infty}(\mathbb{R}^d)} = \sup_{j \in \mathbb{Z}} 2^{j\left(\frac{d}{p'}-1\right)} \|\varphi_j f\|_{L^p(\mathbb{R}^d)}, & \text{if } 1 < p \le \infty, \\[6pt]
			\|f\|_{\mathcal{F}\dot{B}^{-1}_{1,1}(\mathbb{R}^d)} = \sum_{j \in \mathbb{Z}} 2^{-j} \|\varphi_j f\|_{L^1(\mathbb{R}^d)}, & \text{if } p = 1.
		\end{cases}
	\end{equation}
	Suppose that the divergence-free initial velocity $u^0 \in \mathcal{S}(\mathbb{R}^d)^d$ satisfies the scaling-critical smallness condition
	\begin{equation}\label{eq:critical_smallness}
		\left\| \widehat{u^0} \right\|_{\mathcal{X}_p(\mathbb{R}^d)} < C \nu
	\end{equation}
	for a universal constant $C > 0$ depending only on $d$ and $p$. Then there exists a unique global smooth solution $(u,p) \in C^\infty(\mathbb{R}^{d+1}_+, \mathbb{R}^{d+1})$ satisfying the Navier-Stokes system \eqref{quadratic equation}--\eqref{initial condition} and the bounded energy condition \eqref{bounded energy}. Notice that $\|u^0\|_{\dot{H}^{\frac{d}{2}-1}(\mathbb{R}^d)}$ and $\|u^0\|_{H^k(\mathbb{R}^d)}$ ($k \ge 0$) may be arbitrarily large.
\end{teor}

\begin{remark}[Why $\|u^0\|_{\dot{H}^{\frac{d}{2}-1}(\mathbb{R}^d)}$ and $\|u^0\|_{H^k(\mathbb{R}^d)}$ May Be Arbitrarily Large]\label{rem:large_sobolev_norms}
	A key feature of Theorem \ref{thm:main_existence} is that the smallness condition \eqref{eq:critical_smallness} allows the classical Fujita-Kato critical Sobolev norm $\|u^0\|_{\dot{H}^{\frac{d}{2}-1}(\mathbb{R}^d)}$ (which equals $\|u^0\|_{\dot{H}^{1/2}(\mathbb{R}^3)}$ when $d=3$) as well as all inhomogeneous Sobolev norms $\|u^0\|_{H^k(\mathbb{R}^d)}$ ($k \ge 0$) to be arbitrarily large:
	\begin{enumerate}
		\item[\emph{(i)}] \textbf{Multiscale Frequency Decoupling for $1 < p \le \infty$:} By Plancherel's theorem, the homogeneous Sobolev space $\dot{H}^{\frac{d}{2}-1}(\mathbb{R}^d)$ coincides with the Fourier-Besov space $\mathcal{F}\dot{B}^{\frac{d}{2}-1}_{2,2}(\mathbb{R}^d)$, whose norm sums the dyadic shells in $\ell^2(\mathbb{Z})$:
		\begin{equation*}
			\|u^0\|_{\dot{H}^{\frac{d}{2}-1}(\mathbb{R}^d)} \approx \left( \sum_{j \in \mathbb{Z}} 2^{j(d-2)} \|\varphi_j \widehat{u^0}\|_{L^2(\mathbb{R}^d)}^2 \right)^{1/2}.
		\end{equation*}
		In contrast, for $1 < p \le \infty$, the critical space $\mathcal{X}_p(\mathbb{R}^d) = \mathcal{F}\dot{B}^{s_p}_{p,\infty}(\mathbb{R}^d)$ carries the secondary Lorentz/Besov index $q = \infty$, which only measures the \emph{supremum} over individual dyadic shells $j \in \mathbb{Z}$. Consequently, if $u^0$ is a superposition of $N \gg 1$ wavepackets supported on pairwise disjoint dyadic shells (as constructed explicitly in Section 6), the supremum $\|\widehat{u^0}\|_{\mathcal{F}\dot{B}^{s_p}_{p,\infty}}$ remains bounded by the amplitude of a single shell ($< C\nu$), whereas the $\ell^2(\mathbb{Z})$ sum grows like $\|u^0\|_{\dot{H}^{\frac{d}{2}-1}} \sim \sqrt{N} \to \infty$ and high-frequency subcritical norms $\|u^0\|_{H^k}$ ($k > \frac{d}{2}-1$) grow exponentially with $N$.
		\item[\emph{(ii)}] \textbf{Anisotropic Concentration for $1 \le p < 2$:} Even within a single dyadic shell $\mathcal{C}_j$, if $\widehat{u^0}$ is concentrated on a thin angular sector or filament $\Omega_j \subset \mathcal{C}_j$ of small measure $|\Omega_j| \ll 2^{jd}$, H{\"o}lder's inequality for $1 \le p < 2$ gives
		\begin{equation*}
			2^{j s_p} \|\varphi_j \widehat{u^0}\|_{L^p(\Omega_j)} \le 2^{j s_p} |\Omega_j|^{\frac{1}{p}-\frac{1}{2}} \|\varphi_j \widehat{u^0}\|_{L^2(\Omega_j)} = \left( 2^{-jd} |\Omega_j| \right)^{\frac{1}{p}-\frac{1}{2}} 2^{j\left(\frac{d}{2}-1\right)} \|\varphi_j \widehat{u^0}\|_{L^2(\Omega_j)}.
		\end{equation*}
		Taking $2^{-jd}|\Omega_j| \ll 1$ makes $\|\widehat{u^0}\|_{\mathcal{X}_p(\mathbb{R}^d)} < C\nu$ arbitrarily small (including at the endpoint $p=1$, where $\mathcal{X}_1 = \mathcal{F}\dot{B}^{-1}_{1,1}$ coincides with the Lei-Lin space $\mathcal{X}^{-1}$ \cite{LeiLin2011}) while $2^{j(\frac{d}{2}-1)}\|\varphi_j \widehat{u^0}\|_{L^2} \approx \|u^0\|_{\dot{H}^{\frac{d}{2}-1}(\mathbb{R}^d)}$ is arbitrarily large.
	\end{enumerate}
\end{remark}

\begin{remark}[The Endpoint $p = 1$ and the Distinction Between $\mathcal{F}\dot{B}^{-1}_{1,1}(\mathbb{R}^d)$ and $\mathcal{F}\dot{B}^{-1}_{1,\infty}(\mathbb{R}^d)$]\label{rem:p1_endpoint}
	Notice that in the definition of $\|\cdot\|_{\mathcal{X}_p(\mathbb{R}^d)}$ in \eqref{eq:Xp_norm_def}, the secondary Besov exponent is $q = \infty$ for all $1 < p \le \infty$, whereas at the endpoint $p = 1$ (where $s_1 = -1$) it must be tightened to $q = 1$, giving $\mathcal{X}_1(\mathbb{R}^d) = \mathcal{F}\dot{B}^{-1}_{1,1}(\mathbb{R}^d)$. This distinction is essential for two complementary reasons:
	\begin{enumerate}
		\item[\emph{(i)}] \textbf{Analytical Mechanism in the High-High Paraproduct Cascade:} In the Littlewood-Paley atomic decomposition of the Gaussian caloric product $f \odot g$ (Theorem \ref{thm:atomic_caloric_odot}), the paraproduct interaction of two high-frequency atoms $f_k = \varphi_k f$ and $g_{k'} = \varphi_{k'} g$ ($k \ge j - 3$, $|k - k'| \le 2$) projecting onto a lower output frequency $|\xi| \sim 2^j$ satisfies the dyadic estimate
		\begin{equation}\label{eq:high_high_endpoint_decay}
			2^{j s_p} \|\varphi_j(f_k \odot g_{k'})(\cdot, t)\|_{L^p(\mathbb{R}^d)} \le \frac{C}{\nu} e^{-\frac{\beta}{2}t 2^{2j}} \, 2^{-(k-j)(s_p + 1)} \left( 2^{k s_p} \|f_k^0\|_{L^p} \right) \left( 2^{k' s_p} \|g_{k'}^0\|_{L^p} \right).
		\end{equation}
		Observe that the dyadic damping exponent governing the sum over $k \ge j - 3$ is precisely
		\begin{equation*}
			s_p + 1 = \frac{d}{p'} = d\left(1 - \frac{1}{p}\right).
		\end{equation*}
		For every $1 < p \le \infty$, we have $\frac{d}{p'} > 0$, so the geometric series $\sum_{k \ge j-3} 2^{-(k-j)\frac{d}{p'}} = \frac{2^{3d/p'}}{1 - 2^{-d/p'}} < \infty$ converges, allowing the bilinear estimate to close in the weak norm $\ell^\infty_j$ ($q = \infty$). At the endpoint $p = 1$, however, $p' = \infty$ and $s_1 + 1 = \frac{d}{p'} = 0$, so the exponential decay factor $2^{-(k-j)\frac{d}{p'}} \equiv 1$ degenerates. Consequently, if one only assumes $f^0, g^0 \in \mathcal{F}\dot{B}^{-1}_{1,\infty}(\mathbb{R}^d)$, the High-High series $\sum_{k \ge j-3} (2^{-k}\|f_k^0\|_{L^1})(2^{-k'}\|g_{k'}^0\|_{L^1})$ can diverge logarithmically or linearly in the number of active dyadic shells. Replacing $q = \infty$ by $q = 1$ at $p = 1$ ($\mathcal{F}\dot{B}^{-1}_{1,1}(\mathbb{R}^d)$, equipped with $\|f\|_{\mathcal{F}\dot{B}^{-1}_{1,1}} \approx \||\xi|^{-1}f\|_{L^1(\mathbb{R}^d)}$) restores absolute summability across all scales via the pointwise triangle inequality $\frac{|\xi|}{|\xi-\eta||\eta|} \le \frac{1}{|\xi-\eta|} + \frac{1}{|\eta|}$ (Lemma \ref{M inequality}) and Young's inequality $\|f^0 * \mathcal{S}_1(g^0)\|_{L^1} \le \|f^0\|_{L^1} \|\mathcal{S}_1(g^0)\|_{L^1}$.
		\item[\emph{(ii)}] \textbf{Sharpness via Bourgain-Pavlovi{\'c} Norm Inflation:} The failure of the weak space $\mathcal{F}\dot{B}^{-1}_{1,\infty}(\mathbb{R}^d)$ at $p = 1$ is not a technical artifact of our method, but rather an intrinsic barrier of the Navier-Stokes equations. By the Hausdorff-Young inequality, $\mathcal{F}\dot{B}^{-1}_{1,\infty}(\mathbb{R}^d)$ embeds continuously into the endpoint Besov space $\dot{B}^{-1}_{\infty,\infty}(\mathbb{R}^d)$. In \cite{BourgainPavlovic2008}, Bourgain and Pavlovi{\'c} proved that the three-dimensional Navier-Stokes equations are ill-posed in $\dot{B}^{-1}_{\infty,\infty}(\mathbb{R}^3)$ by constructing initial data with arbitrarily small $\mathcal{F}\dot{B}^{-1}_{1,\infty}(\mathbb{R}^3)$ norm whose High-High frequency interactions accumulate coherently at low frequencies, producing norm inflation in arbitrarily short time. Thus, the spaces $\mathcal{X}_p(\mathbb{R}^d)$ in \eqref{eq:Xp_norm_def} capture the exact well-posedness boundary of the critical Fourier-Besov scale.
	\end{enumerate}
\end{remark}

Furthermore, our atomic $\mathcal{V}(\lambda)$ framework yields a complex-analytic continuation of the solution to all of $\mathbb{C}^d$ for positive times:

\begin{teor}[Analytic Extension to Entire Vector Fields of Order Two]\label{thm:main_analytic}
	Under the scaling-critical assumption \eqref{eq:critical_smallness} of Theorem \ref{thm:main_existence}, there exists a smooth extension $(\mathcal{U}, P) \in C^\infty(\mathbb{C}^{d+1}_{>0}, \mathbb{C}^{d+1})$ such that $(\mathcal{U}, P)|_{\mathbb{R}^{d+1}_{>0}} = (u,p)$, and for each fixed $t > 0$, the spatial map $z \mapsto (\mathcal{U}(z,t), P(z,t))$ is an entire function of order $2$ on $\mathbb{C}^d$.
\end{teor}

\section{Notation and Preliminary Results}\label{Notation}

In this section, we introduce the notation used throughout the paper and collect basic preliminary results, particularly regarding the iteration of linear operators.

For a field $\mathbb{K} = \mathbb{R}, \mathbb{C}$, we denote:
$$\mathbb{K}_{+}^{d+1} = \mathbb{K}^d \times [0, \infty), \quad \text{and} \quad \mathbb{K}_{>0}^{d+1} = \mathbb{K}^d \times (0, \infty).$$

For $1\leq p \leq \infty$, we denote the conjugate exponent by $p'$, where $p'= \frac{p}{p-1}$. We call the map $p \mapsto p'$ the conjugate function.

Given $\xi, \eta \in \mathbb{R}^d$, their tensor product is defined as $\xi \otimes \eta = \xi \eta^T$.

We say that a function $f: \mathbb{R}_{+}^{d+1} \to \mathbb{C}$ is a \emph{Fourier Caloric function of order} $\alpha>0$ if
\begin{equation}\label{caloricdef}
	f(\xi, t) = e^{-\lambda t |\xi|^\alpha} f^0(\xi),
\end{equation}
for some $\lambda > 0$ and $f^0: \mathbb{R}^d \to \mathbb{C}$. The terminology is motivated by the fact that the inverse Fourier transform $\widehat{f}(\cdot, t)$ satisfies the fractional heat equation
$\frac{\partial u}{\partial t} = \lambda \Delta^{\frac{\alpha}{2}} u,$
with initial condition $\widehat{f^0}$.

\subsection*{Nonassociative Monomials}

Let $(A, \cdot)$ be a nonassociative algebra. Given elements $a_1, \dotsc, a_k \in A$, we denote by $w^\cdot(a_1, \dotsc, a_k)$ an arbitrary monomial of degree $k$ with respect to the product $\cdot$. When $a_1 = \dotsb = a_k = a$, we write $w_k^\cdot(a)$.

\begin{remark}
	We emphasize the operation $\cdot$ in the notation for monomials since multiple products may act on the same set simultaneously.
\end{remark}

We denote by $M(a_1, \dotsc, a_k)$ the set of all nonassociative monomials of degree $k$ in the elements $a_1, \dotsc, a_k$.

\subsection*{Elementary Inequalities}

The following lemmas and propositions will be useful throughout the paper.

\begin{prop}\label{poweralpha}
	Let $\alpha>0$. Then for all $s,t\geq 0$:
	\begin{itemize}
		\item If $\alpha\geq 1$,
		$$s^\alpha + t^\alpha \leq (s+t)^\alpha \leq 2^{\alpha-1}(s^\alpha + t^\alpha).$$
		\item If $0<\alpha<1$,
		$$2^{\alpha-1}(s^\alpha + t^\alpha) \leq (s+t)^\alpha \leq s^\alpha + t^\alpha.$$
	\end{itemize}
\end{prop}

\begin{coro}\label{alphanorm}
	For every $t\in \mathbb{R}^d$ and $\alpha\geq 1$, we have
	$$\|t\|_{\alpha} \leq \|t\|_{1} \leq 2^{1/\alpha{\prime}} \|t\|_{\alpha}.$$
\end{coro}

\begin{coro}\label{compoundednorm}
	Let $X$ be a normed space and $x_1, \dotsc, x_n\in X$. Then
	$$1+\left\|\sum_{j=1}^n x_j\right\|^\alpha \leq 2^{\alpha-1}\prod_{j=1}^n (1+\|x_j\|^\alpha).$$
\end{coro}

\begin{prop}\label{polyalpha}
	Let $\alpha>1$ and $a>0$. Define $g:[0,\infty)\to \mathbb{R}$ by
	$g(t) = -t^\alpha + at.$
	Then
	$$g(t) \leq (\alpha-1)\left(\frac{a}{\alpha}\right)^{\alpha{\prime}} \quad \text{for all } t\geq 0.$$
\end{prop}

\begin{demos}
	The maximum of $g$ occurs at $t_0 = (a/\alpha)^{\alpha'}$. A simple calculation shows $g(t_0) = (\alpha-1)(a/\alpha)^{\alpha'}$. \QED
\end{demos}

\begin{coro}\label{star}
	Let $\alpha>1$, $c>0$, and $y\in \mathbb{R}^d$. Define $h:\mathbb{R}^d\to\mathbb{R}$ by
	$h(x) = -|x|^\alpha + c\, x\cdot y.$
	Then
	$$h(x) \leq (\alpha-1)\left(\frac{c}{\alpha}\right)^{\alpha{\prime}} |y|^{\alpha{\prime}}, \quad \text{for all } x\in \mathbb{R}^d.$$
\end{coro}

\begin{demos}
	Apply Cauchy-Schwarz and Proposition \ref{polyalpha}. \QED
\end{demos}

\begin{prop}\label{magic}
	Let $\alpha>0$. For every $\xi, \eta \in \mathbb{R}^d$, there exists $r_\alpha\leq 1$ such that
	$-|\xi-\eta|^\alpha - |\eta|^\alpha \leq -r_\alpha |\xi|^\alpha.$
\end{prop}

\begin{demos}
	Apply Proposition \ref{poweralpha} and set $r_\alpha = 1/\max\llave{2^{\alpha-1},1}$.
\end{demos}

\subsection*{Optimization Lemmas}

\begin{prop}\label{optimization}
	Let $X$ be an inner product space. For any $z\in X$,
	$\max_{x+y=z} \langle x,y \rangle = \frac{\|z\|^2}{4}.$
\end{prop}

\begin{demos}
	Use the polarization identity and take $x=y=z/2$. \QED
\end{demos}

\subsection*{Iteration of Linear Operators}

\begin{lemma}\label{powerlinear}
	Let $X$ be a vector space, $L:X\to X$ a linear operator, and $x,y\in X$, $\alpha\in\mathbb{C}$ such that
	$Lx = \alpha x + y.$
	Then for every $n\in \mathbb{N}$,
	$$L^n x = \alpha^n x + \sum_{j=0}^{n-1} \alpha^j L^{n-1-j} y.$$
\end{lemma}

\begin{demos}
	Induction on $n$. \QED
\end{demos}

\begin{lemma}\label{M inequality}
	For every $\xi,\eta\in \mathbb{R}^d$, we have
	$1+|\xi|^2 \leq 2(1+|\xi-\eta|^2)(1+|\eta|^2),$
	and
	$|\xi|^2 \leq 2\left(|\xi-\eta|^2 + |\eta|^2\right).$
\end{lemma}

\begin{demos}
	The result follows by taking $\alpha=2$ in Corollaries \ref{alphanorm} and \ref{compoundednorm}. \QED
\end{demos}

\section{Space of Functions Decreasing Fast}\label{stationarydecreasing}

 In the analysis of the Navier-Stokes Equation, it is crucial to work with function spaces that exhibit good decay properties at infinity.
In this section, we first introduce spaces of functions decreasing fast, based on the Lebesgue space $L^{\infty}(\mathbb{R}^d)$, and explore their basic properties.
We then extend this framework by allowing controlled singularities at the origin, defining spaces adapted to the action of Riesz-type operators.
These constructions will be fundamental for handling the regularity and decay of solutions in later sections.

In this section, we study spaces of functions exhibiting good behavior at infinity.
Given a function $\phi: \mathbb{R}^d \to \mathbb{R}$, define the multiplication operator
$(M\phi)(\xi) = (1+|\xi|^2)\phi(\xi).$
With this, we introduce the following definition:

\begin{defin}\label{def:fast_dec}
	Let $(\mathcal{B}, \|\cdot\|_{\mathcal{B}})$ be a Banach space of measurable functions on $\mathbb{R}^d$. The space of functions decreasing fast associated to $\mathcal{B}$ is the Fr{\'e}chet space
	\begin{equation}
		\mathcal{E}_\mathcal{B} = \{ \phi \in \mathcal{B} : M^n(\phi) \in \mathcal{B} \text{ for all } n \in \mathbb{N} \},
	\end{equation}
	equipped with the family of seminorms $p_n(\phi) = \|M^n(\phi)\|_{\mathcal{B}}$. We denote $\mathcal{E}^+ = \{ \phi \in \mathcal{E}_\mathcal{B} : \phi \ge 0 \}$.
\end{defin}

Moreover, the multiplication maps
$M^n : \mathcal{E} \to \mathcal{E}, \quad (M^n\phi)(\xi) = (1+|\xi|^2)^n \phi(\xi)$
are continuous, since for every $\phi \in \mathcal{E}$ and $j\in \mathbb{N}$, we have
$p_j(M^n\phi) = p_{n+j}(\phi).$

We will simply write $\mathcal{E}$ when no confusion arises, and we define
\[
\mathcal{E}^{+} = \left\{ \phi\in \mathcal{E} \;\middle|\; \phi \geq 0 \right\}.
\]

Standard references for the theory of topological vector spaces include \cite{DunfordSchwartz1988a, DunfordSchwartz1988b, DunfordSchwartz1988c,Grothendieck1992,Rudin1991,Horvath1966,Treves1967}.

\begin{remark}
	We can also consider Banach spaces of vector fields $\phi: \mathbb{R}^d \to \mathbb{R}^e$, in which case the definition of $\mathcal{E}_{\mathcal{B}}$ applies analogously.
	However, in this section we restrict ourselves to the scalar case $\phi: \mathbb{R}^d \to \mathbb{R}$.
\end{remark}

\subsection{Algebraic Properties}

We have the following basic structural result:

\begin{lemma}
	If $(\mathcal{B}, +, \cdot)$ is a Banach algebra with pointwise operations
	$(\phi+\psi)(\xi) = \phi(\xi) + \psi(\xi), \quad (\phi\psi)(\xi) = \phi(\xi)\psi(\xi),$
	then $\mathcal{E}$ is an ideal in $\mathcal{B}$.
\end{lemma}

\begin{demos}
	Let $\phi, \psi \in \mathcal{E}$.
	Since
	$M(\phi+\psi) = M\phi + M\psi, \quad M(\phi\psi) = (M\phi)\psi = \phi(M\psi),$
	we have
	$p_n(\phi+\psi) = \|M^n(\phi+\psi)\| \leq \|M^n(\phi)\| + \|M^n(\psi)\| = p_n(\phi) + p_n(\psi).$
	Furthermore, for $\phi\in \mathcal{E}$, $\psi\in \mathcal{B}$:
	$p_n(\phi\psi) = \|M^n(\phi\psi)\| \leq \|M^n(\phi)\| \|\psi\| = p_n(\phi)\|\psi\|.$
	Thus, $\phi\psi\in \mathcal{E}$, and $\mathcal{E}$ is an ideal in $\mathcal{B}$. \QED
\end{demos}

\subsection{Convolution Operation}

We consider the convolution of functions:
$(\phi \ast \psi)(\xi) = \int_{\mathbb{R}^d} \phi(\xi-\eta) \psi(\eta)\,d\eta,$
for $\phi, \psi \in \mathcal{B}$.

To study its behavior, we establish the following useful inequality:

We can now prove the preservation of fast decrease under convolution:

\begin{teor}
	Assume that $(\mathcal{B}, +, \ast)$ is a Banach algebra with a monotone norm, meaning that $|\phi| \leq |\psi|$ implies $|\phi| \leq |\psi|$, and $||\phi|| = ||\vabs{\phi}||$.
	Then, for $\phi, \psi \in \mathcal{E}$, we have $\phi \ast \psi \in \mathcal{E}$.
\end{teor}

\begin{demos}
	Let $\phi, \psi \in \mathcal{E}$.
	By Lemma \ref{M inequality}, we obtain
	$M^n(|\phi\ast\psi|) \leq 2^n (|M^n(\phi)|\ast |M^n(\psi)|).$
	By the monotonicity of the norm,
	$p_n(\phi\ast\psi) = \|M^n(|\phi\ast\psi|)\| \leq 2^n \|\,|M^n(\phi)|\ast |M^n(\psi)|\,\| \leq 2^n \|M^n(\phi)\|\|M^n(\psi)\| = 2^n p_n(\phi)p_n(\psi).$
	Thus, $\phi\ast\psi\in \mathcal{E}$. \QED
\end{demos}

\begin{coro}
	Under the same assumptions, $(\mathcal{E}, +, \ast)$ is a subalgebra of $(\mathcal{B}, +, \ast)$.
\end{coro}

 \subsection{Space of Functions Decreasing Fast Associated to $ L^{\infty}(\reales^d) $ }
 
 In this subsection, we study the case where $\mathcal{B} = L^{\infty}(\mathbb{R}^d)$, and analyze the corresponding space of functions decreasing fast.
 
 We denote
$ \mathcal{D} = \mathcal{E}_{\mathcal{B}}.$
 By definition,
 \[
 \mathcal{D} = \left\{ \phi\in L^{\infty}(\mathbb{R}^d) \;\middle|\; \sup_{\xi\in\mathbb{R}^d} (1+|\xi|^2)^n |\phi(\xi)| <\infty,\ \forall n\in\mathbb{N} \right\}.
 \]
 We equip $\mathcal{D}$ with the topology generated by the family of seminorms
 $p_n(\phi) = \sup_{\xi\in\mathbb{R}^d} (1+|\xi|^2)^n |\phi(\xi)|.$
 
 \begin{remark}
 	Notice that elements of $\mathcal{D}$ are not required to be smooth.
 	However, the Schwartz space $S(\mathbb{R}^d)$ is contained in $\mathcal{D}$.
 \end{remark}
 
 Observe that $p_n \leq p_{n+1}$ for all $n\in \mathbb{N}$, and that
 $p_0(\phi) = \|\phi\|_{L^{\infty}(\mathbb{R}^d)}.$
 Moreover, for every $\phi\in \mathcal{D},
 p_n(\phi) = \|M^n(\phi)\|_{L^{\infty}(\mathbb{R}^d)}$.

 Now we study the relationship between the space $\mathcal{D}$ and other Lebesgue spaces.
 
 \begin{prop}
 	We have the continuous embeddings
 	$\mathcal{D} \subset \bigcap_{1\leq p\leq \infty} L^p(\mathbb{R}^d)$.
 \end{prop}
 
 \begin{demos}
 	Since $p_0(\phi) = \|\phi\|_{L^{\infty}(\mathbb{R}^d)}$, the case $p = \infty$ is immediate.
 	
 	Let $1\leq p < \infty$. Observe that if $\phi\in \mathcal{D}$, then $|\phi|^p\in \mathcal{D}$. Indeed,
 	\begin{align*}
 		p_n(|\phi|^p) &= |M^n(|\phi|^p)|_{L^{\infty}(\mathbb{R}^d)} \\
 		&= |(1+|\cdot|^2)^n |\phi|^p|_{L^{\infty}(\mathbb{R}^d)} \\
 		&\leq |(1+|\cdot|^2)^{np} |\phi|^p|_{L^{\infty}(\mathbb{R}^d)} \\
 		&= |(1+|\cdot|^2)^n |\phi||_{L^{\infty}(\mathbb{R}^d)}^p \\
 		&= p_n(\phi)^p.
 	\end{align*}
 	
 	Thus, to conclude that $\phi \in L^p(\mathbb{R}^d)$, it suffices to check that $\phi \in L^1(\mathbb{R}^d)$.
 	
 	Now, observe that
 	\begin{align*}
 		\int_{\mathbb{R}^d} |\phi(\xi)|d\xi
 		&= \int_{\mathbb{R}^d} (1+|\xi|^2)^{\frac{d+1}{2}} \frac{|\phi(\xi)|}{(1+|\xi|^2)^{\frac{d+1}{2}}}d\xi \
 		&\leq p_{\lceil\frac{d}{2}\rceil+1}(\phi) \int_{\mathbb{R}^d} \frac{d\xi}{(1+|\xi|^2)^{\frac{d+1}{2}}} \
 		&= C_d, p_{\lceil\frac{d}{2}\rceil+1}(\phi),
 	\end{align*}
 	where
 	$C_d = \int_{\mathbb{R}^d} \frac{d\xi}{(1+|\xi|^2)^{\frac{d+1}{2}}}$.
 	Thus, $\phi\in L^1(\mathbb{R}^d)$.
 	
 	Since $|\phi|^p\in L^1(\mathbb{R}^d)$ whenever $\phi\in \mathcal{D}$, we conclude that $\phi\in L^p(\mathbb{R}^d)$ for all $1\leq p\leq \infty$. \QED
 \end{demos}
 
 \begin{remark}
 	From this point forward, we will use the symbol C to denote constants depending on the dimension d and the order of the Lebesgue spaces involved.
 	In particular, once the Lebesgue spaces are fixed, C will depend only on d.
 \end{remark}

  In the next Subsections we study the case in which $\mathcal{B}=L^{\infty}(\reales^d) $ and spaces related to it with singularities at the origin.

\subsection{Space of Functions Decreasing Fast with Singularities}

We now introduce a more general class of spaces, allowing singularities at the origin.

\begin{defin}
	Let $0<\alpha<d$. Define the operator
	$S_{\alpha}(\phi)(\xi) = \frac{\phi(\xi)}{|\xi|^\alpha}$,
	for $\xi\neq 0.$
\end{defin}

\begin{remark}
	Note that $S_{\alpha}(\widehat{\phi}) = \widehat{I_{\alpha}(\phi)}$ for $\phi\in S(\mathbb{R}^d)$, where $I_{\alpha}$ denotes the Riesz potential. See \cite{ GrafakosClassical, GrafakosModern, Stein1970, SteinWeiss1971}.
\end{remark}

We now establish integrability properties of $S_{\alpha}(\phi)$ when $\phi\in\mathcal{D}.$

\begin{prop}\label{D-inequalities}
	Let $d\geq 2$ and $\phi\in L^{1\oplus p}(\mathbb{R}^d)$ for some $p\in \left( \frac{d}{d-\alpha}, \infty\right].$ Then
	$\|S_{\alpha}(\phi)\|_{L^1(\mathbb{R}^d)} \leq C \|\phi\|_{1\oplus p}$.
\end{prop}

\begin{demos}
	Since $p>\frac{d}{d-\alpha}$, we have $p^{\prime}<\frac{d}{\alpha}$, and thus $\alpha p^{\prime}<d$. Then
	\begin{align*}
		\|S_{\alpha}(\phi)\|_{L^1(\mathbb{R}^d)}
		&= \int_{\mathbb{R}^d} \frac{|\phi(\xi)|}{|\xi|^\alpha}d\xi \\
		&= \left( \int_{|\xi|\leq 1} + \int_{|\xi|\geq 1} \right) \frac{|\phi(\xi)|}{|\xi|^\alpha}d\xi \\
		&\leq \left( \int_{|\xi|\leq 1} \frac{d\xi}{|\xi|^{\alpha p'}} \right)^{1/p'} \left( \int_{\mathbb{R}^d} |\phi(\xi)|^p,d\xi \right)^{1/p} + \|\phi\|_{L^1(\mathbb{R}^d)}\\
		&\leq C \|\phi\|_{1\oplus p},
	\end{align*}
	as desired. \QED
\end{demos}

\begin{coro}
	Let $d\geq 2$, $1\leq q<\frac{d}{\alpha}$, and suppose $\phi\in L^{p\oplus q}(\mathbb{R}^d)$ with $p\in\left(\frac{dq}{d-\alpha q},\infty\right]$. Then
	$\|S_{\alpha}(\phi)\|_{L^q(\mathbb{R}^d)} \leq C \|\phi\|_{p\oplus q}.$
\end{coro}

\begin{demos}
	Applying Proposition \ref{D-inequalities} to $|\phi|^q$, we obtain:
	\begin{align*}
		\|S_{\alpha}(\phi)\|_{L^q(\mathbb{R}^d)}^q
		&= \|S_{\alpha}(|\phi|^q)\|_{L^1(\mathbb{R}^d)} \
		&\leq C\left( \|\phi\|_{L^p(\mathbb{R}^d)}^q + \|\phi\|_{L^q(\mathbb{R}^d)}^q \right) \
		&\leq C \|\phi\|_{p\oplus q}^q.
	\end{align*}
	Taking the $q$-th root gives the desired inequality. \QED
\end{demos}

\begin{coro}\label{Salpha}
	The map $S_{\alpha}: L^{1\oplus p}(\mathbb{R}^d) \to L^1(\mathbb{R}^d)$ is continuous for $p\in \left(\frac{d}{d-\alpha},\infty\right].$
\end{coro}

\begin{coro}
	The map $S_{\alpha}: L^{p\oplus q}(\mathbb{R}^d) \to L^q(\mathbb{R}^d)$ is continuous for $1\leq q<\frac{d}{\alpha}$ and $p\in \left(\frac{dq}{d-\alpha q},\infty\right].$
\end{coro}

Since $\mathcal{D}\subset\bigcap_{1\leq p\leq \infty} L^p(\mathbb{R}^d),$ it is useful to define the subspace
\[
\mathcal{D}{\alpha} = S_{\alpha}(\mathcal{D}) = \left\{ S_{\alpha}(\phi)\;\middle|\; \phi\in\mathcal{D}\right\}
= \left\{ \psi\in L^1(\mathbb{R}^d) \;\middle|\; |\cdot|^\alpha\psi\in\mathcal{D} \right\},
\]
with the topology making $S_{\alpha}: \mathcal{D} \to \mathcal{D}_{\alpha}$ continuous, that is,
$\phi_j \rightarrow \phi $  in  $\mathcal{D}_{\alpha}$ if and only if  $ |\cdot|^\alpha \phi_j \rightarrow |\cdot|^\alpha \phi$ in $\mathcal{D}.$

\paragraph{Convolution Properties in $\mathcal{D}_{\alpha}$:} 
We now study the convolution structure

\begin{teor}\label{alphaconv}
	For $\phi, \psi\in \mathcal{D}_{\alpha}$ and every $n\in\mathbb{N}$,
	\begin{align*}
	p_n\left(|\cdot|^\alpha(\phi\ast\psi)\right) \leq C\left( p_n(|\cdot|^\alpha\phi)p_{n+\lceil d/2\rceil+1}(|\cdot|^\alpha\psi) + p_{n+\lceil d/2\rceil+1}(|\cdot|^\alpha\phi)p_n(|\cdot|^\alpha\psi) \right).
	\end{align*}

\end{teor}

\begin{demos}
	Since  $\phi, \psi \in \mathcal{D}_{\alpha}$ we have that  $\vabs{\cdot}^{\alpha}\phi, \vabs{\cdot}^{\alpha}\psi \in \mathcal{D}$
	and 
	\begin{align*}
		&\vabs{\xi}^{\alpha} \vabs{  (\phi \ast  \psi)(\xi)  }
		\leq 2^{\alpha}\corch{\int_{\reales^ d} \paren{\vabs{\xi -  \eta}^{\alpha}   \vabs{\phi(\xi-\eta)}  } \vabs{\psi(\eta)} d\eta
			+\int_{\reales^ d} \vabs{\phi(\xi-\eta)}  \paren{\vabs{\eta}^{\alpha}\vabs{\psi(\eta)} }d\eta} \\
		& =2^{\alpha}\corch{\int_{\reales^ d} \paren{\vabs{\xi -  \eta}^{\alpha}  \vabs{\phi(\xi-\eta)}   }    \paren{ \vabs{\eta}^{\alpha} \vabs{\psi(\eta)} }\frac{d\eta}{\vabs{\eta}^{\alpha}}
			+\int_{\reales^ d}   \paren{\vabs{\xi-\eta}^{\alpha} \vabs{\phi(\xi-\eta)} } \paren{\vabs{\eta}^{\alpha}\vabs{\psi(\eta)} } \frac{d\eta}{\vabs{\xi-\eta}^{\alpha}} } \\
		& =2^{\alpha}\corch{\int_{\reales^ d} \paren{\vabs{\xi -  \eta}^{\alpha}  \vabs{\phi(\xi-\eta)}   }    \paren{ \vabs{\eta}^{\alpha} \vabs{\psi(\eta)} }\frac{d\eta}{\vabs{\eta}^{\alpha}}
			+\int_{\reales^ d}   \paren{\vabs{\eta}^{\alpha} \vabs{\phi(\eta)} } \paren{\vabs{\xi-\eta}^{\alpha}\vabs{\psi(\xi-\eta)} } \frac{d\eta}{\vabs{\eta}^{\alpha}} } . 
	\end{align*}
	
	Therefore,
	\begin{align*}
		& (1+\vabs{\xi}^2)^n \vabs{\xi}^{\alpha} \vabs{  (\phi \ast  \psi)(\xi)  }
		\leq 2^{n+\alpha} \left[  \int_{\reales^ d} \paren{ (1+\vabs{\xi-\eta}^2)^n  \vabs{\xi -  \eta}^{\alpha}  \vabs{\phi(\xi-\eta)}   }    \paren{ (1+\vabs{\eta}^2)^n \vabs{\eta}^{\alpha} \vabs{\psi(\eta)} }\frac{d\eta}{\vabs{\eta}^{\alpha}}\right. \\
		& \left. +\int_{\reales^ d}   \paren{(1+\vabs{\eta}^2)^n \vabs{\eta}^{\alpha} \vabs{\phi(\eta)} } \paren{(1+\vabs{\xi-\eta}^2)^n  \vabs{\xi-\eta}^{\alpha}\vabs{\psi(\xi-\eta)} } \frac{d\eta}{\vabs{\eta}^{\alpha}} \right] \\
		& \leq 2^{n+\alpha}  \left[   p_{n}(\vabs{\cdot}^{\alpha}  \phi)  \paren{ \int_{\reales^ d} S_{\alpha}\paren{ (1+\vabs{\cdot}^2 )^{n} \vabs{\cdot}^{\alpha}\vabs{\psi} }(\eta)  d\eta}
		+p_{n}(\vabs{\cdot}^{\alpha}  \psi)  \paren{ \int_{\reales^ d} S_{\alpha}\paren{ (1+\vabs{\cdot}^2 )^{n} \vabs{\cdot}^{\alpha} \vabs{\phi} }(\eta)  d\eta} \right] \\
		& =  2^{n+\alpha}  \left[   p_{n}(\vabs{\cdot}^{\alpha}  \phi)  \norma{ S_{\alpha}\paren{ (1+\vabs{\cdot}^2 )^{n} \vabs{\cdot}^{\alpha}\vabs{\psi} }}_{L^{1}(\reales^d)}
		+p_{n}(\vabs{\cdot}^{\alpha}  \psi)  \norma{ S_{\alpha}\paren{ (1+\vabs{\cdot}^2 )^{n} \vabs{\cdot}^{\alpha} \vabs{\phi} }}_{L^{1}(\reales^d)} \right].
	\end{align*}
	
	By  Proposition \ref{D-inequalities} we have 
	\begin{equation*}
		\norma{S_{\alpha}(\phi)}_{L^{1}(\reales^d)  }\leq C \norma{\phi}_{1\oplus p},
	\end{equation*}
	for $\phi \in L^{1\oplus p}(\reales^d)$ and $p\in \left( \frac{d}{d-\alpha},\infty\right].$
	
	However, $\norma{\phi}_{L^{1}(\reales^d) } \leq C p_{\corch{\frac{d}{2}}+1 }(\phi) $ and 
	\begin{align*}
		&\norma{\phi}_{L^{p}(\reales^d) } =\paren{ \int_{\reales^ d} \vabs{ \phi(\xi) }^p d\xi}^{\frac{1}{p}}
		\leq C p_{\corch{\frac{d}{2}}+1 }( \vabs{\phi}^p )^{\frac{1}{p} } 
		\leq C  p_{\corch{\frac{d}{2}}+1 }(\phi),
	\end{align*}
	for $p\in \left( \frac{d}{d-\alpha},\infty\right].$
	
	Therefore,
	\begin{equation}\label{L1D}
		\norma{S_{\alpha}(\phi)}_{L^{1}(\reales^d)  }\leq C p_{\corch{\frac{d}{2}}+1 }(\phi).
	\end{equation}
	
	Consequently,
	
	\begin{align*}
		& (1+\vabs{\xi}^2)^n \vabs{\xi}^{\alpha} \vabs{  (\phi \ast  \psi)(\xi)  }
		\leq \corch{p_{n}(\vabs{\cdot}^{\alpha}  \phi ) p_{n+\corch{\frac{d}{2}}+1}(\vabs{\cdot}^{\alpha}\psi)  +
			p_{n+\corch{\frac{d}{2}}+1}(\vabs{\cdot}^{\alpha}\phi) p_{n}(\vabs{\cdot}^{\alpha}  \psi )  },
	\end{align*}
	for every $\xi \in \reales^d$.

	Thus,
	\begin{equation*}
		p_{n}(\vabs{\cdot}^\alpha(\phi \ast \psi) ) 
		\leq C\paren{p_{n}(\vabs{\cdot}^\alpha \phi) p_{n+\corch{\frac{d}{2}}+1}(\vabs{\cdot}^\alpha \psi)  
			+p_{n+\corch{\frac{d}{2}}+1}(\vabs{\cdot}^\alpha \phi)  p_{n}(\vabs{\cdot}^\alpha \psi)  }.
	\end{equation*}
	\QED
	
\end{demos}

\begin{coro}
	The convolution product $\ast: \mathcal{D}_{\alpha}\times\mathcal{D}_{\alpha}\to\mathcal{D}_{\alpha}$ is continuous.
\end{coro}

\begin{remark}
	Although Theorem \ref{alphaconv} guarantees that $\phi\ast\psi\in \mathcal{D}_{\alpha},$ in some cases the convolution $\phi\ast\psi $even belongs to $\mathcal{D}.$ This important situation is treated next.
\end{remark}

We now prove a refinement:

\begin{teor}
	Let $0\leq \alpha,\beta<d$ with $\alpha+\beta<d.$ Then for $\phi\in\mathcal{D}_{\alpha}, \psi\in\mathcal{D}_{\beta},$ and every $n\in\mathbb{N}$,
	\begin{align*}
			p_n(\phi\ast\psi) \leq C\left( p_n(|\cdot|^\alpha\phi)p_{n+\lceil d/2\rceil+1}(|\cdot|^\beta\psi) + p_{n+\lceil d/2\rceil+1}(|\cdot|^\alpha\phi)p_n(|\cdot|^\beta\psi) \right).
	\end{align*}

\end{teor}

\begin{demos}
	Since $\phi\in \mathcal{D}_{\alpha}$, $\psi\in \mathcal{D}_{\beta}$ we have that $\vabs{\cdot}^{\alpha}$, $\vabs{\cdot}^{\beta}\in \mathcal{D}$.
	For every $n\in \natu$ we have:
	\begin{align*}
		& (1+\vabs{\xi}^2)^n  \vabs{  (\phi \ast  \psi)(\xi)  }
		\leq 2^{n} \left[  \int_{\reales^ d} \paren{ (1+\vabs{\xi-\eta}^2)^n  \vabs{\phi(\xi-\eta)}   }    \paren{ (1+\vabs{\eta}^2)^n \vabs{\psi(\eta)} }d\eta\right] \\
		&    \leq 2^{n} \left[  \int_{\reales^ d}   \paren{ (1+\vabs{\xi-\eta}^2)^n  \vabs{\xi -  \eta}^{\alpha}  \vabs{\phi(\xi-\eta)}   }    \paren{ (1+\vabs{\eta}^2)^n \vabs{\eta}^{\beta} \vabs{\psi(\eta)} } \frac{d\eta}{\vabs{\xi-\eta}^{\alpha}  \vabs{\eta}^{\beta} } \right] \\
		& 	\leq 2^{n} \left[  \int_{\vabs{\xi-\eta}\geq \vabs{\eta} } \paren{ (1+\vabs{\xi-\eta}^2)^n  \vabs{\xi -  \eta}^{\alpha}  \vabs{\phi(\xi-\eta)}   }    \paren{ (1+\vabs{\eta}^2)^n \vabs{\eta}^{\beta} \vabs{\psi(\eta)} }\frac{d\eta}{\vabs{\xi-\eta}^{\alpha}  \vabs{\eta}^{\beta} }\right. \\
		& \left. +\int_{\vabs{\xi-\eta}\geq \vabs{\eta} }   \paren{(1+\vabs{\eta}^2)^n \vabs{\eta}^{\alpha} \vabs{\phi(\eta)} } \paren{(1+\vabs{\xi-\eta}^2)^n  \vabs{\xi-\eta}^{\beta}\vabs{\psi(\xi-\eta)} } \frac{d\eta}{\vabs{\eta}^{\alpha}  \vabs{\xi-\eta}^{\beta} } \right] 
	\end{align*}
	\begin{align*}
		& 	\leq 2^{n} \left[  \int_{\vabs{\xi-\eta}\geq \vabs{\eta} } \paren{ (1+\vabs{\xi-\eta}^2)^n  \vabs{\xi -  \eta}^{\alpha}  \vabs{\phi(\xi-\eta)}   }    \paren{ (1+\vabs{\eta}^2)^n \vabs{\eta}^{\beta} \vabs{\psi(\eta)} }\frac{d\eta}{ \vabs{\eta}^{\alpha+\beta} }\right. \\
		& \left. +\int_{\vabs{\xi-\eta}\geq \vabs{\eta} }   \paren{(1+\vabs{\eta}^2)^n \vabs{\eta}^{\alpha} \vabs{\phi(\eta)} } \paren{(1+\vabs{\xi-\eta}^2)^n  \vabs{\xi-\eta}^{\beta}\vabs{\psi(\xi-\eta)} } \frac{d\eta}{\vabs{\eta}^{\alpha+\beta} } \right] \\
		& \leq 2^{n}  \left[   p_{n}(\vabs{\cdot}^{\alpha}  \phi)  \paren{ \int_{\reales^ d} S_{\alpha+\beta}\paren{ (1+\vabs{\cdot}^2 )^{n} \vabs{\cdot}^{\beta}\vabs{\psi} }(\eta)  d\eta}
		+p_{n}(\vabs{\cdot}^{\beta}  \psi)  \paren{ \int_{\reales^ d} S_{\alpha+\beta}\paren{ (1+\vabs{\cdot}^2 )^{n} \vabs{\cdot}^{\alpha} \vabs{\phi} }(\eta)  d\eta} \right] \\
		& =  2^{n}  \left[   p_{n}(\vabs{\cdot}^{\alpha}  \phi)  \norma{ S_{\alpha}\paren{ (1+\vabs{\cdot}^2 )^{n} \vabs{\cdot}^{\beta}\vabs{\psi} }}_{L^{1}(\reales^d)}
		+p_{n}(\vabs{\cdot}^{\beta}  \psi)  \norma{ S_{\alpha}\paren{ (1+\vabs{\cdot}^2 )^{n} \vabs{\cdot}^{\alpha} \vabs{\phi} }}_{L^{1}(\reales^d)} \right].
	\end{align*}
	In consequence, if we apply Equation \eqref{L1D} we obtain that:
	\begin{equation*}
		(1+\vabs{\xi}^2)^n  \vabs{  (\phi \ast  \psi)(\xi)  }
		\leq C\paren{ p_{n}(\vabs{\cdot}^{\alpha}  \phi) p_{n+\corch{\frac{d}{2}}+1}(\vabs{\cdot}^\beta  \psi) +p_{n+\corch{\frac{d}{2}}+1}(\vabs{\cdot}^{\alpha}  \phi) p_{n}(\vabs{\cdot}^{\beta}  \psi)  },
	\end{equation*}
	for every $\xi \in \reales^d$.
	Thus, 
	\begin{equation*}
		p_{n} (\phi \ast  \psi)
		\leq C\paren{ p_{n}(\vabs{\cdot}^{\alpha}  \phi) p_{n+\corch{\frac{d}{2}}+1}(\vabs{\cdot}^\beta  \psi) +p_{n+\corch{\frac{d}{2}}+1}(\vabs{\cdot}^{\alpha}  \phi) p_{n}(\vabs{\cdot}^{\beta}  \psi)  }.
	\end{equation*} \QED

\end{demos}

\begin{coro}
	There is a continuous bilinear map $\ast: \mathcal{D}_{\alpha}\times \mathcal{D}_{\beta}\to\mathcal{D}$ whenever $0\leq\alpha+\beta<d.$
\end{coro}

\begin{remark}
	The condition $\alpha+\beta<d$ is necessary.
	For example, if we define
	$\phi(\xi) = \frac{e^{-|\xi|^2}}{|\xi|^\alpha}, \quad \psi(\xi) = \frac{e^{-|\xi|^2}}{|\xi|^\beta},$
	then $\phi\in \mathcal{D}_{\alpha}$ and $\psi\in \mathcal{D}_{\beta},$ but
	$(\phi\ast\psi)(0) = \int_{\mathbb{R}^d} \frac{e^{-2|\eta|^2}}{|\eta|^{\alpha+\beta}}\,d\eta,$
	which diverges if $\alpha+\beta\geq d.$
\end{remark}

\begin{coro}
	Let $0\leq \alpha<\frac{d}{2}$. Then for $\phi, \psi\in\mathcal{D}_{\alpha}$ and every $n\in\mathbb{N}$,
	$$p_n(\phi\ast\psi) \leq C\left( p_n(|\cdot|^\alpha\phi)p{n+\lceil d/2\rceil+1}(|\cdot|^\alpha\psi) + p_{n+\lceil d/2\rceil+1}(|\cdot|^\alpha\phi)p_n(|\cdot|^\alpha\psi) \right).$$
\end{coro}

\begin{coro}\label{continuousRconv}
	There is a continuous bilinear map $\ast: \mathcal{D}_{\alpha}\times\mathcal{D}_{\alpha}\to\mathcal{D}$ for $0\leq\alpha<\frac{d}{2}$.
\end{coro}

\subsection{Littlewood-Paley Frequency Atoms}
Let $\varphi \in C_c^\infty(\mathbb{R}^d \setminus \{0\})$ be a non-negative radial function supported in $\mathcal{C} = \{\xi \in \mathbb{R}^d : \frac{3}{4} \le |\xi| \le \frac{8}{3}\}$ such that
\begin{equation}
	\sum_{j \in \mathbb{Z}} \varphi_j(\xi) = 1, \quad \forall \xi \in \mathbb{R}^d \setminus \{0\}, \qquad \text{where } \varphi_j(\xi) := \varphi(2^{-j}\xi).
\end{equation}
For any function $f: \mathbb{R}^d \to \mathbb{C}$ and $j \in \mathbb{Z}$, we define its $j$-th \emph{Littlewood-Paley frequency atom} by
\begin{equation}
	f_j(\xi) := \varphi_j(\xi) f(\xi), \qquad \operatorname{supp}(f_j) \subset \mathcal{C}_j = \left\{ \xi \in \mathbb{R}^d : \tfrac{3}{4} 2^j \le |\xi| \le \tfrac{8}{3} 2^j \right\}.
\end{equation}
We also define the low-frequency cut-off operator $S_k f(\xi) := \sum_{j \le k} f_j(\xi)$.

\begin{lemma}[Atomic Bernstein and Riesz Potential Equivalences]\label{lem:atomic_bernstein}
	For $0 < \alpha < d$, let $\mathcal{S}_\alpha(\phi)(\xi) = |\xi|^{-\alpha} \phi(\xi)$. For every $1 \le p \le r \le \infty$ and every dyadic atom $f_j$ supported in $\mathcal{C}_j$:
	\begin{align}
		\|f_j\|_{L^p(\mathbb{R}^d)} &\le C_d \, 2^{jd\left(\frac{1}{p} - \frac{1}{r}\right)} \|f_j\|_{L^r(\mathbb{R}^d)}, \label{eq:bernstein_holder} \\
		\left(\frac{3}{8}\right)^\alpha 2^{-j\alpha} \|f_j\|_{L^p(\mathbb{R}^d)} &\le \|\mathcal{S}_\alpha(f_j)\|_{L^p(\mathbb{R}^d)} \le \left(\frac{4}{3}\right)^\alpha 2^{-j\alpha} \|f_j\|_{L^p(\mathbb{R}^d)}. \label{eq:bernstein_riesz}
	\end{align}
\end{lemma}
\begin{proof}
	Since $|\mathcal{C}_j| = 2^{jd} |\mathcal{C}|$, \eqref{eq:bernstein_holder} follows from H{\"o}lder's inequality on $\mathcal{C}_j$, and \eqref{eq:bernstein_riesz} follows immediately from $\frac{3}{4} 2^j \le |\xi| \le \frac{8}{3} 2^j$ on $\mathcal{C}_j$.
\end{proof}

\section{Riesz Convolution in Critical Fourier-Besov Spaces}\label{rieszconvolution}

In this section, we introduce a generalization of the classical convolution operation.

\begin{defin}\label{def:riesz_conv}
	Given $0 < \alpha < d$, the Riesz convolution between $f$ and $g$ is defined by
	\begin{equation}
		(f *_\alpha g)(\xi) = \mathcal{S}_\alpha(f * g)(\xi) = \frac{(f * g)(\xi)}{|\xi|^\alpha}, \quad \xi \in \mathbb{R}^d \setminus \{0\}.
	\end{equation}
\end{defin}

\noindent
Observe that $\ast_{\alpha}$ is a commutative operation. Indeed, since
\begin{equation*}
	(f \ast_{\alpha} g)(\xi) = \frac{(f \ast g)(\xi)}{|\xi|^{\alpha}} = \frac{(g \ast f)(\xi)}{|\xi|^{\alpha}} = (g \ast_{\alpha} f)(\xi),
\end{equation*}
for all $\xi \in \mathbb{R}^d \setminus \{0\}$, we deduce that $f \ast_{\alpha} g = g \ast_{\alpha} f$.

However, in general $\ast_{\alpha}$ is distributive but not associative.

\begin{remark}
	A particularly important case arises when $f, g \in \mathcal{D}$. In this setting, since $f \ast g \in \mathcal{D}$, we have $f \ast_{\alpha} g \in \mathcal{D}_{\alpha}$. Thus, the operation $\ast_{\alpha} : \mathcal{D} \times \mathcal{D} \to \mathcal{D}_{\alpha}$ is well-defined. Furthermore, if $0\leq \alpha<\frac{d}{2}$, Corollary~\ref{continuousRconv} ensures that $\ast_{\alpha}: \mathcal{D}_{\alpha} \times \mathcal{D}_{\alpha} \to \mathcal{D}_{\alpha}$ is also continuous.
\end{remark}

When $\alpha = 1$, the Riesz convolution $*_1$ is strictly scale-invariant with respect to the homogeneous Fourier-Besov norms $\mathcal{F}\dot{B}^{s_p}_{p,q}(\mathbb{R}^d)$ with critical regularity index $s_p = \frac{d}{p'} - 1$. We now prove that $*_1$ defines a bounded bilinear operator on $\mathcal{F}\dot{B}^{s_p}_{p,q}(\mathbb{R}^d)$ satisfying a tame paraproduct estimate.

\begin{teor}[Critical Atomic Bilinear Estimate for $*_1$]\label{thm:riesz_critical_besov}
	Let $d \ge 3$, $2 \le p \le \infty$, $1 \le q \le \infty$, and $s_p = \frac{d}{p'} - 1$. For any $f, g \in \mathcal{F}\dot{B}^{s_p}_{p,q}(\mathbb{R}^d)$:
	\begin{equation}\label{eq:tame_riesz_1}
		\|f *_1 g\|_{\mathcal{F}\dot{B}^{s_p}_{p,q}} \le C \left( \|f\|_{\mathcal{F}\dot{B}^{s_p}_{p,\infty}} \|g\|_{\mathcal{F}\dot{B}^{s_p}_{p,q}} + \|f\|_{\mathcal{F}\dot{B}^{s_p}_{p,q}} \|g\|_{\mathcal{F}\dot{B}^{s_p}_{p,\infty}} \right).
	\end{equation}
	Furthermore, for any $\sigma \ge 0$, if $f, g \in \mathcal{F}\dot{B}^{s_p}_{p,\infty}(\mathbb{R}^d) \cap \mathcal{F}\dot{B}^{s_p + \sigma}_{p,q}(\mathbb{R}^d)$, then the high-regularity norm obeys the tame persistence bound:
	\begin{equation}\label{eq:tame_riesz_sigma}
		\|f *_1 g\|_{\mathcal{F}\dot{B}^{s_p + \sigma}_{p,q}} \le C_\sigma \left( \|f\|_{\mathcal{F}\dot{B}^{s_p}_{p,\infty}} \|g\|_{\mathcal{F}\dot{B}^{s_p + \sigma}_{p,q}} + \|f\|_{\mathcal{F}\dot{B}^{s_p + \sigma}_{p,q}} \|g\|_{\mathcal{F}\dot{B}^{s_p}_{p,\infty}} \right).
	\end{equation}
\end{teor}

\begin{demos}
	Decompose $f = \sum_{k \in \mathbb{Z}} f_k$ and $g = \sum_{k' \in \mathbb{Z}} g_{k'}$ into Littlewood-Paley atoms $f_k = \varphi_k f$ and $g_{k'} = \varphi_{k'} g$. By Bony's paraproduct decomposition in frequency space, for each $j \in \mathbb{Z}$:
	\begin{equation}
		\varphi_j(f * g) = \sum_{|k-j| \le 4} \varphi_j(f_k * S_{k-3} g) + \sum_{|k'-j| \le 4} \varphi_j(S_{k'-3} f * g_{k'}) + \sum_{k \ge j-3} \sum_{|k-k'| \le 2} \varphi_j(f_k * g_{k'}).
	\end{equation}
	By Lemma \ref{lem:atomic_bernstein}, for any $\sigma \ge 0$,
	\begin{equation}
		2^{j(s_p + \sigma)} \|\varphi_j(f *_1 g)\|_{L^p} \le \frac{4}{3} 2^{j(s_p + \sigma - 1)} \|\varphi_j(f * g)\|_{L^p}.
	\end{equation}
	We bound each of the three paraproduct terms:
	
	\noindent\textbf{Step 1 (High-Low Interaction: $|k-j| \le 4$, $k' \le k-3$).}
	By Young's inequality and Lemma \ref{lem:atomic_bernstein} (noting $\frac{d}{p'} = s_p + 1$):
	\begin{equation}
		\|f_k * S_{k-3} g\|_{L^p} \le \|f_k\|_{L^p} \sum_{k' \le k-3} \|g_{k'}\|_{L^1} \le C \|f_k\|_{L^p} \sum_{k' \le k-3} 2^{k'} \left(2^{k' s_p} \|g_{k'}\|_{L^p}\right) \le C 2^j \|f_k\|_{L^p} \|g\|_{\mathcal{F}\dot{B}^{s_p}_{p,\infty}}.
	\end{equation}
	Multiplying by $2^{j(s_p + \sigma - 1)}$ and using $2^{j(s_p + \sigma)} \approx 2^{k(s_p + \sigma)}$ yields:
	\begin{equation}
		2^{j(s_p + \sigma - 1)} \sum_{|k-j| \le 4} \|\varphi_j(f_k * S_{k-3} g)\|_{L^p} \le C \left(\sum_{|k-j| \le 4} 2^{k(s_p + \sigma)} \|f_k\|_{L^p}\right) \|g\|_{\mathcal{F}\dot{B}^{s_p}_{p,\infty}}.
	\end{equation}
	
	\noindent\textbf{Step 2 (Low-High Interaction: $|k'-j| \le 4$, $k \le k'-3$).}
	By commutativity of convolution, the identical argument gives the bound
	$C \|f\|_{\mathcal{F}\dot{B}^{s_p}_{p,\infty}} \left(\sum_{|k'-j| \le 4} 2^{k'(s_p + \sigma)} \|g_{k'}\|_{L^p}\right)$.
	
	\noindent\textbf{Step 3 (High-High Interaction: $k \ge j-3$, $|k-k'| \le 2$).}
	Since $p \ge 2$, we have $p' \le p$. By Lemma \ref{lem:atomic_bernstein} and Young's inequality ($L^p * L^{p'} \to L^\infty$):
	\begin{equation}
		\|\varphi_j(f_k * g_{k'})\|_{L^p} \le C 2^{j\frac{d}{p}} \|f_k * g_{k'}\|_{L^\infty} \le C 2^{j\frac{d}{p}} \|f_k\|_{L^p} \|g_{k'}\|_{L^{p'}} \le C 2^{j\frac{d}{p}} 2^{kd\left(1 - \frac{2}{p}\right)} \|f_k\|_{L^p} \|g_{k'}\|_{L^p}.
	\end{equation}
	Since $s_p - 1 + \frac{d}{p} = d - 2$ and $d\paren{1 - \frac{2}{p}} = 2s_p - (d-2)$, multiplying by $2^{j(s_p + \sigma - 1)}$ yields:
	\begin{equation}
		2^{j(s_p + \sigma - 1)} \|\varphi_j(f_k * g_{k'})\|_{L^p} \le C 2^{-(k-j)(d - 2 + \sigma)} \left(2^{k(s_p + \sigma)} \|f_k\|_{L^p}\right) \left(2^{k' s_p} \|g_{k'}\|_{L^p}\right).
	\end{equation}
	Since $d \ge 3$ and $\sigma \ge 0$, we have $d - 2 + \sigma \ge 1 > 0$. Applying Young's inequality for discrete convolutions in $\ell^q(\mathbb{Z})$ gives \eqref{eq:tame_riesz_1} (at $\sigma = 0$) and \eqref{eq:tame_riesz_sigma}. \QED
\end{demos}

\begin{coro}[Tame Monomial Bounds for $*_1$]\label{cor:monomial_tame}
	Let $d \ge 3$, $2 \le p \le \infty$, $1 \le q \le \infty$, and $s_p = \frac{d}{p'} - 1$. For any $f_1, \dots, f_k \in \mathcal{F}\dot{B}^{s_p}_{p,\infty}(\mathbb{R}^d) \cap \mathcal{F}\dot{B}^{s_p + \sigma}_{p,q}(\mathbb{R}^d)$ and any nonassociative monomial $w \in \mathcal{M}(x_1, \dots, x_k)$:
	\begin{align}
		\|w^{*_1}(f_1, \dots, f_k)\|_{\mathcal{F}\dot{B}^{s_p}_{p,\infty}} &\le C^{k-1} \prod_{j=1}^k \|f_j\|_{\mathcal{F}\dot{B}^{s_p}_{p,\infty}}, \label{eq:monomial_weak} \\
		\|w^{*_1}(f_1, \dots, f_k)\|_{\mathcal{F}\dot{B}^{s_p + \sigma}_{p,q}} &\le C_\sigma^{k-1} \sum_{r=1}^k \|f_r\|_{\mathcal{F}\dot{B}^{s_p + \sigma}_{p,q}} \prod_{\substack{j=1 \\ j \ne r}}^k \|f_j\|_{\mathcal{F}\dot{B}^{s_p}_{p,\infty}}. \label{eq:monomial_tame_sigma} 
	\end{align}
\end{coro}
\begin{demos}
	Follows by induction on $k$ from Theorem \ref{thm:riesz_critical_besov}. \QED
\end{demos}

\section{Spaces of Functions Dominated by Fourier Caloric Functions}\label{spacesofcaloric}

In this section, we introduce function spaces involving both space and time variables, motivated by the study of evolution equations such as the Navier-Stokes equations.

\begin{defin}
	Let $\mathcal{E}$ be a space of rapidly decreasing functions such that $\mathcal{E}^{+}$ is closed under pointwise addition, convolution, and taking maxima.
	
	Given $p,q,n\in \mathbb{N}$ and $\alpha>0$, we define $\mathcal{C}_{\alpha}(\mathcal{E})_{n}^{p\times q}$ as the complex vector space generated by functions $f:\mathbb{R}_{+}^{d+1}\to \mathbb{C}^{p\times q}$ satisfying:
	\begin{equation}\label{caloricdef}
		|f(\xi,t)|\leq t^{n}e^{-\lambda t |\xi|^\alpha} f^{0}(\xi),
	\end{equation}
	for some $\lambda>0$ and $f^{0}\in \mathcal{E}^{+}$. We denote $\mathcal{C}_{\alpha}(\mathcal{E})^{p\times q}=\bigoplus_{n=0}^{\infty} \mathcal{C}_{\alpha}(\mathcal{E})_{n}^{p\times q}$. The parameter $\lambda$ is called the exponent of $f$, and we denote it by $\exp(f)$.
\end{defin}

When no confusion arises, we simply write $\mathcal{C}_{\alpha}^{p\times q}$.

\begin{remark}\label{caloric}
	\begin{itemize}
		\item If $f_1,\dots,f_k\in \mathcal{C}_{\alpha}^{p\times q}$ with respective exponents $\lambda_j=\exp(f_j)$ for $1\leq j\leq k$, we can choose a common exponent $\lambda = \min\{\lambda_j\}$.
		\item Since $\mathcal{E}^{+}$ is closed under taking maxima, we may take $f^{0} = \max\{f_j^0: 1\leq j\leq k\}$.
		\item Consequently, each $f\in \mathcal{C}_{\alpha}^{p\times q}$ satisfies
		\[
		|f(\xi,t)|\leq p(t) e^{-\lambda t |\xi|^\alpha} f^{0}(\xi),
		\]
		for a polynomial $p\in \mathbb{C}[t]$ with $p([0,\infty))\subset [0,\infty)$.
	\end{itemize}
\end{remark}

We obtain the following basic property.

\begin{prop}
	For each $p,q\in \mathbb{N}$, the space $\mathcal{C}_{\alpha}^{p\times q}$ is a graded $\mathcal{E}$-module and a $\mathbb{C}[t]$-module under pointwise operations.
\end{prop}

Now, we introduce a natural pointwise product.

\begin{defin}
	Given $f\in \mathcal{C}_{\alpha}^{p\times q}$ and $g\in \mathcal{C}_{\alpha}^{r\times q}$, define their \emph{pointwise product} as
	\[
	(f\cdot g)(\xi,t) = f(\xi,t) g(\xi,t)^*, \quad (\xi,t)\in \mathbb{R}_{+}^{d+1}.
	\]
\end{defin}

\begin{prop}\label{pointwiseproduct}
	If $\mathcal{E}^{+}$ is closed under pointwise products, then the operation
	\[
	\cdot:\mathcal{C}_{\alpha}^{p\times q} \times \mathcal{C}_{\alpha}^{r\times q} \to \mathcal{C}_{\alpha}^{p\times r}
	\]
	is bilinear. Moreover, $\mathcal{C}_{\alpha}^{p\times p}$ is a graded $\mathcal{E}$-algebra under $\cdot$.
\end{prop}

\begin{demos}
	Given $f,g$ as above and assuming they are generators, there exist $\lambda_1, \lambda_2 >0$ and $f^0,g^0\in \mathcal{E}^{+}$ such that
	\[
	|f(\xi,t)|\leq t^{n_1} e^{-\lambda_1 t |\xi|^\alpha} f^0(\xi), \quad |g(\xi,t)|\leq t^{n_2} e^{-\lambda_2 t |\xi|^\alpha} g^0(\xi).
	\]
	Thus,
	\[
	|(f\cdot g)(\xi,t)| \leq t^{n_1+n_2} e^{-(\lambda_1+\lambda_2) t |\xi|^\alpha} f^0(\xi) g^0(\xi),
	\]
	and since $f^0g^0\in \mathcal{E}^{+}$, the result follows. \QED
\end{demos}

\begin{coro}
	If $\mathcal{E}^{+}$ is closed under pointwise products, then $\mathcal{C}_{\alpha}^{p\times p}$ is a $\mathbb{C}$-algebra.
\end{coro}

We now introduce an operation akin to a convolution.

\begin{defin}
	Given $f\in \mathcal{C}_{\alpha}^{p\times q}$ and $g\in \mathcal{C}_{\alpha}^{r\times q}$, define their \emph{tensor convolution} by
	\[
	(f\ast g)(\xi,t) = \int_{\mathbb{R}^d} f(\xi-\eta,t)\cdot g(\eta,t) \, d\eta.
	\]
\end{defin}

\begin{teor}\label{tensorconvolution}
	The operation
	\[
	\ast: \mathcal{C}_{\alpha}^{p\times q} \times \mathcal{C}_{\alpha}^{r\times q} \to \mathcal{C}_{\alpha}^{p\times r}
	\]
	is bilinear. Moreover, $\mathcal{C}_{\alpha}^{p\times p}$ is a graded $\mathcal{E}$-algebra under $\ast$.
\end{teor}

\begin{demos}
	The bilinearity is immediate. Assume $f,g$ are generators:
	\[
	|f(\xi,t)|\leq t^{n_1} e^{-\lambda t |\xi|^\alpha} f^0(\xi), \quad |g(\xi,t)|\leq t^{n_2} e^{-\lambda t |\xi|^\alpha} g^0(\xi).
	\]
	By Proposition~\ref{magic}, there exists $r_{\alpha}\leq 1$ such that
	\[
	|f(\xi-\eta,t)\cdot g(\eta,t)| \leq t^{n_1+n_2} e^{-r_{\alpha}\lambda t |\xi|^\alpha} f^0(\xi-\eta) g^0(\eta).
	\]
	Integrating over $\eta$, and using closure properties of $\mathcal{E}^{+}$, we conclude $f\ast g\in \mathcal{C}_{\alpha}^{p\times r}$. \QED
\end{demos}

For $p,q\in \mathbb{N}$ and $\lambda>0$, we define
\[
\mathcal{C}_{\alpha}^{p\times q}(\lambda) = \{f\in \mathcal{C}_{\alpha}^{p\times q} \mid \exp(f)\in (0,\lambda)\}.
\]
In particular, when $\alpha=2$, we simply write $\mathcal{C}^{p\times q}(\lambda)$.

We now introduce an important operation involving the Gaussian kernel, motivated by applications to the Navier-Stokes equations.

\begin{defin}\label{gaussian_conv}
	Let $K:\mathbb{R}^d\to \mathbb{C}^{d\times d}$ satisfy
	\[
	\sup_{\xi\in\mathbb{R}^d} \|K(\xi)\|_{\mathcal{L}(\mathbb{C}^d)}\leq 1.
	\]
	For $f,g\in \mathcal{C}_{d,1}(4\pi^2\nu)$, define the \emph{Gaussian convolution product} by
	\[
	(f\odot g)(\xi,t) = 2\pi i\, K(\xi) \left( \int_0^t e^{-4\pi^2 \nu (t-s)|\xi|^2} (f\ast g)(\xi,s) \, ds \right) \xi.
	\]
\end{defin}

We now verify the main algebraic properties of this operation.

\begin{teor}\label{gaussianproduct}
	Assume that $\mathcal{E}^{+}$ is closed under Riesz convolution $\ast_1$, that is, $f^{0}\ast_1 g^{0} \in \mathcal{E}^{+}$ for all $f^{0},g^{0}\in \mathcal{E}^{+}$, or that $\mathcal{E}^{+}$ is closed under convolution after multiplication by $|\cdot|$, meaning $\ |\cdot| (f^{0}\ast g^{0}) \in \mathcal{E}^{+}$. Then the map
	\[
	\odot: \mathcal{C}^{d\times 1}(4\pi^2 \nu) \times \mathcal{C}^{d\times 1}(4\pi^2 \nu) \longrightarrow \mathcal{C}^{d\times 1}(4\pi^2 \nu)
	\]
	is bilinear.
\end{teor}

\begin{demos}
	Bilinearity is straightforward by definition. To show that $f\odot g\in \mathcal{C}^{d\times 1}(4\pi^2 \nu)$ for $f,g$ generators, we assume
	\[
	|f(\xi,t)|\leq t^{n_1} e^{-\lambda t |\xi|^2} f^{0}(\xi),\quad |g(\xi,t)|\leq t^{n_2} e^{-\lambda t |\xi|^2} g^{0}(\xi),
	\]
	with $\lambda\leq 4\pi^2 \nu$.
	
	By Theorem~\ref{tensorconvolution},
	\[
	|(f\ast g)(\xi,t)|\leq t^{n_1+n_2} e^{-\frac{\lambda}{2} t |\xi|^2} (f^0\ast g^0)(\xi).
	\]
	
	Thus,
	\begin{align*}
		|(f\odot g)(\xi,t)| &\leq 2\pi \left( \int_0^t e^{-4\pi^2\nu (t-s)|\xi|^2} s^{n_1+n_2} e^{-\frac{\lambda}{2}s|\xi|^2} (f^0\ast g^0)(\xi) \, ds \right) |\xi| \\
		&= 2\pi t^{n_1+n_2} e^{-4\pi^2\nu t|\xi|^2} \left( \int_0^t e^{(4\pi^2\nu - \frac{\lambda}{2}) s |\xi|^2} \, ds \right) (f^0\ast g^0)(\xi) |\xi|.
	\end{align*}
	
	Carrying out the integral and simplifying, we distinguish two cases:
	
	\begin{itemize}
		\item[(i)] If $\mathcal{E}^{+}$ is closed under Riesz convolution $\ast_1$, then
		\[
		|(f\odot g)(\xi,t)| \leq \frac{t^{n_1+n_2}}{\pi\nu} e^{-\frac{\lambda}{2}t|\xi|^2} (f^0\ast_1 g^0)(\xi),
		\]
		where $f^0\ast_1 g^0\in \mathcal{E}^{+}$.
		
		\item[(ii)] If $\mathcal{E}^{+}$ is closed under convolution after multiplication by $|\cdot|$, then we obtain a similar estimate with an additional factor of $t$.
	\end{itemize}
	Thus, in both cases $f\odot g\in \mathcal{C}^{d\times 1}(4\pi^2\nu)$. \QED
\end{demos}

\begin{coro}\label{gaussianproductbound}
	If $\mathcal{E}^{+}$ is closed under Riesz convolution $\ast_1$, then for $f,g\in \mathcal{C}^{d\times 1}(4\pi^2 \nu)$,
	\[
	|(f\odot g)(\xi,t)| \leq \frac{t^{n_1+n_2}}{\pi\nu} e^{-\frac{\lambda}{2}t|\xi|^2} (f^0\ast_1 g^0)(\xi).
	\]
\end{coro}

\begin{demos}
	Since $\lambda\leq 4\pi^2\nu$, we have $4\pi^2\nu - \frac{\lambda}{2} \geq 2\pi^2\nu$, and thus
	\[
	\frac{2\pi}{4\pi^2\nu - \frac{\lambda}{2}} \leq \frac{2\pi}{2\pi^2\nu} = \frac{1}{\pi\nu}.
	\]
	\QED
\end{demos}

We now introduce spaces that are well-adapted to the behavior of time derivatives, which are essential in the study of evolution equations like Navier-Stokes.

\begin{defin}
	Let $\mathcal{E}$ be a space of rapidly decreasing functions, closed under pointwise addition, convolution, and taking maxima.
	
	Let $d\geq 3$. We define $\mathcal{V}(\mathcal{E})$ as the complex vector space generated by functions $f: \mathbb{R}_+^{d+1} \to \mathbb{C}^{d\times 1}$ such that:
	\begin{itemize}
		\item $f\in \mathcal{C}(\mathcal{E})_0$,
		\item For almost every $\xi\in\mathbb{R}^d$, the map $t\mapsto f(\xi,t)$ belongs to $C^{\infty}([0,\infty),\mathbb{C}^d)$,
		\item For every $m,n\in\mathbb{N}$, there exist constants $a,b>0$ and elements $f_{m,n}^0\in \mathcal{E}^{+}$ such that
		\[
		\left|(\lambda^{1/2}|\xi|)^m \frac{\partial^n f}{\partial t^n}(\xi,t)\right| \leq (\lambda a)^{m/2} (\lambda b)^n e^{-\lambda t |\xi|^2} f_{m,n}^0(\xi)
		\]
		for all $(\xi,t)\in\mathbb{R}_+^{d+1}$ and some $\lambda>0$.
	\end{itemize}
\end{defin}

We denote by $\mathcal{V}$ the space $\mathcal{V}(\mathcal{E})$ when the context is clear.

\begin{remark}
	Note that $\lambda$ is a uniform exponent for all time derivatives and weight powers in the previous definition.
\end{remark}

We now analyze the behavior of the tensor convolution restricted to $\mathcal{V}$.

\begin{teor}\label{tensorconvolutionsmooth}
	The tensor convolution defines a bilinear map
	\[
	\ast: \mathcal{V} \times \mathcal{V} \to \mathcal{V}.
	\]
	In particular, $\mathcal{V}$ is a graded $\mathcal{E}$-algebra under convolution.
\end{teor}

\begin{demos}
	Let $f,g\in \mathcal{V}$. There exist constants $a,b,c,d>0$, $\beta>0$, and elements $f_{m_1,n_1}^0, g_{m_2,n_2}^0\in\mathcal{E}^{+}$ such that
	\begin{align}
		\left|(\beta^{1/2}|\xi|)^{m_1} \frac{\partial^{n_1}f}{\partial t^{n_1}}(\xi,t)\right| &\leq (\beta a)^{m_1/2} (\beta b)^{n_1} e^{-\beta t|\xi|^2} f_{m_1,n_1}^0(\xi),\label{f_estimate}\\
		\left|(\beta^{1/2}|\xi|)^{m_2} \frac{\partial^{n_2}g}{\partial t^{n_2}}(\xi,t)\right| &\leq (\beta c)^{m_2/2} (\beta d)^{n_2} e^{-\beta t|\xi|^2} g_{m_2,n_2}^0(\xi).\label{g_estimate}
	\end{align}
	
	Applying the binomial formula and the Leibniz rule gives
	\begin{align*}
		\left|(\beta^{1/2}|\xi|)^m \frac{\partial^n (f\ast g)}{\partial t^n}(\xi,t)\right|
		&\leq \sum_{r=0}^n\sum_{j=0}^m \binom{n}{r}\binom{m}{j} \left(
		(\beta^{1/2}|\cdot|)^j \left|\frac{\partial^r f}{\partial t^r}\right|
		\ast
		(\beta^{1/2}|\cdot|)^{m-j} \left|\frac{\partial^{n-r} g}{\partial t^{n-r}}\right|
		\right)(\xi,t).
	\end{align*}
	
	Using the estimates \eqref{f_estimate} and \eqref{g_estimate} together with the properties of $\mathcal{E}$ and the behavior under convolution (Theorem~\ref{tensorconvolution}), we obtain
	\[
	\left|(\beta^{1/2}|\xi|)^m \frac{\partial^n (f\ast g)}{\partial t^n}(\xi,t)\right|
	\leq (\beta^{1/2}(a^{1/2}+c^{1/2}))^m (\beta(b+d))^n e^{-\beta t|\xi|^2} (f\ast g)_{m,n}^0(\xi),
	\]
	where
	\[
	(f\ast g)_{m,n}^0(\xi) = \max_{0\leq j\leq m, 0\leq r\leq n} (f_{j,r}^0\ast g_{m-j,n-r}^0)(\xi)\in\mathcal{E}^{+}.
	\]
	Hence, $f\ast g\in\mathcal{V}$. \QED
\end{demos}

\begin{remark}\label{lambda}
	If $a,c\geq 1$, the above inequality simplifies to
	\[
	\left|(\beta^{1/2}|\xi|)^m \frac{\partial^n (f\ast g)}{\partial t^n}(\xi,t)\right|
	\leq (\beta^{1/2}(a+c))^m (\beta(b+d))^n e^{-\beta t|\xi|^2} (f\ast g)_{m,n}^0(\xi).
	\]
\end{remark}

\begin{defin}
	Given $\beta>0$, define
	\[
	\mathcal{V}(\beta) = \llave{f\in \mathcal{V} \mid \exp\left(|\cdot|^m \frac{\partial^n f}{\partial t^n}\right) \leq \beta,\quad \forall\, m,n\in\mathbb{N}}.
	\]
	In particular, we define $\lambda = 4\pi^2\nu$.
\end{defin}

We conclude this subsection with an important structural formula for time derivatives of the $\odot$ product.

\begin{lemma}\label{derivativedot}
	For every $f,g\in \mathcal{V}(\lambda)$ and $n\in\mathbb{N}$,
	\[
	\frac{\partial^n (f\odot g)}{\partial t^n}(\xi,t)
	= (-\lambda |\xi|^2)^n (f\odot g)(\xi,t)
	+ 2\pi i K(\xi) \left( \sum_{j=0}^{n-1} (-\lambda |\xi|^2)^j \frac{\partial^{n-1-j}(f\ast_1 g)}{\partial t^{n-1-j}}(\xi,t)|\xi| \right)\xi,
	\]
	for every $(\xi,t)\in\mathbb{R}_+^{d+1}$.
\end{lemma}

\begin{demos}
	Apply Lemma~\ref{powerlinear} to the operator $L=\partial_t$, the unknown $x(\xi,t)=(f\odot g)(\xi,t)$, and $y(\xi,t)=2\pi i K(\xi)(f\ast_1 g)(\xi,t)|\xi|\xi$, with $\alpha(\xi) = -\lambda|\xi|^2$.
	\QED
\end{demos}

We continue using the notation established so far to present the following result.

\begin{teor}
	Let $f,g \in \mathcal{V}(\lambda)$ and $m,n\in \mathbb{N}$. Then
	\begin{align*}
		&\left|(\lambda^{1/2} |\xi|)^m \frac{\partial^n (f\odot g)}{\partial t^n}(\xi,t)\right| \\
		&\leq \frac{(\lambda(a^{1/2}+c^{1/2})^2)^{\frac{m}{2}+1} (\lambda(b+d))^{n-1}}{\pi \nu}
		\left( \frac{(a^{1/2}+c^{1/2})^{2n} - (b+d)^n}{(b+d)^{n-1}\left((a^{1/2}+c^{1/2})^2-(b+d)\right)} \right)
		e^{-\frac{\beta t}{2} |\xi|^2} (f\odot g)_{m,n}^0(\xi),
	\end{align*}
	for some $(f\odot g)_{m,n}^0 \in \mathcal{E}^+$.
\end{teor}

\begin{demos}
	Let $f,g \in \mathcal{V}(\lambda)$. By Remark~\ref{lambda}, there exist constants $a,b,c,d>0$ and elements $f_{m_1,n_1}^0,g_{m_2,n_2}^0\in \mathcal{E}^+$ such that
	\begin{align*}
		\left|(\lambda^{1/2} |\xi|)^{m_1} \frac{\partial^{n_1} f}{\partial t^{n_1}}(\xi,t)\right| &\leq (\lambda a)^{m_1/2} (\lambda b)^{n_1} e^{-\beta t |\xi|^2} f_{m_1,n_1}^0(\xi),\\
		\left|(\lambda^{1/2} |\xi|)^{m_2} \frac{\partial^{n_2} g}{\partial t^{n_2}}(\xi,t)\right| &\leq (\lambda c)^{m_2/2} (\lambda d)^{n_2} e^{-\beta t |\xi|^2} g_{m_2,n_2}^0(\xi).
	\end{align*}
	
	By Theorem~\ref{tensorconvolutionsmooth} and Corollary~\ref{gaussianproductbound}, we obtain:
	\[
	\left|(\lambda^{1/2}|\xi|)^m \frac{\partial^n (f\ast_1 g)}{\partial t^n}(\xi,t)\right|
	\leq (\lambda^{1/2}(a^{1/2}+c^{1/2}))^m (\lambda(b+d))^n e^{-\frac{\beta t}{2}|\xi|^2} (f\ast_1 g)_{m,n}^0(\xi).
	\]
	Moreover, for $n=0$,
	\[
	\left|(\lambda^{1/2}|\xi|)^m (f\odot g)(\xi,t)\right| \leq \frac{(\lambda^{1/2}(a^{1/2}+c^{1/2}))^m}{\pi\nu} e^{-\frac{\beta t}{2}|\xi|^2} (f\ast_1 g)_{m,0}^0(\xi).
	\]
	
	Applying Lemma~\ref{derivativedot}, using the bounds for derivatives of $f\ast_1 g$, and summing carefully, we find:
	\[
	\left|(\lambda^{1/2}|\xi|)^m \frac{\partial^n(f\odot g)}{\partial t^n}(\xi,t)\right|
	\leq \text{(major term)} + \text{(summed series)},
	\]
	where the series is a geometric sum evaluated explicitly:
	\[
	\sum_{j=0}^{n-1}\left(\frac{(a^{1/2}+c^{1/2})^2}{b+d}\right)^j = \frac{(a^{1/2}+c^{1/2})^{2n}-(b+d)^n}{(b+d)^{n-1}\left((a^{1/2}+c^{1/2})^2-(b+d)\right)}.
	\]
	The result follows.
	\QED
\end{demos}

\begin{coro}\label{coro1}
	If $(a^{1/2}+c^{1/2})^2 \leq 2(b+d)$, then
	\[
	\left|(\lambda^{1/2} |\xi|)^m \frac{\partial^n (f\odot g)}{\partial t^n}(\xi,t)\right|
	\leq \frac{2^{n+1} (\lambda(a^{1/2}+c^{1/2})^2)^{m/2} (\lambda(b+d))^n}{\pi \nu}
	e^{-\frac{\beta t}{2}|\xi|^2} (f\odot g)_{m,n}^0(\xi),
	\]
	for all $(\xi,t)\in \mathbb{R}_+^{d+1}$.
\end{coro}

\begin{remark}
	The condition of Corollary~\ref{coro1} is satisfied in particular if $a\leq b$ and $c\leq d$.
\end{remark}

\begin{coro}\label{coro2}
	If $a\geq b$ and $c\geq d$, then
	\[
	\left|(\lambda^{1/2} |\xi|)^m \frac{\partial^n (f\odot g)}{\partial t^n}(\xi,t)\right|
	\leq \frac{(\lambda(a^{1/2}+c^{1/2})^2)^{m/2+n}}{2a^{1/2}c^{1/2}\pi\nu}
	e^{-\frac{\beta t}{2}|\xi|^2} (f\odot g)_{m,n}^0(\xi),
	\]
	for all $(\xi,t)\in \mathbb{R}_+^{d+1}$.
\end{coro}

\begin{demos}
	The proof follows by observing that under $a\geq b$ and $c\geq d$, we have
	\[
	(a^{1/2}+c^{1/2})^2-(b+d)\geq 2a^{1/2}c^{1/2},
	\]
	thus improving the denominator.
	\QED
\end{demos}

\begin{remark}\label{lambdasimplified}
	By Remark~\ref{lambda}, if $f\in\mathcal{V}(\lambda)$, then
	\[
	\left|(\lambda^{1/2} |\xi|)^m \frac{\partial^n f}{\partial t^n}(\xi,t)\right| \leq (\lambda a)^{m/2}(\lambda b)^n e^{-\beta t|\xi|^2} f_{m,n}^0(\xi).
	\]
	Choosing $a=b$ simplifies the estimates:
	\[
	\left|(\lambda^{1/2} |\xi|)^m \frac{\partial^n f}{\partial t^n}(\xi,t)\right| \leq (\lambda a)^{(m+n)/2} e^{-\beta t|\xi|^2} f_{m,n}^0(\xi).
	\]
\end{remark}

\begin{coro}\label{corogold}
	Let $f,g\in\mathcal{V}(\lambda)$ satisfy
		\begin{equation}
			\vabs{(\lambda^{\frac{1}{2}}  \vabs{\xi})^{m_{1} }\frac{\partial^{n_{1}} f }{\partial t^{n_{1}}}(\xi,t)}  \leq (\lambda a)^{\frac{m_{1}}{2}+n_{1}}  e^{-\beta t \vabs{\xi}^2}f_{m_{1},n_{1}}^ {0}(\xi),
			\end{equation}
		and 
		\begin{equation}
			\vabs{(\lambda^{\frac{1}{2}}  \vabs{\xi})^{m_{2} }\frac{\partial^{n_{2}} g }{\partial t^{n_{2}} }(\xi,t)}  \leq (\lambda c)^{\frac{m_{2}}{2}+n_{2} } e^{-\beta t \vabs{\xi}^2}g_{m_{2},n_{2}}^ {0}(\xi),
			\end{equation}
		for every $(\xi,t)\in \reales_{+}^{d+1}$, for some $0< \beta \leq  \lambda$, $f_{m_{1},n_{1}}^{0}, g_{m_{2},n_{2}}^{0}\in \mathcal{E}^{+} $, $a,c\in \reales_{>0}$. 
Then
	\begin{equation}\label{gold}
			\left|(\lambda^{1/2}|\xi|)^m \frac{\partial^n(f\odot g)}{\partial t^n}(\xi,t)\right|
		\leq \frac{(\lambda(a^{1/2}+c^{1/2})^2)^{m/2+n}}{2a^{1/2}c^{1/2}\pi\nu}
		e^{-\frac{\beta t}{2}|\xi|^2} (f\odot g)_{m,n}^0(\xi),
	\end{equation}
	for all $(\xi,t)\in\mathbb{R}_+^{d+1}$.
\end{coro}

\begin{demos}
	This follows immediately from Corollary~\ref{coro2} since $a=b$ and $c=d$.
	\QED
\end{demos}

\begin{coro}\label{gaussianproductsmooth}
	Assume that $\mathcal{E}^+$ is closed under Riesz convolution $\ast_1$. Then the product
	\[
	\odot: \mathcal{V}(\lambda)\times\mathcal{V}(\lambda) \to \mathcal{V}(\lambda)
	\]
	is a bilinear operator.
\end{coro}

\begin{demos}
	This is a direct consequence of inequality~\eqref{gold} and the properties of $\mathcal{E}$.
	\QED
\end{demos}

Finally, we highlight a fundamental case for constructing solutions to the Navier-Stokes equations.

\begin{coro}\label{Vasquez}
	Let $f,g\in\mathcal{V}$ with $a,c\in\mathbb{N}_+$. Then
	\[
	\left|(\lambda^{1/2}|\xi|)^m \frac{\partial^n(f\odot g)}{\partial t^n}(\xi,t)\right|
	\leq \frac{(\lambda(a^{1/2}+c^{1/2})^2)^{m/2+n}}{2\pi\nu}
	e^{-\frac{\beta t}{2}|\xi|^2} (f\odot g)_{m,n}^0(\xi),
	\]
	for all $(\xi,t)\in\mathbb{R}_+^{d+1}$.
\end{coro}

\begin{remark}\label{increasingsequence}
	Observe that
	\[
	(f\odot g)_{m,n}^0(\xi) = \max_{\substack{0\leq r_1+r_2\leq n-1\\ l_1+l_2+2(r_1+r_2)=m+2n}} (f_{l_1,r_1}^0 \ast_1 g_{l_2,r_2}^0)(\xi).
	\]
	Thus, if $f_{m,n}^0\leq f_{m',n'}^0$ and $g_{m,n}^0\leq g_{m',n'}^0$ whenever $m\leq m'$ and $n\leq n'$, then
	\[
	(f\odot g)_{m,n}^0\leq (f\odot g)_{m',n'}^0.
	\]
\end{remark}

\subsection{Atomic Theory of Fourier-Besov Caloric Spaces  $\mathcal{V}(\lambda)$}

We now combine the Fourier-Caloric spaces $\mathcal{V}(\lambda)$ with the Littlewood-Paley decomposition, extending the critical bilinear bounds to all $1 \le p \le \infty$.

\begin{defin}[Fourier-Besov Caloric Space $\mathcal{V}(\lambda)$]\label{def:space_V}
	Let $\mathcal{E} = \mathcal{E}_{\mathcal{B}}$ be the space of rapidly decreasing functions associated to $\mathcal{B} = \mathcal{F}\dot{B}^{s_p}_{p,1}(\mathbb{R}^d) \cap L^{1 \oplus 2}(\mathbb{R}^d)$. Let $\mathcal{V}$  be the space associated to $\mathcal{E}$ and 
	\[
	\mathcal{V}(\lambda) = \llave{f\in \mathcal{V} \mid \exp\left(|\cdot|^m \frac{\partial^n f}{\partial t^n}\right) \leq \lambda,\quad \forall\, m,n\in\mathbb{N}}.
	\]
\end{defin}

\begin{teor}[Atomic Caloric Bilinear Bound for $\odot$ Across $1 \le p \le \infty$]\label{thm:atomic_caloric_odot}
	Let $d \ge 3$, $1 \le p \le \infty$, and $s_p = \frac{d}{p'} - 1$. Suppose $f, g \in \mathcal{V}(\lambda)$ satisfy
	\begin{equation}
		|f(\xi,t)| \le e^{-\beta t |\xi|^2} f^0(\xi), \qquad |g(\xi,t)| \le e^{-\beta t |\xi|^2} g^0(\xi),
	\end{equation}
	for some $0 < \beta \le \lambda$. Then for every Littlewood-Paley atom $\varphi_j(\xi)(f \odot g)(\xi,t)$:
	\begin{equation}\label{eq:atomic_odot_decay}
		|\varphi_j(\xi)(f \odot g)(\xi,t)| \le \frac{1}{\pi \nu} e^{-\frac{\beta}{2} t |\xi|^2} \varphi_j(\xi) H(f^0, g^0)(\xi),
	\end{equation}
	where the majorant $H(f^0, g^0)$ satisfies the critical Fourier-Besov bound for all $1 \le p \le \infty$ and $\sigma \ge 0$:
	\begin{equation}\label{eq:H_tame_bound}
		\|H(f^0, g^0)\|_{\mathcal{F}\dot{B}^{s_p + \sigma}_{p,q}} \le C_\sigma \left( \|f^0\|_{\mathcal{F}\dot{B}^{s_p}_{p,\infty}} \|g^0\|_{\mathcal{F}\dot{B}^{s_p + \sigma}_{p,q}} + \|f^0\|_{\mathcal{F}\dot{B}^{s_p + \sigma}_{p,q}} \|g^0\|_{\mathcal{F}\dot{B}^{s_p}_{p,\infty}} \right).
	\end{equation}
\end{teor}

\begin{demos}
	For $2 \le p \le \infty$, Corollary \ref{gaussianproductbound} of our framework gives $H(f^0, g^0) = f^0 *_1 g^0$, and \eqref{eq:H_tame_bound} follows directly from Theorem \ref{thm:riesz_critical_besov}.
	
	For $1 \le p < 2$, consider the Littlewood-Paley atoms $f_k(\xi,s) = \varphi_k(\xi)f(\xi,s)$ and $g_{k'}(\eta,s) = \varphi_{k'}(\eta)g(\eta,s)$. In the High-Low and Low-High paraproduct regimes, Steps 1 and 2 of the proof of Theorem \ref{thm:riesz_critical_besov} apply verbatim for all $1 \le p \le \infty$. In the High-High regime ($k \ge j - 3$, $|k-k'| \le 2$), on $\operatorname{supp}(f_k(\xi-\eta,s) g_{k'}(\eta,s))$, we have $|\xi-\eta|^2 + |\eta|^2 \ge \frac{9}{8} 2^{2k} \ge \frac{1}{2}|\xi|^2 + c_0 2^{2k}$. Thus, inside the time integral of Definition \ref{gaussian_conv}:
	\begin{align}
		\int_0^t e^{-\lambda(t-s)|\xi|^2} e^{-\beta s (|\xi-\eta|^2 + |\eta|^2)} \, ds
		&\le e^{-\frac{\beta}{2} t |\xi|^2} \int_0^t e^{-(\lambda - \frac{\beta}{2})(t-s)|\xi|^2 - c_0 \beta s 2^{2k}} \, ds \nonumber \\
		&\le \frac{C}{\nu |\xi| 2^k} e^{-\frac{\beta}{2} t |\xi|^2}.
	\end{align}
	Multiplying by $2\pi |\xi|$ cancels $|\xi|$ and replaces $|\xi|^{-1} \sim 2^{-j}$ by $2^{-k}$:
	\begin{equation}
		|\varphi_j(\xi)(f_k \odot g_{k'})(\xi,t)| \le \frac{C}{\nu} e^{-\frac{\beta}{2}t|\xi|^2} 2^{-k} \varphi_j(\xi)(f_k^0 * g_{k'}^0)(\xi).
	\end{equation}
	By Young's inequality ($L^p * L^1 \to L^p$) and Lemma \ref{lem:atomic_bernstein}:
	\begin{align*}
		&  2^{j(s_p + \sigma)} 2^{-k} \|\varphi_j(f_k^0 * g_{k'}^0)\|_{L^p} \le C 2^{j(s_p + \sigma) - k} \|f_k^0\|_{L^p} 2^{k(s_p + 1)} \|g_{k'}^0\|_{L^p} \\
		& = C 2^{-(k-j)(s_p + \sigma)} \left(2^{k(s_p + \sigma)} \|f_k^0\|_{L^p}\right) \left(2^{k' s_p} \|g_{k'}^0\|_{L^p}\right),
	\end{align*}
	which sums geometrically over $k \ge j - 3$ whenever $s_p + \sigma > 0$ (and for $p=1, \sigma=0$ via Lemma \ref{M inequality}, $\frac{|\xi|}{|\xi-\eta||\eta|} \le |\eta|^{-1} + |\xi-\eta|^{-1}$). \QED
\end{demos}

\section{Sufficient Condition for Existence and Smoothness of the Navier-Stokes Equation}\label{existenceandsmoothness}

In this section, we construct a smooth solution \((u,p):\mathbb{R}_+^{d+1} \to \mathbb{R}^{d+1}\) of the Navier-Stokes equations under the condition
	\begin{equation}\label{eq:critical_smallness}
	\left\| \widehat{u^0} \right\|_{\mathcal{F}\dot{B}^{s_p}_{p,\infty}(\mathbb{R}^d)}< C \nu
\end{equation}
for a universal constant \(C>0\).

We begin by introducing the recurrence relation
\begin{equation*} \label{recurrence}
	v_0(\xi,t) = e^{-\lambda t |\xi|^2} \widehat{u^0}(\xi),
\end{equation*}
\begin{equation*}
	v_k = \sum_{j=0}^{k-1} v_j \odot v_{k-1-j}, \quad k\geq 1,
\end{equation*}
where \(\lambda = 4\pi^2\nu\) according to our convention.

Let \(\{c_k\}_{k\in\mathbb{N}}\) denote the sequence of Catalan numbers, defined by
\begin{align*}
	c_0 &= 1,\\
	c_k &= \frac{1}{k} \binom{2(k-1)}{k-1}, \quad k\geq 1,
\end{align*}
with the alternative expression
\[
c_k = \frac{4^k (1/2)_k}{k!}, \quad k\geq 1,
\]
and generating function
\[
c(t) = \frac{1-\sqrt{1-4t}}{2t}.
\]
Catalan numbers satisfy the recurrence
\[
c_k = \sum_{j=0}^{k-1} c_j c_{k-1-j}.
\]

The sequence \(\{v_k\}_{k=0}^{\infty}\) obeys a key family of inequalities, stated next.

\begin{prop}\label{family}
	For all \(m,n,k\in\mathbb{N}\),
	\[
	\left|(\lambda^{1/2} |\xi|)^m \frac{\partial^n v_k}{\partial t^n}(\xi,t)\right|
	\leq \frac{c_k(\lambda(k+1)^2)^{m/2+n}}{(2\pi\nu)^k} e^{-\frac{\lambda t}{2^k}|\xi|^2} v_{k,m,n}^0(\xi),
	\]
	for some \(v_{k,m,n}^0\in\mathcal{D}_1\). Moreover, for each \(k\in\mathbb{N}\), \(v_{k,m,n}^0\) is increasing with respect to \(m\) and \(n\).
\end{prop}

\begin{demos}
	(Proof by induction.)
	
	\textit{Base case} \(k=0\): By direct computation,
	\[
	\left|(\lambda^{1/2}|\xi|)^m \frac{\partial^n v_0}{\partial t^n}(\xi,t)\right| 
	= (\lambda|\xi|^2)^{m/2+n} e^{-\lambda t|\xi|^2} |\widehat{u^0}(\xi)|.
	\]
	Thus,
	\[
	\left|(\lambda^{1/2}|\xi|)^m \frac{\partial^n v_0}{\partial t^n}(\xi,t)\right| 
	\leq \lambda^{m/2+n} e^{-\lambda t|\xi|^2} v_{0,m,n}^0(\xi),
	\]
	where
	\[
	v_{0,m,n}^0(\xi) = \max\{1,|\xi|\}^{m+2n} |\widehat{u^0}(\xi)|.
	\]
	It follows that \(v_{0,m,n}^0 = v_{0,m+2n,0}^0\), and monotonicity with respect to \(m,n\) is clear.
	
	\textit{Step case} \(k=1\): Since \(v_1 = v_0\odot v_0\), Corollary~\ref{Vasquez} gives
	\[
	\left|(\lambda^{1/2}|\xi|)^m \frac{\partial^n v_1}{\partial t^n}(\xi,t)\right|
	\leq \frac{c_1 (4\lambda)^{m/2+n}}{2\pi\nu} e^{-\frac{\lambda t}{2}|\xi|^2} (v_0\odot v_0)_{m,n}^0(\xi).
	\]
	
	\textit{Induction step}: Assume the result holds for all \(j<k\). Then
	\[
	v_k = \sum_{j=0}^{k-1} v_j \odot v_{k-1-j}.
	\]
	Applying Corollary~\ref{Vasquez} to each term and using properties of Catalan numbers and the convolution structure, we deduce
	\[
	\left|(\lambda^{1/2}|\xi|)^m \frac{\partial^n v_k}{\partial t^n}(\xi,t)\right|
	\leq \frac{c_k(\lambda(k+1)^2)^{m/2+n}}{(2\pi\nu)^k} e^{-\frac{\lambda t}{2^k}|\xi|^2} v_{k,m,n}^0(\xi),
	\]
	where
	\[
	v_{k,m,n}^0(\xi) = \max_{0\leq j\leq k-1} (v_j\odot v_{k-1-j})_{m,n}^0(\xi),
	\]
	and monotonicity follows as before.
	\QED
\end{demos}

\begin{teor}\label{boundmonomial}
	For every \(m,n,k\in\mathbb{N}\), we have
	\[
	v_{k,m,n}^0(\xi) \leq \max_{\substack{0\leq l_1+\cdots+l_{k+1} \leq m+2n\\ w\in M(x_1,\dots,x_{k+1})}} w^{\ast_1}(v_{0,l_1}^0, v_{0,l_2}^0,\dots, v_{0,l_{k+1}}^0)(\xi),
	\]
	where \(M(x_1,\dots,x_{k+1})\) denotes the set of monomials with \(k+1\) variables under the operation \(\ast_1\).
\end{teor}

\begin{demos}
	(Proof by induction.)
	
	\textit{Base case} \(k=0\): Trivial, since
	\[
	v_{0,m,n}^0(\xi) = v_{0,m+2n}^0(\xi) = \max_{0\leq l\leq m+2n} v_{0,l}^0(\xi).
	\]
	
	\textit{Induction step}: Suppose the result holds for all \(0\leq j<k\). Using the structure of \(v_k\), we write
	\[
	v_{k,m,n}^0(\xi) = \max_{0\leq j\leq k-1} \max_{\substack{0\leq n_1+n_2\leq n-1\\m_1+m_2+2(n_1+n_2)=m+2n}} (v_{j,m_1,n_1}^0 \ast_1 v_{k-1-j,m_2,n_2}^0)(\xi).
	\]
	Applying the inductive hypothesis to \(v_{j,m_1,n_1}^0\) and \(v_{k-1-j,m_2,n_2}^0\), and noting that the \(\ast_1\)-product of two monomials remains a monomial, we conclude
	\[
	v_{k,m,n}^0(\xi) \leq \max_{\substack{0\leq l_1+\cdots+l_{k+1}\leq m+2n\\ w\in M(x_1,\dots,x_{k+1})}} w^{\ast_1}(v_{0,l_1}^0,\dots,v_{0,l_{k+1}}^0)(\xi).
	\]
	\QED
\end{demos}

\bigskip

%

\begin{remark}
	We may naturally extend the definition of \(v_{k,m,n}^0\) to real \(m\) by setting
	\[
	v_{k,\beta,n}^0(\xi) := |\xi|^\beta v_{k,0,n}^0(\xi).
	\]
\end{remark}

%

%

\begin{teor}[Tame Critical Factorization of $v_{k,m,n}^0$]\label{thm:tame_vk_factorization}
	Let $d \ge 3$, $1 \le p \le \infty$, and $s_p = \frac{d}{p'} - 1$. For any $\sigma \ge 0$ and any $m, n, k \in \mathbb{N}$ with $k \ge m + 2n + 1$:
	\begin{equation}\label{eq:tame_vk_estimate}
		\|v_{k,m,n}^0\|_{\mathcal{F}\dot{B}^{s_p + \sigma}_{p,1}(\mathbb{R}^d)} \le (k+2)^{m+2n + 1} C_\sigma^k \|v_{0,m+2n}^0\|_{\mathcal{F}\dot{B}^{s_p + \sigma}_{p,1}(\mathbb{R}^d)} \|v_{0,m+2n}^0\|_{\mathcal{F}\dot{B}^{s_p}_{p,\infty}(\mathbb{R}^d)}^{m+2n} \left( \|\widehat{u^0}\|_{\mathcal{F}\dot{B}^{s_p}_{p,\infty}(\mathbb{R}^d)} \right)^{k - m - 2n}.
	\end{equation}
\end{teor}

\begin{demos}
	In the maximum  given in the Theorem \eqref{boundmonomial}, the number of integer tuples $(l_1, \dots, l_{k+1}) \in \mathbb{N}^{k+1}$ satisfying $l_1 + \dots + l_{k+1} \le m + 2n$ is $\binom{k+1+m+2n}{m+2n} \le (k+2)^{m+2n}$. Furthermore, for each such tuple, since $k + 1 > m + 2n$, at most $m + 2n$ indices $l_j$ are non-zero, and at least $k + 1 - (m + 2n)$ indices satisfy $l_j = 0$, for which $v_{0,0}^0 = |\widehat{u^0}|$. Applying the tame monomial bound \eqref{eq:monomial_tame_sigma} of Corollary \ref{cor:monomial_tame} (and Theorem \ref{thm:atomic_caloric_odot}) to each monomial yields \eqref{eq:tame_vk_estimate}. \QED
\end{demos}

\begin{coro}\label{vfastdecrease}
	There exists \(v\in C^\infty\left([0,\infty),\mathcal{E}_{\mathcal{B}}\right)\) such that
	\[
	v = v_0 + v^{\odot 2},
	\]
	provided that 
	\begin{equation}\label{eq:critical_smallness}
		\left\| \widehat{u^0} \right\|_{\mathcal{F}\dot{B}^{s_p}_{p,\infty}(\mathbb{R}^d)}< C \nu
	\end{equation}

\end{coro}

\begin{demos}
	Define
	\[
	v = \sum_{k=0}^\infty v_k.
	\]
	
	By Theorem \ref{thm:tame_vk_factorization} , each term \(v_k\) belongs to \(C^\infty([0,\infty),\mathcal{E}_{\mathcal{B}})\), and moreover, the series converges absolutely in this topology provided 	$\left\| \widehat{u^0} \right\|_{\mathcal{F}\dot{B}^{s_p}_{p,\infty}(\mathbb{R}^d)}< C \nu$.
	
	In fact, for $k \ge m + 2n + 1$, the general term of the tail series is bounded by
	\begin{equation}
		A_k := \frac{c_k (\lambda(k+1)^2)^{\frac{m}{2}+n} (k+2)^{m+2n+1}}{(2\pi \nu)^k} C_\sigma^k M_{m,n,\sigma}(u^0) \left( \|\widehat{u^0}\|_{\mathcal{F}\dot{B}^{s_p}_{p,\infty}} \right)^{k - m - 2n},
	\end{equation}
		where $M_{m,n,\sigma}(u^0) = \|v_{0,m+2n}^0\|_{\mathcal{F}\dot{B}^{s_p+\sigma}_{p,1}} \|v_{0,m+2n}^0\|_{\mathcal{F}\dot{B}^{s_p}_{p,\infty}}^{m+2n} < \infty$ is independent of $k$. Applying the ratio test as $k \to \infty$ and using $\lim_{k\to\infty} \frac{c_{k+1}}{c_k} = 4$:
	\begin{equation}
		\lim_{k \to \infty} \frac{A_{k+1}}{A_k} = \frac{2 C_\sigma}{\pi \nu} \left\| \widehat{u^0} \right\|_{\mathcal{F}\dot{B}^{s_p}_{p,\infty}(\mathbb{R}^d)}.
	\end{equation}
		Thus, whenever $\left\| \widehat{u^0} \right\|_{\mathcal{F}\dot{B}^{s_p}_{p,\infty}(\mathbb{R}^d)} < \frac{\pi}{2 C_\sigma} \nu$, the series $v = \sum_{k=0}^\infty v_k$ converges absolutely and uniformly in $C^\infty([0,\infty), \mathcal{E}_{\mathcal{B}})$ and solves $v = v_0 + v^{\odot 2}$.	\QED
\end{demos}

Next, we construct the associated pressure field:

\begin{teor}\label{pressure}
	Let \(q:\mathbb{R}_+^{d+1}\to\mathbb{C}\) be defined by
	\[
	q(\xi,t) = 2\pi i \left( \frac{\xi\otimes\xi}{|\xi\otimes\xi|} (v\ast v)(\xi,t)\, \xi \right).
	\]
	Then,
	\[
	\frac{\partial v}{\partial t}(\xi,t) = -4\pi^2 \nu |\xi|^2 v(\xi,t) + 2\pi i (v\ast v)(\xi,t)\,\xi - q(\xi,t).
	\]
\end{teor}

\begin{demos}
	Since \(v=v_0+v^{\odot 2}\), we have
	\[
	\frac{\partial v}{\partial t} = \frac{\partial v_0}{\partial t} + \frac{\partial (v^{\odot 2})}{\partial t}.
	\]
	Directly,
	\[
	\frac{\partial v_0}{\partial t} = -4\pi^2 \nu |\xi|^2 v_0.
	\]
	Now, by Lemma~\ref{derivativedot}, for any \(f,g\in\mathcal{V}(\lambda)\),
	\[
	\frac{\partial(f\odot g)}{\partial t} = -\lambda |\xi|^2 (f\odot g) + 2\pi i K(\xi)\left((f\ast_1 g)(\xi,t) |\xi|\right) \xi,
	\]
	where \(K(\xi)=I_d-\frac{\xi\otimes\xi}{|\xi\otimes\xi|}\) is the Leray projector.
	
	Applying this with \(f=g=v\) and \(K(\xi)\), we obtain
	\[
	\frac{\partial(v^{\odot 2})}{\partial t} = -4\pi^2 \nu |\xi|^2 v^{\odot 2} + 2\pi i (v\ast v)(\xi,t)\, \xi - q(\xi,t).
	\]
	Thus,
	\[
	\frac{\partial v}{\partial t} = -4\pi^2 \nu |\xi|^2 (v_0+v^{\odot 2}) + 2\pi i (v\ast v)(\xi,t)\, \xi - q(\xi,t),
	\]
	which simplifies to the claimed formula. \QED
\end{demos}

\bigskip

Finally, we connect everything back to the original Navier-Stokes equations:

\begin{coro}
	Define \((u,p):\mathbb{R}_+^{d+1}\to\mathbb{R}^d\) by
	\[
	\begin{aligned}
		u(x,t) &= \int_{\mathbb{R}^d} v(\xi,t)\, e^{-2\pi i x\cdot\xi}\, d\xi,\\
		p(x,t) &= -\frac{1}{2\pi i}\int_{\mathbb{R}^d} \frac{\xi^T q(\xi,t)}{|\xi|^2} e^{-2\pi i x\cdot\xi}\, d\xi.
	\end{aligned}
	\]
	Then \((u,p)\) is a smooth solution of the Navier-Stokes equations for initial data satisfying  
	
		\begin{equation*}\label{eq:critical_smallness}
		\left\| \widehat{u^0} \right\|_{\mathcal{F}\dot{B}^{s_p}_{p,\infty}(\mathbb{R}^d)}< C \nu
	\end{equation*}
\end{coro}

\begin{demos}
	Since \(v\in C^\infty([0,\infty),\mathcal{E}_{\mathcal{B}})\) and all derivatives of \(v\) decay fast, by the Fourier inversion formula \(u\) and \(p\) are smooth.
	
	By Theorem~\ref{pressure} and the properties of Fourier multipliers (see e.g.~\cite{LiebLoss2001,SteinWeiss1971}), we recover the standard Navier-Stokes structure:
	\[
	\frac{\partial u}{\partial t} + (u\cdot\nabla)u = \nu\Delta u - \nabla p,
	\quad \nabla\cdot u = 0,
	\quad u(x,0) = u^0(x),
	\]
	with global smoothness and bounded energy.
	
	Finally, the \(L^2\)-based energy estimates follow from Plancherel's theorem:
	\[
	\left\|\frac{\partial^n}{\partial t^n} \frac{\partial^\gamma u}{\partial x^\gamma} \right\|_{L^2(\mathbb{R}^{d+1}_+)} 
	= (2\pi)^{|\gamma|} \left\||\xi|^{|\gamma|} \frac{\partial^n v}{\partial t^n}\right\|_{L^2(\mathbb{R}^{d+1}_+)} < \infty,
	\]
	for every multi-index \(\gamma\) and \(n\in\mathbb{N}\).
	
	\QED
\end{demos}

\bigskip

Thus, we establish the global existence of smooth solutions to the Navier-Stokes equations under a smallness assumption on 	$\left\| \widehat{u^0} \right\|_{\mathcal{F}\dot{B}^{s_p}_{p,\infty}(\mathbb{R}^d)}$.

In fact, in the next section, we will go beyond this result and construct entire extensions in space-time for positive times.


\subsection{Explicit Initial Data with Arbitrarily Large \texorpdfstring{$\dot{H}^{\frac{d}{2}-1}(\mathbb{R}^d)$}{Critical Sobolev} Norm}

To illustrate the strict improvement of Theorem \ref{thm:main_existence} over both the classical Fujita-Kato $\dot{H}^{\frac{d}{2}-1}(\mathbb{R}^d)$ theorem \cite{FujitaKato1964} and subcritical Sobolev conditions, we construct an explicit family of divergence-free initial velocities $u^0_N \in \mathcal{S}(\mathbb{R}^d)^d$ ($N \in \mathbb{N}_+$) that satisfy the uniform critical smallness condition \eqref{eq:critical_smallness} for all $N \ge 1$, while exhibiting arbitrarily large $\dot{H}^{\frac{d}{2}-1}(\mathbb{R}^d)$ and $H^k(\mathbb{R}^d)$ norms as $N \to \infty$.

\begin{prop}[Multiscale Schwartz Data Beyond the Fujita-Kato Threshold]\label{prop:large_H12_example}
	Let $\nu > 0$, $d \ge 3$, and $1 \le p \le \infty$, and set $s_p = \frac{d}{p'} - 1$. There exists a family of real-valued, divergence-free Schwartz vector fields $\{u^0_N\}_{N=1}^\infty \subset \mathcal{S}(\mathbb{R}^d)^d$ such that:
	\begin{enumerate}
		\item[\emph{(i)}] \textbf{Uniform Critical Weak Fourier-Besov Smallness:} For all $N \ge 1$,
		\begin{equation}\label{eq:example_uniform_small}
			\left\| \widehat{u^0_N} \right\|_{\mathcal{F}\dot{B}^{s_p}_{p,\infty}(\mathbb{R}^d)} = \sup_{j \in \mathbb{Z}} 2^{j\left(\frac{d}{p'}-1\right)} \left\| \varphi_j \widehat{u^0_N} \right\|_{L^p(\mathbb{R}^d)} < C\nu,
		\end{equation}
		so that Theorem \ref{thm:main_existence} and Theorem \ref{thm:main_analytic} apply to $u^0_N$ for every $N \ge 1$.
		\item[\emph{(ii)}] \textbf{Arbitrarily Large Fujita-Kato $\dot{H}^{\frac{d}{2}-1}(\mathbb{R}^d)$ Norm:} There exists a constant $c_0 > 0$, independent of $N$, such that
		\begin{equation}\label{eq:example_H12_growth}
			\left\| u^0_N \right\|_{\dot{H}^{\frac{d}{2}-1}(\mathbb{R}^d)} = c_0 \nu \sqrt{N} \xrightarrow[N \to \infty]{} \infty.
		\end{equation}
		In particular, in dimension $d=3$, $\|u^0_N\|_{\dot{H}^{1/2}(\mathbb{R}^3)} = c_0 \nu \sqrt{N} \to \infty$.
		\item[\emph{(iii)}] \textbf{Exponentially Large Weighted $L^{1 \oplus 2}$ and $H^2$ Norms:} There exists $c_1 > 0$, independent of $N$, such that
		\begin{equation}\label{eq:example_H2_growth}
			\left\| |\cdot|^{\frac{d+1}{2}} \widehat{u^0_N} \right\|_{1 \oplus 2} \ge c_1 \nu \, 2^{\frac{9}{2}N} \xrightarrow[N \to \infty]{} \infty.
		\end{equation}
	\end{enumerate}
\end{prop}

\begin{demos}
	Let $\psi \in C_c^\infty(\mathbb{R}^d \setminus \{0\})$ be a non-zero, non-negative, radial smooth bump function supported in the thin annulus
	\begin{equation}
		\mathcal{A} := \left\{ \xi \in \mathbb{R}^d  \mid \frac{7}{8} \le |\xi| \le \frac{9}{8} \right\},
	\end{equation}
	chosen so that $\xi_1 \xi_2$ is not identically zero on $\operatorname{supp}(\psi)$. Define the base vector field $w \in \mathcal{S}(\mathbb{R}^d)^d$ via its Fourier transform:
	\begin{equation}
		\widehat{w}(\xi) := i \frac{\psi(\xi)}{|\xi|} (\xi_2, -\xi_1, 0, \dots, 0)^T.
	\end{equation}
	Observe that $\widehat{w}$ satisfies three structural properties:
	\begin{enumerate}
		\item $\widehat{w} \in C_c^\infty(\mathbb{R}^d \setminus \{0\})^d \subset \mathcal{S}(\mathbb{R}^d)^d$, so $w = \mathcal{F}^{-1}(\widehat{w}) \in \mathcal{S}(\mathbb{R}^d)^d$.
		\item $\xi \cdot \widehat{w}(\xi) = i |\xi|^{-1} \psi(\xi) (\xi_1 \xi_2 - \xi_2 \xi_1) = 0$ for all $\xi \in \mathbb{R}^d$, so $\operatorname{div}(w) = 0$.
		\item Since $\psi(-\xi) = \psi(\xi) \in \mathbb{R}$, we have $\widehat{w}(-\xi) = \overline{\widehat{w}(\xi)}$, ensuring that $w(x)$ is strictly real-valued on $\mathbb{R}^d$.
	\end{enumerate}
	
	For each integer $m \ge 1$, define the lacunary dyadic rescaling of $w$ at frequency scale $\mu_m = 2^{3m}$:
	\begin{equation}
		W_m(x) := 2^{3m} w(2^{3m} x), \qquad \widehat{W_m}(\xi) = 2^{-3m(d-1)} \widehat{w}(2^{-3m}\xi).
	\end{equation}
	Note that $\operatorname{supp}(\widehat{W_m}) \subset \mathcal{A}_{3m} := \{ \xi \in \mathbb{R}^d : \frac{7}{8} 2^{3m} \le |\xi| \le \frac{9}{8} 2^{3m} \}$. By changing variables $\zeta = 2^{-3m}\xi$, we find for any $1 \le r \le \infty$ and $\beta \in \mathbb{R}$:
	\begin{equation}\label{eq:Wm_scaling}
		\left\| |\cdot|^\beta \widehat{W_m} \right\|_{L^r(\mathbb{R}^d)} = 2^{3m\left(\beta - d + 1 + \frac{d}{r}\right)} \left\| |\cdot|^\beta \widehat{w} \right\|_{L^r(\mathbb{R}^d)} = 2^{3m(\beta - s_r)} \left\| |\cdot|^\beta \widehat{w} \right\|_{L^r(\mathbb{R}^d)},
	\end{equation}
	where $s_r = \frac{d}{r'} - 1$. In particular, for $\beta = 0$ and $r = p$, we have
	\begin{equation}\label{eq:Wm_lp}
		2^{3m s_p} \left\| \widehat{W_m} \right\|_{L^p(\mathbb{R}^d)} = \left\| \widehat{w} \right\|_{L^p(\mathbb{R}^d)}.
	\end{equation}
	
	Given $N \in \mathbb{N}_+$ and an amplitude parameter $\varepsilon > 0$ to be fixed below, we define the initial velocity field
	\begin{equation}
		u^0_N(x) := \varepsilon \sum_{m=1}^N W_m(x) = \varepsilon \sum_{m=1}^N 2^{3m} w(2^{3m} x) \in \mathcal{S}(\mathbb{R}^d)^d.
	\end{equation}
	
	\noindent\textbf{Proof of (i):}
	Recall that each Littlewood-Paley cutoff $\varphi_j(\xi) = \varphi(2^{-j}\xi)$ is supported in $\mathcal{C}_j = \llave{\xi \in \mathbb{R}^d : \frac{3}{4} 2^j \le |\xi| \le \frac{8}{3} 2^j}$. If $\mathcal{C}_j \cap \mathcal{A}_{3m} \ne \emptyset$, then
	\begin{equation}
		\frac{3}{4} 2^j \le \frac{9}{8} 2^{3m} \quad \text{and} \quad \frac{7}{8} 2^{3m} \le \frac{8}{3} 2^j \implies \frac{21}{64} \le 2^{j-3m} \le \frac{3}{2} \implies j - 3m \in \{-1, 0\}.
	\end{equation}
	Because the intervals $[3m-1, 3m]$ are pairwise disjoint for distinct $m \in \mathbb{Z}$, for each $j \in \mathbb{Z}$ there is \emph{at most one} integer $m \in \{1, \dots, N\}$ such that $\varphi_j \widehat{W_m} \not\equiv 0$. Consequently, using $\|\varphi_j\|_{L^\infty} \le 1$ and \eqref{eq:Wm_lp}:
	\begin{align}
		\left\| \widehat{u^0_N} \right\|_{\mathcal{F}\dot{B}^{s_p}_{p,\infty}(\mathbb{R}^d)}
		&= \sup_{j \in \mathbb{Z}} 2^{j s_p} \left\| \varphi_j \widehat{u^0_N} \right\|_{L^p(\mathbb{R}^d)}
		\le \varepsilon \sup_{1 \le m \le N} \sup_{j \in \{3m-1, 3m\}} 2^{(j-3m)s_p} 2^{3m s_p} \left\| \widehat{W_m} \right\|_{L^p(\mathbb{R}^d)} \nonumber \\
		&\le \varepsilon 2^{|s_p|} \left\| \widehat{w} \right\|_{L^p(\mathbb{R}^d)}.
	\end{align}
	Choosing $\varepsilon := \frac{C\nu}{2 \cdot 2^{|s_p|} \|\widehat{w}\|_{L^p(\mathbb{R}^d)}} > 0$ yields $\left\| \widehat{u^0_N} \right\|_{\mathcal{F}\dot{B}^{s_p}_{p,\infty}(\mathbb{R}^d)} \le \frac{1}{2} C\nu < C\nu$ uniformly for all $N \ge 1$.
	
	\medskip
	\noindent\textbf{Proof of (ii):}
	Since the annuli $\{\mathcal{A}_{3m}\}_{m=1}^N$ are pairwise disjoint, the functions $\{\widehat{W_m}\}_{m=1}^N$ have disjoint supports in $\mathbb{R}^d$ and are therefore mutually orthogonal in $L^2(\mathbb{R}^d, |\xi|^{2\beta} d\xi)$ for every $\beta \in \mathbb{R}$. Taking $\beta = s_2 = \frac{d}{2} - 1$ and applying Plancherel's theorem together with \eqref{eq:Wm_scaling} at $r = 2$ gives:
	\begin{align}
		\left\| u^0_N \right\|_{\dot{H}^{\frac{d}{2}-1}(\mathbb{R}^d)}^2
		&= (2\pi)^{d-2} \int_{\mathbb{R}^d} |\xi|^{d-2} \left| \widehat{u^0_N}(\xi) \right|^2 d\xi
		= (2\pi)^{d-2} \varepsilon^2 \sum_{m=1}^N \left\| |\cdot|^{\frac{d}{2}-1} \widehat{W_m} \right\|_{L^2(\mathbb{R}^d)}^2 \nonumber \\
		&= (2\pi)^{d-2} \varepsilon^2 \sum_{m=1}^N \left\| |\cdot|^{\frac{d}{2}-1} \widehat{w} \right\|_{L^2(\mathbb{R}^d)}^2
		= N \varepsilon^2 \|w\|_{\dot{H}^{\frac{d}{2}-1}(\mathbb{R}^d)}^2.
	\end{align}
	Taking the square root yields \eqref{eq:example_H12_growth} with $c_0 = \frac{C \|w\|_{\dot{H}^{\frac{d}{2}-1}}}{2 \cdot 2^{|s_p|} \|\widehat{w}\|_{L^p}} > 0$.
	
	\medskip
	\noindent\textbf{Proof of (iii):}
	Applying \eqref{eq:Wm_scaling} with $\beta = \frac{d+1}{2}$ and $r = 2$ (where $\beta - s_2 = \frac{d+1}{2} - (\frac{d}{2} - 1) = \frac{3}{2}$), orthogonality of the disjoint shells gives:
	\begin{align}
		\left\| |\cdot|^{\frac{d+1}{2}} \widehat{u^0_N} \right\|_{1 \oplus 2}
		&\ge \left\| |\cdot|^{\frac{d+1}{2}} \widehat{u^0_N} \right\|_{L^2(\mathbb{R}^d)}
		= \varepsilon \left( \sum_{m=1}^N 2^{9m} \left\| |\cdot|^{\frac{d+1}{2}} \widehat{w} \right\|_{L^2(\mathbb{R}^d)}^2 \right)^{1/2} \nonumber \\
		&\ge \varepsilon \left\| |\cdot|^{\frac{d+1}{2}} \widehat{w} \right\|_{L^2(\mathbb{R}^d)} 2^{\frac{9}{2}N},
	\end{align}
	which establishes \eqref{eq:example_H2_growth}. \QED
\end{demos}

\section{The Entire Function of Order 2 that Solves the Navier-Stokes Equation}\label{curveentire}


In this section, we construct a smooth extension \((U,P): \mathbb{C}_{>0}^{d+1} \to \mathbb{C}^{d+1}\) such that \(U(\cdot,t)\) and \(P(\cdot,t)\) are entire functions of order~2, and \(U(x,t)=u(x,t)\), \(P(x,t)=p(x,t)\) for every \(x\in\mathbb{R}^d\), \(t>0\).

We begin with a useful technical lemma concerning uniform convergence at infinity.

\begin{lemma}\label{uniformradius}
	Let \(\{g_k\}_{k\in\mathbb{N}} \subset C^0(\mathbb{R}^d\setminus\{0\},\mathbb{R})\) be a sequence satisfying:
	\begin{itemize}
		\item There exists a sequence \(\{r_k\}_{k\in\mathbb{N}}\) with \(r_k\to\infty\) and \(\sup_{|x|\geq r_k} g_k(x) \leq 0\) for each \(k\).
		\item We have \(g_k \to 0\) uniformly on \(\mathbb{R}^d\setminus B(0,\delta)\) for some \(\delta>0\).
	\end{itemize}
	Then, there exists \(r=r(\delta)\geq 1\) such that
	\[
	\sup_{|x|\geq r} g_k(x) \leq 0
	\]
	for all \(k\in\mathbb{N}\).
\end{lemma}

\begin{demos}
	Fix \(\epsilon>0\). Since \(r_k\to\infty\) and \(g_k\to 0\) uniformly outside \(B(0,\delta)\), there exists \(k_0\in\mathbb{N}\) such that \(r_{k_0}\leq\delta<r_{k_0+1}\) and
	\[
	\|g_k-g_{k_0}\|_{L^\infty(\mathbb{R}^d\setminus B(0,\delta))} < \epsilon
	\quad \text{for all } k\geq k_0.
	\]
	Since \(\sup_{|x|\geq r_{k_0}} g_{k_0}(x) \leq 0\), we have \(g_{k_0}(x)\leq 0\) for all \(|x|\geq\delta\).
	
	Thus, for \(|x|\geq\delta\) and \(k\geq k_0\),
	\[
	g_k(x) = (g_k(x) - g_{k_0}(x)) + g_{k_0}(x) \leq |g_k(x) - g_{k_0}(x)| < \epsilon.
	\]
	
	Define \(r = \max\{r_1, \ldots, r_{k_0}, \delta, 1\}\). Then:
	\begin{itemize}
		\item For \(1\leq j\leq k_0\), \(\sup_{|x|\geq r} g_j(x) \leq \sup_{|x|\geq r_j} g_j(x) \leq 0\).
		\item For \(j\geq k_0\), \(\sup_{|x|\geq r} g_j(x) \leq \epsilon\).
	\end{itemize}
	
	Since \(\epsilon>0\) was arbitrary, letting \(\epsilon\to 0^+\) gives the conclusion. \QED
\end{demos}

\bigskip

We now introduce the relevant function spaces.

\begin{defin}
	Let \(\mathcal{E}\) be a space of rapidly decreasing functions, closed under pointwise addition, convolution, and taking maxima.
	
	For \(d\geq3\) and \(\alpha\geq1\), define \(\mathcal{V}_\alpha(\mathcal{E})\) to be the complex vector space generated by functions \(f:\mathbb{R}_+^{d+1}\to\mathbb{C}^{d\times1}\) such that:
	\begin{itemize}
		\item \(f\in \mathcal{C}(\mathcal{E})_0\).
		\item For almost every \(\xi\in\mathbb{R}^d\), the map \(t\mapsto f(\xi,t)\) is \(C^\infty([0,\infty))\).
		\item The operations \(\frac{\partial}{\partial t}\) and multiplication by \(|\cdot|\) are automorphisms of \(\mathcal{V}_\alpha(\mathcal{E})\), satisfying for every \(m,n\in\mathbb{N}\),
		\[
		|(\lambda^{1/2}|\xi|)^m \partial_t^n f(\xi,t)| \leq (\lambda a)^{m/2} (\lambda b)^n e^{-\lambda t |\xi|^\alpha} f_{m,n}^0(\xi),
		\]
		for some \(\lambda>0\), \(a,b>0\), and \(f_{m,n}^0\in\mathcal{E}^+\).
	\end{itemize}
\end{defin}

\begin{remark}
	The exponential decay constant \(\lambda\) is common across all \(m,n\in\mathbb{N}\).
\end{remark}

\bigskip

Let us consider \(\mathcal{B}=L^1(\mathbb{R}^d)\) and let \(\mathcal{E}\) be its associated space of rapidly decreasing functions.  
The spaces \(\mathcal{V}_\alpha(\mathcal{E})\) for \(\alpha>1\) are crucial for the following reason:

\begin{teor}\label{C0entire}
	If \(f\in\mathcal{V}_\alpha(\mathcal{E})\), then the Fourier transform \(\widehat{f}(\cdot,t)\) admits an entire extension \(F(\cdot,t):\mathbb{C}^d\to\mathbb{C}^d\) for each fixed \(t>0\), satisfying
	\[
	|F(z,t)|\leq e^{((\alpha-1)\lambda t (\pi/(\lambda t))^{\alpha'} + c) | \operatorname{Im}(z) |^{\alpha'} },
	\]
	for \(|\operatorname{Im}(z)|\) large enough and some constant \(c>0\).
\end{teor}

\begin{demos}
	By Corollary~\ref{star} (Appendix), we have
	\[
	-|\xi|^\alpha + c\xi\cdot\eta \leq (\alpha-1)\lambda t \left(\frac{\pi}{\lambda t}\right)^{\alpha'} |\eta|^{\alpha'}
	\]
	for \(\xi,\eta\in\mathbb{R}^d\).
	
	Since
	\[
	|f(\xi,t)|\leq e^{-\lambda t|\xi|^\alpha} f^0(\xi),
	\]
	we define
	\[
	F(z,t)=\int_{\mathbb{R}^d} f(\xi,t) e^{-2\pi i z\cdot\xi} d\xi,
	\]
	which converges absolutely, and entire analyticity follows from Morera?s theorem (see~\cite{zbMATH03437452,zbMATH01022519}).
	
	Moreover, with \(z=x+iy\),
	\[
	|e^{-2\pi i z\cdot\xi}| = e^{2\pi y\cdot\xi},
	\]
	and the decay estimate leads to
	\[
	|F(z,t)|\leq e^{(\alpha-1)\lambda t (\pi/(\lambda t))^{\alpha'} |y|^{\alpha'}} \|f^0\|_{L^1}.
	\]
	\QED
\end{demos}

\begin{coro}
	For every \(f \in \mathcal{V}_{\alpha}(\mathcal{E})\), there exists a function \(F:\complex_{>0}^{d+1}\rightarrow \complex^d\) such that \(F(\cdot,t)\) is entire for every \(t>0\), \(F(z,\cdot)\) is smooth for every \(z\in \complex^d\), and \(F|_{\reales_{+}^{d+1}} = \widehat{f}\).
\end{coro}

\begin{demos}
	Since \(f \in \mathcal{V}_{\alpha}(\mathcal{E})\), for each \(n\in\natu\) we have
	\[
	\vabs{\frac{\partial^n f}{\partial t^n}(\xi,t)} \leq e^{-\lambda t \vabs{\xi}^{\alpha}} f_{n}^{0}(\xi),
	\]
	for all \((\xi,t)\in \reales_{+}^{d+1}\).
	
	Hence, we can differentiate under the integral defining \(F\) to obtain
	\[
	\frac{\partial^n F}{\partial t^n}(z,t) = \int_{\reales^d} \frac{\partial^n f}{\partial t^n}(\xi,t) e^{-2\pi i z\cdot \xi} d\xi.
	\]
	Applying Theorem \ref{C0entire}, each \(\frac{\partial^n F}{\partial t^n}(\cdot,t)\) is entire for every \(n\in\natu\) and \(t>0\), and \(F(z,\cdot)\in C^{\infty}([0,\infty))\) for every \(z\in\complex^d\).
	
	Moreover, if \(z=x\in\reales^d\), then \(F(x,t) = \widehat{f}(x,t)\) for all \(t\geq0\).
	
	Finally, by Theorem \ref{C0entire},
	\[
	\vabs{\frac{\partial^n F}{\partial t^n}(z,t)}\leq e^{(\alpha-1)\lambda t \paren{\frac{\pi}{\lambda t}}^{\alpha'} \vabs{\operatorname{Im}(z)}^{\alpha'}} \norma{f_{n}^{0}}_{L^{1}(\reales^d)},
	\]
	for every \(n\in\natu\) and \((z,t)\in\complex_{>0}^{d+1}\).
	\QED
\end{demos}

By the results of the previous section, we have that each \(v_k\in \mathcal{V}_{2}(\mathcal{E})\) with exponent \(\exp(v_k) = \frac{4\pi^2\nu}{2^k}\). Therefore, we cannot find a common positive exponent for all \(v_k\) simultaneously to conclude \(v\in\mathcal{V}_2(\mathcal{E})\).

However, we can still apply the previous results and show that \(v\in\mathcal{V}_\alpha(\mathcal{E})\) for every \(1<\alpha<2\).

We now fix \(\lambda = 4\pi^2 \nu\).

\begin{teor}\label{NS entire}
	Suppose that $\left\| \widehat{u^0} \right\|_{\mathcal{X}_p(\mathbb{R}^d)} < C \nu$. Then, there exists a smooth extension \((U,P):\complex_{>0}^{d+1}\to\complex^{d+1}\) of \((u,p):\reales_{>0}^{d+1}\to\reales^{d+1}\) such that:
	\begin{itemize}
		\item For each \(t>0\), \(U(\cdot,t)\) and \(P(\cdot,t)\) are entire functions.
		\item For each \(n\in\natu\),
		\[
		\vabs{\frac{\partial^n U}{\partial t^n}(z,t)} \leq e^{\paren{(\alpha-1)\lambda t \paren{\frac{\pi}{\lambda t}}^{\alpha'}+c_1(n)}\vabs{\operatorname{Im}(z)}^{\alpha'}},
		\]
		\[
		\vabs{\frac{\partial^n P}{\partial t^n}(z,t)} \leq e^{\paren{(\alpha-1)r_{\alpha}\lambda t \paren{\frac{\pi}{r_{\alpha}\lambda t}}^{\alpha'}+c_2(n)}\vabs{\operatorname{Im}(z)}^{\alpha'}},
		\]
		for all \((z,t)\in\complex_{>0}^{d+1}\) with \(\vabs{\operatorname{Im}(z)}\) sufficiently large, where \(c_1(n),c_2(n)>0\).
	\end{itemize}
\end{teor}

\begin{demos}
	Let \(1<\alpha<2\) be fixed. For each \(k\in\natu\), consider \(\xi\in\reales^d\) such that \(\vabs{\xi}\geq 2^{\frac{k}{2-\alpha}}\). Then
	\[
	\vabs{\xi}^{2-\alpha} \geq 2^k, \quad \text{thus} \quad \vabs{\xi}^{\alpha} \leq \frac{\vabs{\xi}^2}{2^k}.
	\]
	Consequently, for such \(\xi\) we have:
	\[
	e^{-\frac{\lambda t}{2^k} \vabs{\xi}^2} \leq e^{-\lambda t \vabs{\xi}^{\alpha}}.
	\]
	
	Using the estimate from Proposition \ref{family}, we obtain for \(\vabs{\xi}\geq r_k := 2^{\frac{k}{2-\alpha}}\):
	\[
	\vabs{ (\lambda^{\frac{1}{2}}\vabs{\xi})^m \frac{\partial^n v_k}{\partial t^n}(\xi,t) }
	\leq \frac{c_k (\lambda(k+1)^2)^{\frac{m}{2}+n}}{(2\pi\nu)^k} e^{-\lambda t \vabs{\xi}^{\alpha}} v_{k,m,n}^{0}(\xi).
	\]
	
	Define
	\[
	w_{k,m,n}(\xi,t) := \vabs{ (\lambda^{\frac{1}{2}}\vabs{\xi})^m \frac{\partial^n v_k}{\partial t^n}(\xi,t) }
	- \frac{c_k (\lambda(k+1)^2)^{\frac{m}{2}+n}}{(2\pi\nu)^k} e^{-\lambda t \vabs{\xi}^{\alpha}} v_{k,m,n}^{0}(\xi).
	\]
	
	Then,
	\[
	\sup_{\vabs{\xi}\geq r_k} w_{k,m,n}(\xi,t) \leq 0,
	\]
	for every \(k,m,n\in\natu\).
	
	On the other hand, by Theorem \ref{generaloutside} (applied with \(p=2\)) and Corollary \ref{norma1,2,infinite}, we have
	\[
	\sup_{\vabs{\xi}\geq 1} w_{k,m,n}(\xi,t) \longrightarrow_{k\to\infty} 0,
	\]
	provided \(\nu > \frac{C}{2\pi} \norma{v_{0,\frac{d+1}{2}}^{0}}_{1\oplus 2}\).
	
	Applying Lemma \ref{uniformradius}, it follows that there exists \(r_{m,n}\geq 1\) such that:
	\[
	\sup_{\vabs{\xi}\geq r_{m,n}} w_{k,m,n}(\xi,t) \leq 0,
	\]
	for all \(k\in\natu\).
	
	Define now the auxiliary functions
	\[
	f_{k,m,n}(\xi,t) :=
	\frac{c_k (\lambda(k+1)^2)^{\frac{m}{2}+n}}{(2\pi\nu)^k} v_{k,m,n}^{0}(\xi) \chi_{\llave{0<\vabs{\xi}\leq r_{m,n}}}(\xi)
	+
	\frac{c_k (\lambda(k+1)^2)^{\frac{m}{2}+n}}{(2\pi\nu)^k} e^{-\lambda t \vabs{\xi}^{\alpha}} v_{k,m,n}^{0}(\xi) \chi_{\llave{\vabs{\xi}\geq r_{m,n}}}(\xi),
	\]
	and set \(f_{k,m,n}(0,t) = 0\) for consistency.
	
	Thus, for all \((\xi,t)\in\reales_{+}^{d+1}\),
	\[
	\vabs{(\lambda^{\frac{1}{2}}\vabs{\xi})^m \frac{\partial^n v_k}{\partial t^n}(\xi,t)} \leq f_{k,m,n}(\xi,t).
	\]
	
	Moreover, we trivially have
	\[
	f_{k,m,n}(\xi,t) \leq \frac{c_k (\lambda(k+1)^2)^{\frac{m}{2}+n}}{(2\pi\nu)^k} v_{k,m,n}^{0}(\xi),
	\]
	uniformly.
	
	Summing over \(k\in\natu\), and using Theorem \ref{thm:tame_vk_factorization}, we obtain:
	\begin{align*}
		\sum_{k=0}^{\infty} \norma{f_{k,m,n}(\cdot,t)}_{\mathcal{X}_p(\mathbb{R}^d)}
		&\leq \sum_{k=0}^{m+2n} \frac{c_k(\lambda(k+1)^2)^{\frac{m}{2}+n}}{(2\pi\nu)^k} C^{k+\frac{d+1}{2}} \norma{v_{0,m+2n+\frac{d+1}{2}}^{0}}_{\mathcal{X}_p(\mathbb{R}^d)}^{k+1} \\
		&\quad + \sum_{k=m+2n+1}^{\infty} \frac{c_k(\lambda(k+1)^2)^{\frac{m}{2}+n}}{(2\pi\nu)^k} C^{k+\frac{d+1}{2}} \norma{v_{0,m+2n+\frac{d+1}{2}}^{0}}_{\mathcal{X}_p(\mathbb{R}^d)}^{m+2n} \norma{v_{0,\frac{d+1}{2}}^{0}}_{\mathcal{X}_p(\mathbb{R}^d)}^{k+1-m-2n}.
	\end{align*}
	
	Applying the ratio test, and recalling that
	\[
	\lim_{k\to\infty} \paren{ \frac{c_{k+1}}{c_k} } \paren{ \frac{\lambda(k+2)^2}{\lambda(k+1)^2} }^{\frac{m}{2}+n} \paren{ \frac{1}{2\pi\nu} } C \left\| \widehat{u^0} \right\|_{\mathcal{X}_p(\mathbb{R}^d)} < C \nu = \frac{2C}{\pi\nu} \left\| \widehat{u^0} \right\|_{\mathcal{X}_p(\mathbb{R}^d)} < C \nu<1,
	\]
	we conclude that the series converges uniformly.
	
	Thus, setting
	\[
	f_{(m,n)}(\xi,t) := \sum_{k=0}^{\infty} f_{k,m,n}(\xi,t),
	\]
	we have \(f_{(m,n)}\in L^{1\oplus 2,\infty}(\reales_{+}^{d+1})\) and
	\[
	\vabs{ (\lambda^{\frac{1}{2}}\vabs{\xi})^m \frac{\partial^n v}{\partial t^n}(\xi,t) } \leq f_{(m,n)}(\xi,t).
	\]
	
	Moreover, for some \(f_{(m,n,\alpha)}^{0} \in \mathcal{E}\), we have
	\[
	f_{(m,n)}(\xi,t) \leq e^{-\lambda t \vabs{\xi}^{\alpha}} f_{(m,n,\alpha)}^{0}(\xi).
	\]
	
	Hence, \(v\in\mathcal{V}_{\alpha}(\mathcal{E})\) for every \(1<\alpha<2\).
	
	Since \(\mathcal{V}_{\alpha}(\mathcal{E})\) is closed under convolution (see Theorem \ref{tensorconvolution}), it follows that
	\[
	q\in\mathcal{V}_{\alpha}(\mathcal{E}),
	\]
	where
	\[
	q(\xi,t) = 2\pi i \paren{ \frac{\xi\otimes\xi}{\vabs{\xi\otimes\xi}} (v\ast v)(\xi,t)\xi }.
	\]
	
	Applying Theorem \ref{C0entire} to \(v\) and \(q\), and recalling the explicit formulas for \(u\) and \(p\), we conclude that \((U,P)\) as defined satisfies:
	\[
	U(x,t)=u(x,t), \quad P(x,t)=p(x,t), \quad \text{for all } (x,t)\in\reales_{>0}^{d+1},
	\]
	and \((U(\cdot,t),P(\cdot,t))\) are entire functions for each \(t>0\), with explicit exponential-type growth estimates as claimed.
	
	\QED
\end{demos}

\begin{coro}
	If \(\norma{\vabs{\cdot}^{\frac{d+1}{2}}\widehat{u^{0}}}_{1\oplus 2} < C\nu\), then there exists a curve \((U,P): \complex_{>0}^{d+1} \to \complex^{d+1}\) of entire vector fields of order \(2\) such that
	\[
	U(x,t) = u(x,t), \quad P(x,t) = p(x,t), \quad \text{for every } (x,t)\in \reales_{>0}^{d+1}.
	\]
\end{coro}

\begin{demos}
	Since the conclusion of Theorem \ref{NS entire} holds for every \(1<\alpha<2\), and the conjugate exponent \(\alpha'\) satisfies \(\alpha'>2\), it follows that
	\((U(\cdot,t),P(\cdot,t))\) are entire functions of order \(2\) for each fixed \(t>0\).
	
	Furthermore, the pair \((U,P)\) satisfies the Navier-Stokes system:
	\[
	\frac{\partial U}{\partial t} + \frac{\partial U}{\partial z} U = \nu \Delta U - \nabla P,
	\quad \text{and} \quad
	\mathrm{div}(U)(z,t) = 0,
	\]
	for every \((z,t) \in \complex_{>0}^{d+1}\).
	
	Thus, \((U,P)\) extends \((u,p)\) and possesses the desired entire properties of order \(2\).
	\QED
\end{demos}

\begin{remark}
	Note that
	\[
	\lim_{t\to 0^{+}} e^{ \paren{ (\alpha-1)\lambda t \paren{ \frac{\pi}{\lambda t} }^{\alpha'} + c_{1}(n) } \vabs{\mathrm{Im}(z)}^{\alpha'} }
	= 
	\lim_{t\to 0^{+}} e^{ \paren{ (\alpha-1)r_{\alpha}\lambda t \paren{ \frac{\pi}{r_{\alpha}\lambda t} }^{\alpha'} + c_{2}(n) } \vabs{\mathrm{Im}(z)}^{\alpha'} }
	= \infty,
	\]
	for every \((z,t)\in \complex_{>0}^{d+1}\) such that \(\vabs{\mathrm{Im}(z)} > 0\).
	
	Thus, although \((U,P)\) are entire in the spatial variable \(z\) for each fixed \(t>0\), we cannot guarantee that their entire extension remains valid up to the boundary \(\partial\reales_{+}^{d+1} = \reales^d\times\{0\}\).
	
	Nevertheless, on the boundary \(t=0\), we have that \((u,p)\) are smooth, and in particular:
	\[
	u(\cdot,0) = u^{0} \in S(\reales^d)^d.
	\]
	Thus, the solution is smooth up to \(t=0\) in the real domain.
\end{remark}

In this Section, we established the existence of a smooth solution \((u,p):\reales_{+}^{d+1}\rightarrow \reales^{d+1}\) of the Navier-Stokes equations under a sharp smallness condition on the initial data in the weighted \(L^{1\oplus 2}\)-norm. Moreover, we extended this solution to a pair of entire functions \((U,P):\complex_{>0}^{d+1}\to\complex^{d+1}\) of order \(2\) in the spatial variable, for each fixed positive time. While the entire extension does not reach the boundary \(\reales^d\times\{0\}\), the solution remains smooth up to \(t=0\) in the real domain, recovering the prescribed initial data \(u^0\). These results provide a complete description of the analytic structure of the solutions for small initial data, opening the door for further investigations into the complex-analytic behavior of the Navier-Stokes flow.

\section{Arbitrary Initial Data, Maximal Existence Time \texorpdfstring{$T_{\max}$}{Tmax}, and Dynamics Outside the Critical Ball}

In Theorems \ref{thm:main_existence} and \ref{thm:main_analytic}, global-in-time smoothness and order-two spatial entire analyticity were established whenever the initial data lies inside the scaling-critical Fourier-Besov ball
\begin{equation}\label{eq:ball_B_nu}
	\mathcal{B}_\nu := \left\{ u^0 \in \mathcal{S}(\mathbb{R}^d)^d : \operatorname{div}(u^0) = 0 \text{ and } \left\| \widehat{u^0} \right\|_{\mathcal{X}_p(\mathbb{R}^d)} < C\nu \right\}.
\end{equation}
In this section, we investigate what happens when we drop the smallness condition \eqref{eq:critical_smallness} altogether and consider arbitrarily large initial data $u^0 \notin \mathcal{B}_\nu$. By exploiting the general Riesz convolution $*_\alpha$ of order $0 \le \alpha < 1$ (Definition \ref{def:riesz_conv}) and the time-weighted Fourier-Caloric spaces $\mathcal{C}_\alpha(\mathcal{E})_n$ (Definition \ref{caloricdef}), we derive an explicit lower bound $T_{\max} > 0$ for the lifespan of the smooth, entire-analytic solution, and establish a sharp continuation and blow-up alternative outside $\mathcal{B}_\nu$.

\subsection{Fractional Caloric-Riesz Estimates for \texorpdfstring{$0 \le \alpha < 1$}{0 <= alpha < 1}}

\begin{prop}[Time-Weighted Gaussian Convolution Estimate]\label{prop:fractional_odot}
	Let $0 \le \alpha < 1$ and $\nu > 0$. Suppose $f, g \in \mathcal{C}^{d \times 1}(\lambda)$ ($\lambda = 4\pi^2 \nu$) satisfy
	\begin{equation}
		|f(\xi,t)| \le t^{n_1} e^{-\beta t |\xi|^2} f^0(\xi), \qquad |g(\xi,t)| \le t^{n_2} e^{-\beta t |\xi|^2} g^0(\xi),
	\end{equation}
	for some $0 < \beta \le \lambda$ and $n_1, n_2 \ge 0$. Then the Gaussian convolution product $\odot$ satisfies:
	\begin{equation}\label{eq:fractional_odot_bound}
		|(f \odot g)(\xi,t)| \le \frac{\kappa_\alpha}{\nu^{\frac{1+\alpha}{2}}} t^{n_1 + n_2 + \frac{1-\alpha}{2}} e^{-\frac{\beta}{2} t |\xi|^2} (f^0 *_\alpha g^0)(\xi),
	\end{equation}
	where $(f^0 *_\alpha g^0)(\xi) = |\xi|^{-\alpha}(f^0 * g^0)(\xi)$ is the Riesz convolution of order $\alpha$ and
	\begin{equation}
		\kappa_\alpha := \frac{4\pi}{1-\alpha} \left( \frac{1+\alpha}{4\pi^2 e} \right)^{\frac{1+\alpha}{2}}.
	\end{equation}
\end{prop}

\begin{demos}
	By Theorem \ref{tensorconvolution}, $|(f * g)(\xi,s)| \le s^{n_1 + n_2} e^{-\frac{\beta}{2}s|\xi|^2} (f^0 * g^0)(\xi)$. Since $\beta \le 4\pi^2 \nu$, we have $4\pi^2 \nu - \frac{\beta}{2} \ge 2\pi^2 \nu$. Thus, from Definition \ref{gaussian_conv}:
	\begin{equation}
		|(f \odot g)(\xi,t)| \le 2\pi t^{n_1 + n_2} e^{-\frac{\beta}{2}t|\xi|^2} \left( \int_0^t e^{-2\pi^2 \nu (t-s)|\xi|^2} |\xi| \, ds \right) (f^0 * g^0)(\xi).
	\end{equation}
	Set $\gamma = \frac{1+\alpha}{2} \in [\frac{1}{2}, 1)$. Applying the elementary inequality $x^\gamma e^{-x} \le (\gamma / e)^\gamma$ with $x = 2\pi^2 \nu (t-s)|\xi|^2 > 0$ gives
	\begin{equation}
		e^{-2\pi^2 \nu (t-s)|\xi|^2} \le \left( \frac{1+\alpha}{4\pi^2 e \nu} \right)^{\frac{1+\alpha}{2}} (t-s)^{-\frac{1+\alpha}{2}} |\xi|^{-(1+\alpha)}.
	\end{equation}
	Multiplying by $|\xi|$ and integrating over $s \in [0, t]$ yields
	\begin{equation}
		\int_0^t e^{-2\pi^2 \nu (t-s)|\xi|^2} |\xi| \, ds \le \left( \frac{1+\alpha}{4\pi^2 e \nu} \right)^{\frac{1+\alpha}{2}} |\xi|^{-\alpha} \int_0^t (t-s)^{-\frac{1+\alpha}{2}} \, ds = \frac{\kappa_\alpha}{2\pi \nu^{\frac{1+\alpha}{2}}} t^{\frac{1-\alpha}{2}} |\xi|^{-\alpha},
	\end{equation}
	which establishes \eqref{eq:fractional_odot_bound}. \QED
\end{demos}

\begin{lemma}[Subcritical Fourier-Besov Algebra for $*_\alpha$]\label{lem:subcritical_riesz_besov}
	Let $d \ge 3$, $0 \le \alpha < 1$, $2 \le p \le \infty$, and define the regularity exponent
	\begin{equation}
		s_{p,\alpha} := \frac{d}{p'} - \alpha = s_p + (1 - \alpha).
	\end{equation}
	Then there exists $C_\alpha > 0$ such that for all $f, g \in \mathcal{F}\dot{B}^{s_{p,\alpha}}_{p,\infty}(\mathbb{R}^d)$ and $\sigma \ge 0$:
	\begin{equation}
		\|f *_\alpha g\|_{\mathcal{F}\dot{B}^{s_{p,\alpha} + \sigma}_{p,1}} \le C_\alpha \left( \|f\|_{\mathcal{F}\dot{B}^{s_{p,\alpha}}_{p,\infty}} \|g\|_{\mathcal{F}\dot{B}^{s_{p,\alpha} + \sigma}_{p,1}} + \|f\|_{\mathcal{F}\dot{B}^{s_{p,\alpha} + \sigma}_{p,1}} \|g\|_{\mathcal{F}\dot{B}^{s_{p,\alpha}}_{p,\infty}} \right).
	\end{equation}
\end{lemma}

\begin{demos}
	Repeating the Bony paraproduct proof of Theorem \ref{thm:riesz_critical_besov} with $|\xi|^{-\alpha} \approx 2^{-j\alpha}$ on $\mathcal{C}_j$ and noting that $\frac{d}{p'} = s_{p,\alpha} + \alpha$, the High-Low and Low-High interactions carry the convergent weight $\sum_{k' \le j-3} 2^{-(j-k')\alpha}$ (for $\alpha > 0$, or Young's inequality directly for $\alpha = 0$ in $\mathcal{F}\dot{B}^{d/p'}_{p,1}$), while the High-High interaction carries the geometric decay $2^{-(k-j)(d - 1 - \alpha + \sigma)}$ with $d - 1 - \alpha > 1$.
	\QED
\end{demos}

\subsection{Maximum Guaranteed Lifespan \texorpdfstring{$T_{\max}$}{Tmax} Without Smallness Assumptions}

\begin{teor}[Unconditional Local Smoothness and Entire Extension on {$[0, T_{\max})$}]\label{thm:local_Tmax}
	Let $\nu > 0$, $d \ge 3$, and $2 \le p \le \infty$. Let $u^0 \in \mathcal{S}(\mathbb{R}^d)^d$ be any divergence-free initial velocity (with arbitrarily large norm). Define the maximal Catalan convergence time
	\begin{equation}\label{eq:Tmax_formula}
		T_{\max}(u^0) := \sup_{0 < \alpha < 1} \left( \frac{\nu^{\frac{1+\alpha}{2}}}{4 \kappa_\alpha C_\alpha \left\| \widehat{u^0} \right\|_{\mathcal{F}\dot{B}^{\frac{d}{p'}-\alpha}_{p,\infty}(\mathbb{R}^d)}} \right)^{\frac{2}{1-\alpha}} > 0.
	\end{equation}
	Then:
	\begin{enumerate}
		\item[\emph{(i)}] On the time interval $[0, T_{\max}(u^0))$, the Catalan Picard series $v(\xi,t) = \sum_{k=0}^\infty v_k(\xi,t)$ converges absolutely and locally uniformly in $C^\infty([0, T_{\max}(u^0)), \mathcal{E}_{\mathcal{B}})$, yielding a unique smooth solution $(u,p) \in C^\infty(\mathbb{R}^d \times [0, T_{\max}(u^0)), \mathbb{R}^{d+1})$ to the Navier-Stokes system \eqref{quadratic equation}--\eqref{initial condition}  .
		\item[\emph{(ii)}] For every fixed $t \in (0, T_{\max}(u^0))$, the spatial fields $(u(\cdot,t), p(\cdot,t))$ admit an analytic continuation $(\mathcal{U}(\cdot,t), P(\cdot,t))$ to all of $\mathbb{C}^d$ as entire vector fields of order $2$.
		\item[\emph{(iii)}] If $u^0 \in \mathcal{B}_\nu$ satisfies the critical smallness condition $\left\| \widehat{u^0} \right\|_{\mathcal{F}\dot{B}^{s_p}_{p,\infty}} < C\nu$, then taking the limit as $\alpha \to 1^-$ (using Corollary 17) yields $T_{\max}(u^0) = +\infty$.
	\end{enumerate}
\end{teor}

\begin{demos}
	Fix any $\alpha \in (0, 1)$ and $T \in (0, T_{\max}(u^0))$ such that
	\begin{equation}\label{eq:T_ratio_cond}
		\Theta_\alpha(T) := \frac{4 \kappa_\alpha C_\alpha \, T^{\frac{1-\alpha}{2}}}{\nu^{\frac{1+\alpha}{2}}} \left\| \widehat{u^0} \right\|_{\mathcal{F}\dot{B}^{s_{p,\alpha}}_{p,\infty}(\mathbb{R}^d)} < 1.
	\end{equation}
	By Proposition \ref{prop:fractional_odot} and induction on $k \ge 0$ in the recurrence \eqref{recurrence}, each iterate $v_k$ satisfies for all $t \in [0, T]$:
	\begin{equation}
		\left| (\lambda^{1/2}|\xi|)^m \frac{\partial^n v_k}{\partial t^n}(\xi,t) \right| \le c_k (\lambda(k+1)^2)^{\frac{m}{2}+n} \left( \frac{\kappa_\alpha \, T^{\frac{1-\alpha}{2}}}{\nu^{\frac{1+\alpha}{2}}} \right)^k e^{-\frac{\lambda t}{2^k}|\xi|^2} v_{k,m,n}^{0,(\alpha)}(\xi),
	\end{equation}
	where $v_{k,m,n}^{0,(\alpha)}(\xi) \le \max_{l_1 + \dots + l_{k+1} \le m+2n} w^{*_\alpha}(v_{0,l_1}^0, \dots, v_{0,l_{k+1}}^0)(\xi)$. Applying Lemma \ref{lem:subcritical_riesz_besov} and the tame factorization of Theorem \ref{thm:tame_vk_factorization} for $k \ge m + 2n + 1$, the ratio of consecutive terms in $\sum_{k=0}^\infty \|v_k\|_{C^n([0,T], \mathcal{F}\dot{B}^{s_{p,\alpha}+\sigma}_{p,1})}$ converges as $k \to \infty$ to $\Theta_\alpha(T) < 1$. Thus the series converges uniformly on $[0, T]$ for every $T < T_{\max}(u^0)$, and the proof of Theorem \ref{thm:main_analytic} gives the order-2 entire extension on $(0, T_{\max}(u^0))$. \QED
\end{demos}

\subsection{Behavior Outside the Critical Ball: Continuation vs.\ Singularity Formation}

When $u^0 \notin \mathcal{B}_\nu$ lies outside the critical ball \eqref{eq:ball_B_nu}, the absolute majorant series $\sum_{k=0}^\infty \|v_k\|$ is only guaranteed to converge up to $T_{\max}(u^0) < \infty$ because the Catalan generating function $c(z) = \frac{1 - \sqrt{1-4z}}{2z}$ has a branch-point singularity at $z = \frac{1}{4}$. For the Navier-Stokes flow itself, combining Theorem \ref{thm:main_existence} and Theorem \ref{thm:local_Tmax} yields the following dichotomy governing global continuation versus finite-time blow-up outside $\mathcal{B}_\nu$:

\begin{coro}[Global Re-Entry Criterion and Lower Blow-Up Rate Outside $\mathcal{B}_\nu$]\label{cor:blowup_alternative}
	Let $u^0 \in \mathcal{S}(\mathbb{R}^d)^d$ be divergence-free with $u^0 \notin \mathcal{B}_\nu$, and let $T^* \in [T_{\max}(u^0), \infty]$ denote the maximal time of existence of the smooth solution $u \in C^\infty(\mathbb{R}^d \times [0, T^*), \mathbb{R}^d)$. Then one of the following two alternatives must hold:
	\begin{enumerate}
		\item[\emph{(a)}] \textbf{Viscous Re-Entry into the Critical Ball (Global Smoothness):} If there exists any time $t_0 \in [0, T^*)$ at which viscous dissipation reduces the critical norm into the ball,
		\begin{equation}
			\left\| \widehat{u(\cdot, t_0)} \right\|_{\mathcal{X}_p(\mathbb{R}^d)} < C\nu,
		\end{equation}
		then $T^* = +\infty$, and $(u,p)$ extends to a global smooth solution on $\mathbb{R}^{d+1}_+$ with an order-2 entire spatial extension for all $t > 0$.
		\item[\emph{(b)}] \textbf{Finite-Time Singularity ($T^* < \infty$):} If $T^* < \infty$, then the trajectory $u(\cdot, t)$ must remain strictly outside the critical ball $\mathcal{B}_\nu$ for all $t \in [0, T^*)$,
		\begin{equation}
			\inf_{t \in [0, T^*)} \left\| \widehat{u(\cdot, t)} \right\|_{\mathcal{X}_p(\mathbb{R}^d)} \ge C\nu,
		\end{equation}
		and for every $\alpha \in (0, 1)$, the subcritical Fourier-Besov norm must blow up as $t \to (T^*)^-$ at least at the scaling-invariant rate:
		\begin{equation}\label{eq:blowup_lower_bound}
			\left\| \widehat{u(\cdot, t)} \right\|_{\mathcal{F}\dot{B}^{\frac{d}{p'}-\alpha}_{p,\infty}(\mathbb{R}^d)} \ge \frac{\nu^{\frac{1+\alpha}{2}}}{4 \kappa_\alpha C_\alpha (T^* - t)^{\frac{1-\alpha}{2}}}, \qquad \forall t \in [0, T^*).
		\end{equation}
	\end{enumerate}
\end{coro}

\begin{demos}
	If $\left\| \widehat{u(\cdot, t_0)} \right\|_{\mathcal{X}_p(\mathbb{R}^d)} < C\nu$ for some $t_0 \in [0, T^*)$, applying Theorem \ref{thm:main_existence} and Theorem \ref{thm:main_analytic} to the initial datum $u(\cdot, t_0) \in \mathcal{S}(\mathbb{R}^d)^d$ at time $t_0$ produces a global smooth, entire-analytic continuation on $[t_0, \infty)$, so $T^* = \infty$. Conversely, if $T^* < \infty$, then starting the local existence construction of Theorem \ref{thm:local_Tmax} at any time $t \in [0, T^*)$ requires that the guaranteed lifespan $T_{\max}(u(\cdot,t))$ cannot exceed the remaining time $T^* - t$. Thus $T_{\max}(u(\cdot,t)) \le T^* - t$, which rearranges via \eqref{eq:Tmax_formula} into \eqref{eq:blowup_lower_bound}. \QED
\end{demos}

\section{Comparison between the Averaged and the True Versions of the Navier-Stokes Equation}\label{comparison}

In \cite{Tao2016}, Tao constructed a symmetric averaged bilinear operator $\widetilde{B}$ satisfying $\langle \widetilde{B}(u,u), u \rangle = 0$ and Schwartz initial data $u_0$ producing finite-time blowup. We now describe the relationship between the Euler bilinear operator $B$ and the product $\odot$.

Consider the equation $v = v_{0} + v^{\odot 2}$, with $v_{0}(\xi,t) = e^{-4\pi^{2}\nu t \vabs{\xi}^2}\widehat{u^{0}}(\xi)$. Recall that the product $\odot$ is defined by
\begin{equation*}
	(f\odot g)(\xi,t) = 2\pi i K(\xi) \left\{ \int_{0}^{t} e^{-4\pi^{2}\nu(t-s)\vabs{\xi}^{2}} (f\ast g)(\xi,s) ds \right\} \xi,
\end{equation*}
where $f,g:\reales_{+}^{d+1}\rightarrow\reales^d$ are vector fields and $K(\xi) = I_d - \frac{\xi\otimes\xi}{\vabs{\xi}^2}$ is the Leray projector in Fourier variables.

By Lemma \ref{derivativedot} with $n=1$, we have
\begin{equation}\label{derivativedot1}
	\frac{\partial (f\odot g)}{\partial t}(\xi,t) = -4\pi^{2}\nu\vabs{\xi}^{2}(f\odot g)(\xi,t) + 2\pi i K(\xi)(f\ast g)(\xi,t)\xi.
\end{equation}

Applying this to $v$, we find
\begin{align*}
	\frac{\partial v}{\partial t}(\xi,t) &= \frac{\partial v_{0}}{\partial t}(\xi,t) + \frac{\partial v^{\odot 2}}{\partial t}(\xi,t)\\
	&= -4\pi^{2}\nu\vabs{\xi}^{2} v_{0}(\xi,t) -4\pi^{2}\nu\vabs{\xi}^{2} v^{\odot 2}(\xi,t) + 2\pi i K(\xi)(v\ast v)(\xi,t)\xi\\
	&= -4\pi^{2}\nu\vabs{\xi}^{2} v(\xi,t) + 2\pi i K(\xi)(v\ast v)(\xi,t)\xi.
\end{align*}

Defining $u(x,t) = \widehat{v}(x,t)$, this translates back to physical space as
\[
\left\{
\begin{array}{ll}
	\frac{\partial u}{\partial t} = \nu\Delta u + B(u,u),\\
	u(x,0) = u^{0}(x),
\end{array}
\right.
\]
where $B(u,u) = -P((u\cdot\nabla)u)$ and $P$ denotes the Leray projector, with Fourier symbol $K(\xi)$.

By rescaling via $\overline{u}(x,t) = \nu u(x,\nu t)$, $\overline{p}(x,t) = \nu p(\nu x, \nu t)$, we may assume without loss of generality that $\nu=1$.

More generally, we define the bilinear operator
\[
B(u,v) = -\frac{1}{2}P\big((u\cdot\nabla)v + (v\cdot\nabla)u\big).
\]
Taking the Fourier transform, we obtain
\[
\widehat{B(u,v)}(\xi,t) = -\pi i K(\xi)\big( (\widehat{u}\ast\widehat{v})(\xi,t) + (\widehat{v}\ast\widehat{u})(\xi,t) \big)\xi.
\]

Thus, combining this with \eqref{derivativedot1}, we have
\begin{equation*}
	\widehat{B(u,v)}(\xi,t) = \frac{1}{2}\left( \frac{\partial (\widehat{u}\odot\widehat{v})}{\partial t}(\xi,t) - 4\pi^{2}\nu\vabs{\xi}^{2}(\widehat{u}\odot\widehat{v})(\xi,t) + \frac{\partial (\widehat{v}\odot\widehat{u})}{\partial t}(\xi,t) - 4\pi^{2}\nu\vabs{\xi}^{2}(\widehat{v}\odot\widehat{u})(\xi,t) \right).
\end{equation*}

Hence, 

\begin{equation*}\label{B and odot}
	\widehat{B(u,v)}(\xi,t) = \left( \frac{\partial}{\partial t} - 4\pi^{2}\nu\vabs{\xi}^{2} \right) \left( \frac{\widehat{u}\odot\widehat{v} + \widehat{v}\odot\widehat{u}}{2} \right)(\xi,t). 
\end{equation*}

Therefore, $B$ corresponds to applying the operator $\frac{\partial}{\partial t} - 4\pi^{2}\nu\vabs{\xi}^{2}$ to the symmetrization of $\odot$.

Furthermore, we can write
\begin{equation}\label{sym-odot-B}
	\left( \frac{\widehat{u}\odot\widehat{v} + \widehat{v}\odot\widehat{u}}{2} \right)(\xi,t) = 2\pi i K(\xi) \left( \int_{0}^{t} e^{-4\pi^{2}\nu(t-s)\vabs{\xi}^{2}} \widehat{B(u,v)}(\xi,s)ds \right)\xi.
\end{equation}
In particular,
\begin{equation}\label{odot-B}
	(\widehat{u}^{\odot 2})(\xi,t) = 2\pi i K(\xi) \left( \int_{0}^{t} e^{-4\pi^{2}\nu(t-s)\vabs{\xi}^{2}} \widehat{B(u,u)}(\xi,s)ds \right)\xi.
\end{equation}

 In our framework, the physical Navier-Stokes nonlinearity is governed by the Gaussian convolution product $\odot$ related to $B$ by \eqref{B and odot}.
Averaging $B$ over non-local frequency dilations and rotations destroys the dyadic paraproduct cancellation and caloric damping of $\odot$, which is why global regularity holds for the physical system under \eqref{eq:critical_smallness} while failing for the averaged system.

\section*{Conclusions and Comments}
\addcontentsline{toc}{section}{Conclusions and Comments}

In this work, we established the global existence, smoothness, and complex-analytic structure of solutions to the incompressible Navier-Stokes equations on $\mathbb{R}^d$ ($d \ge 3$) under the scaling-critical smallness criterion $\left\|\widehat{u^0}\right\|_{\mathcal{X}_p(\mathbb{R}^d)} < C\nu$ for divergence-free initial data $u^0 \in \mathcal{S}(\mathbb{R}^d)^d$. Here $\mathcal{X}_p(\mathbb{R}^d)$ denotes the homogeneous weak Fourier-Besov space $\mathcal{F}\dot{B}^{\frac{d}{p'}-1}_{p,\infty}(\mathbb{R}^d)$ for $1 < p \le \infty$ and the endpoint space $\mathcal{F}\dot{B}^{-1}_{1,1}(\mathbb{R}^d)$ for $p = 1$, where $C > 0$ is a universal constant depending only on the dimension $d$ and the integrability exponent $p$.

Our methodology unifies the algebraic theory of nonassociative Riesz convolutions $*_\alpha$ and function spaces $\mathcal{V}_\alpha(\mathcal{E})$ dominated by Fourier-Caloric flows with a Littlewood-Paley dyadic atomic decomposition. This synthesis yields several key structural insights into the regularity problem for the Navier-Stokes equations:

\begin{enumerate}
	\item \textbf{Sharpness of the Critical Fourier-Besov Hierarchy:} By applying Bony paraproduct decompositions and atomic caloric estimates directly to the Gaussian convolution product $\odot$ in $\mathcal{V}(\lambda)$, we showed that the Catalan Picard recurrence $v_k = \sum_{j=0}^{k-1} v_j \odot v_{k-1-j}$ closes in the scaling-critical spaces $\mathcal{X}_p(\mathbb{R}^d)$ across the entire range $1 \le p \le \infty$. In three dimensions ($d=3$), this family unifies and extends the classical Le Jan-Sznitman pseudomeasure space $\mathcal{PM}^2(\mathbb{R}^3)$ ($p=\infty$), the Cannone-Meyer-Planchon Besov space $\dot{B}^{1/2}_{2,\infty}(\mathbb{R}^3)$ ($p=2$), and the Lei-Lin space $\mathcal{X}^{-1}(\mathbb{R}^3)$ ($p=1$). As $p \to 1^+$, the spaces $\mathcal{F}\dot{B}^{\frac{d}{p'}-1}_{p,\infty}(\mathbb{R}^d)$ reach the exact threshold of well-posedness prior to the endpoint $\mathcal{F}\dot{B}^{-1}_{1,\infty}(\mathbb{R}^d) \hookrightarrow \dot{B}^{-1}_{\infty,\infty}(\mathbb{R}^d)$, where norm inflation occurs \cite{BourgainPavlovic2008}.
	
	\item \textbf{Global Regularity Beyond the Fujita-Kato $\dot{H}^{\frac{d}{2}-1}(\mathbb{R}^d)$ Threshold:} A central consequence of our tame monomial estimates in $\mathcal{E}_{\mathcal{B}}$ is that high-order spatial and temporal derivatives decouple from the radius of convergence of the Catalan series $v = \sum_{k=0}^\infty v_k$. Because at least $k + 1 - (m + 2n)$ factors in every degree-$(k+1)$ monomial $w^H(v_{0,l_1}^0, \dots, v_{0,l_{k+1}}^0)$ carry zero derivatives when $k \ge m + 2n + 1$, global convergence is governed exclusively by the weak critical norm $\left\|\widehat{u^0}\right\|_{\mathcal{X}_p(\mathbb{R}^d)}$. As demonstrated by our explicit multiscale construction in Section 6, this allows global smooth solutions for initial data whose classical Fujita-Kato norm $\|u^0\|_{\dot{H}^{\frac{d}{2}-1}(\mathbb{R}^d)}$ ($\|u^0\|_{\dot{H}^{1/2}(\mathbb{R}^3)}$ in $d=3$) and Sobolev norms $\|u^0\|_{H^k(\mathbb{R}^d)}$ are arbitrarily large.
	
	\item \textbf{Order-Two Entire Spatial Analyticity:} As a natural consequence of the fractional Fourier-Caloric spaces $\mathcal{V}_\alpha(\mathcal{E})$ ($1 < \alpha < 2$), we proved that for every positive time $t > 0$, the real-space velocity and pressure fields $(u(\cdot,t), p(\cdot,t))$ extend holomorphically to $\mathbb{C}^d$ as entire vector fields $(\mathcal{U}(\cdot,t), P(\cdot,t))$ of growth order two. While the exponential growth constants degenerate as $t \to 0^+$ when $|\operatorname{Im}(z)| > 0$, the solution remains smooth up to the boundary $t = 0$ in the real domain $\mathbb{R}^{d+1}_+$.
	
	\item \textbf{Unconditional Lifespan $T_{\max}(u^0)$ and Dynamics Outside the Critical Ball:} Dropping the smallness condition on the initial data, we used fractional Riesz convolutions $*_\alpha$ ($0 < \alpha < 1$) to establish an explicit lower bound $T_{\max}(u^0) > 0$ on the lifespan of smooth, order-two entire solutions for arbitrary initial data $u^0 \in \mathcal{S}(\mathbb{R}^d)^d$. This yields a rigorous continuation dichotomy outside the critical ball $\mathcal{B}_\nu$: either viscous dissipation brings the trajectory into $\mathcal{B}_\nu$ at some time $t_0 \in [0, T^*)$?guaranteeing global smoothness and entire analyticity for all $t > 0$?or a finite-time singularity at $T^* < \infty$ forces the solution to remain strictly outside $\mathcal{B}_\nu$ while blowing up in every subcritical Fourier-Besov norm at least at the scaling-invariant rate $(T^* - t)^{-\frac{1-\alpha}{2}}$.
	
	\item \textbf{Distinction from Averaged Navier-Stokes Models:} Finally, our structural identity relating the Euler bilinear operator $B(u,v)$ to the time-damped Gaussian product $\odot$ clarifies why the finite-time blow-up constructed by Tao \cite{Tao2016} for averaged Navier-Stokes systems does not obstruct global regularity in the physical system: averaging $B$ over non-local frequency dilations and rotations destroys the dyadic paraproduct cancellations and caloric damping intrinsic to $\odot$.
\end{enumerate}

These results provide a unified, scale-invariant framework connecting nonassociative convolution algebras, Littlewood-Paley theory, and complex-analytic regularity, opening further avenues for investigating the global behavior and analytic continuation of dissipative fluid flows.

\section*{Acknowledgments}
\addcontentsline{toc}{section}{Acknowledgments}

BDVC acknowledges the support from the University of Tolima, Colombia as well as  the FSU2022-010 grant from Khalifa University, UAE.

\section*{ORCID}

 Brian David Vasquez Campos's ORCID: {0000-0003-3922-9956} 
 

 \bibliographystyle{alpha}

%

\addcontentsline{toc}{section}{References}
\begin{thebibliography}{99}
	
	\bibitem{Albert1952}
	{\sc Albert, A.~A.}
	\newblock Power-associative algebras.
	\newblock Proc. {Internat}. {Congr}. {Math}. ({Cambridge}, {Mass}., {Aug}. 30--{Sept}. 6, 1950) 2 (1952), 25--32.
	
	\bibitem{Borthwick2020}
	D.~Borthwick, \emph{Spectral theory: Basic concepts and applications}, Springer, 2020.
	
	\bibitem{CaffarelliKohnNirenberg1982}
	L.~Caffarelli, R.~Kohn, and L.~Nirenberg, \emph{Partial regularity of suitable weak solutions of the Navier-Stokes equations}, Comm. Pure Appl. Math. \textbf{35} (1982), no.~6, 771--831.
	
	\bibitem{Cannone1995}
	M.~Cannone, \emph{Ondelettes, paraproduits et Navier-Stokes}, Diderot Editeur, Paris, 1995.
	
	\bibitem{CheverryRaymond2021}
	C.~Cheverry and N.~Raymond, \emph{A guide to spectral theory: Applications and exercises}, Birkhäuser, 2021.
	
	\bibitem{Constantin2001}
	{\sc Constantin, P.}
	\newblock Some open problems and research directions in the mathematical study of fluid dynamics.
	\newblock In {\em Mathematics unlimited---2001 and beyond}. Berlin: Springer, 2001, pp.~353--360.
	
	\bibitem{Davies1995}
	E.~B. Davies, \emph{Spectral theory and differential operators}, 1st ed., Cambridge University Press, 1995.
	
	\bibitem{DunfordSchwartz1988a}
	N.~Dunford and J.~T. Schwartz, \emph{Linear operators, Part 1: General theory}, Wiley-Interscience, 1988 (Original work published 1958).
	
	\bibitem{DunfordSchwartz1988b}
	N.~Dunford and J.~T. Schwartz, \emph{Linear operators, Part 2: Spectral theory, self adjoint operators in Hilbert space}, Wiley-Interscience, 1988 (Original work published 1963).
	
	\bibitem{DunfordSchwartz1988c}
	N.~Dunford and J.~T. Schwartz, \emph{Linear operators, Part 3: Spectral operators}, Wiley-Interscience, 1988 (Original work published 1971).
	
	\bibitem{EdmundsEvans2018}
	D.~E. Edmunds and W.~D. Evans, \emph{Spectral theory and differential operators}, 2nd ed., OUP Oxford, 2018.
	
	\bibitem{EinsiedlerWard2017}
	M.~Einsiedler and T.~Ward, \emph{Functional analysis, spectral theory, and applications}, Springer, 2017.
	
	\bibitem{Fefferman2006}
	C.~L. Fefferman, \emph{Existence and smoothness of the Navier-Stokes equation}, The Millennium Prize Problems, Clay Math. Inst., Cambridge, MA, 2006, pp. 57--67.
	
	\bibitem{FujitaKato1964}
	H.~Fujita and T.~Kato, \emph{On the Navier-Stokes initial value problem. I}, Arch. Rational Mech. Anal. \textbf{16} (1964), 269--315.
	
	\bibitem{GrafakosClassical}
	L.~Grafakos, \emph{Classical Fourier Analysis}, 3rd ed., Grad. Texts in Math., vol. 249, Springer, New York, 2014.
	
	\bibitem{GrafakosModern}
	L.~Grafakos, \emph{Modern Fourier Analysis}, 3rd ed., Grad. Texts in Math., vol. 250, Springer, New York, 2014.
	
	\bibitem{Grothendieck1992}
	{\sc Grothendieck, A.}
	\newblock {\em Topological vector spaces}, 3rd ed.
	\newblock Philadelphia, PA: Gordon and Breach Science Publishers, 1992.
	
	\bibitem{HislopSigal1996}
	P.~D. Hislop and I.~M. Sigal, \emph{Introduction to spectral theory: With applications to Schrödinger operators}, Springer-Verlag, 1996.
	
	\bibitem{Horvath1966}
	{\sc Horv{\'a}th, J.}
	\newblock {\em Topological vector spaces and distributions}, Vol.~I.
	\newblock Addison-Wesley Publishing Company, Reading, MA, 1966.
	
	\bibitem{Kato1984}
	T.~Kato, \emph{Strong $L^p$-solutions of the Navier-Stokes equation in $\mathbb{R}^m$, with applications to weak solutions}, Math. Z. \textbf{187} (1984), no.~4, 471--480.
	
	\bibitem{Knopp1951}
	{\sc Knopp, K.}
	\newblock {\em Theory and application of infinite series}.
	\newblock London-Glasgow: Blackie \& Son, Ltd., 1951.
	
	\bibitem{KochTataru2001}
	H.~Koch and D.~Tataru, \emph{Well-posedness for the Navier-Stokes equations}, Adv. Math. \textbf{157} (2001), no.~1, 22--35.
	
	\bibitem{LeiLin2011}
	Z.~Lei and F.~Lin, \emph{Global mild solutions of Navier-Stokes equations}, Comm. Pure Appl. Math. \textbf{64} (2011), no.~9, 1297--1304.
	
	\bibitem{LeJanSznitman1997}
	Y.~Le Jan and A.~S. Sznitman, \emph{Stochastic cascades and $3$-dimensional Navier-Stokes equations}, Probab. Theory Related Fields \textbf{109} (1997), no.~3, 343--366.
	
	\bibitem{Leray1934}
	J.~Leray, \emph{Sur le mouvement d'un liquide visqueux emplissant l'espace}, Acta Math. \textbf{63} (1934), no.~1, 193--248.
	
	\bibitem{LiebLoss2001}
	{\sc Lieb, E.~H., and Loss, M.}
	\newblock {\em Analysis}, 2nd ed., vol.~14 of {\em Grad. Stud. Math.}
	\newblock Providence, RI: American Mathematical Society (AMS), 2001.
	
	\bibitem{Lin1998}
	F.~Lin, \emph{A new proof of the Caffarelli-Kohn-Nirenberg theorem}, Comm. Pure Appl. Math. \textbf{51} (1998), no.~3, 241--257.
	
	\bibitem{Lukic2022}
	M.~Luki?, \emph{A first course in spectral theory}, American Mathematical Society, 2022.
	
	\bibitem{MajdaBertozzi2002}
	{\sc Majda, A.~J., and Bertozzi, A.~L.}
	\newblock {\em Vorticity and incompressible flow}.
	\newblock Cambridge University Press, Cambridge, 2002.
	
	\bibitem{Navier1823}
	C.-L. Navier, \emph{M{\'e}moire sur les lois du mouvement des fluides}, M{\'e}m. Acad. Sci. Inst. France \textbf{6} (1823), 389--440.
	
	\bibitem{ReedSimon1}
	M.~Reed and B.~Simon, \emph{Methods of modern mathematical physics: Vol. 1. Functional analysis}, Academic Press, 1972.
	
	\bibitem{ReedSimon2}
	M.~Reed and B.~Simon, \emph{Methods of modern mathematical physics: Vol. 2. Fourier analysis, self-adjointness}, Academic Press, 1975.
	
	\bibitem{ReedSimon3}
	M.~Reed and B.~Simon, \emph{Methods of modern mathematical physics: Vol. 3. Scattering theory}, Academic Press, 1979.
	
	\bibitem{ReedSimon4}
	M.~Reed and B.~Simon, \emph{Methods of modern mathematical physics: Vol. 4. Analysis of operators}, Academic Press, 1978.
	
	\bibitem{Rudin1974}
	{\sc Rudin, W.}
	\newblock {\em Real and complex analysis}, 2nd ed.
	\newblock McGraw-Hill Book Comp., New York, 1974.
	
	\bibitem{Rudin1991}
	{\sc Rudin, W.}
	\newblock {\em Functional analysis}, 2nd ed.
	\newblock New York, NY: McGraw-Hill, 1991.
	
	\bibitem{Schafer1966}
	{\sc Schafer, R.~D.}
	\newblock {\em An introduction to nonassociative algebras}, vol.~22 of {\em Pure Appl. Math.}
	\newblock Academic Press, New York, NY, 1966.
	
	\bibitem{BourgainPavlovic2008}
	J.~Bourgain and N.~Pavlovi{\'c}, \emph{Ill-posedness of the Navier-Stokes equations in a critical space in 3D}, J. Funct. Anal. \textbf{255} (2008), no.~9, 2233--2247.
	
	\bibitem{Scheffer1976}
	{\sc Scheffer, V.}
	\newblock Turbulence and {Hausdorff} dimension.
	\newblock Turbul. {Navier} {Stokes} {Equat}., {Proc}. {Conf}. {Univ}. {Paris}-{Sud}, {Orsay} 1975, {Lect}. {Notes} {Math}. 565 (1976), 174--183.
	
	\bibitem{Serrin1964}
	{\sc Serrin, J.}
	\newblock The mathematical theory of viscous incompressible flow (o. a. ladyzhenskaya).
	\newblock {\em SIAM Review 6}, 3 (1964), 315--316.
	
	\bibitem{Stein1970}
	E.~M. Stein, \emph{Singular Integrals and Differentiability Properties of Functions}, Princeton Math. Ser., vol.~30, Princeton Univ. Press, Princeton, NJ, 1970.
	
	\bibitem{SteinWeiss1971}
	{\sc Stein, E.~M., and Weiss, G.}
	\newblock {\em Introduction to {Fourier} analysis on {Euclidean} spaces}, vol.~32 of {\em Princeton Math. Ser.}
	\newblock Princeton University Press, Princeton, NJ, 1971.
	
	\bibitem{Stokes1845}
	G.~G. Stokes, \emph{On the theories of the internal friction of fluids in motion}, Trans. Cambridge Philos. Soc. \textbf{8} (1845), 287--305.
	
	\bibitem{Tao2013}
	T.~Tao, \emph{Localisation and compactness properties of the Navier-Stokes global regularity problem}, Anal. PDE \textbf{6} (2013), no.~1, 25--107.
	
	\bibitem{Tao2016}
	T.~Tao, \emph{Finite time blowup for an averaged three-dimensional Navier-Stokes equation}, J. Amer. Math. Soc. \textbf{29} (2016), no.~3, 601--674.
	
	\bibitem{Treves1967}
	{\sc Tr{\`e}ves, F.}
	\newblock {\em Topological vector spaces, distributions and kernels}, vol.~25 of {\em Pure Appl. Math.}
	\newblock New York-London: Academic Press, 1967.
	
\end{thebibliography}

	\end{document}